\documentclass[11pt]{article}

\usepackage[english]{babel}

\usepackage{amsmath,amsfonts,amssymb,amsthm,enumerate,graphicx}
\usepackage{tikz,authblk}
\usepackage{subfig}
\usepackage{url}
\usepackage{gastex}
\usepackage{gastex}
\newtheorem{theo}{{Theorem}}[section]
\newtheorem{lemma}{{Lemma}}[section]

\newtheorem{cor}{{Corollary}}[section]
\newtheorem{defn}{{Definition}}[section]

\begin{document}

\author[1] {Dipendu Maity}
\author[2] { Ashish Kumar Upadhyay}

\affil[1,2]{Department of Mathematics, Indian Institute of Technology Patna, Bihta,  801\,103, India. ${}^1$dipendumaity@gmail.com, ${}^2$upadhyay@iitp.ac.in.}

\title{On enumeration of a class of toroidal graphs}


\date{April 21, 2017}

\maketitle

\vspace{-10mm}


\begin{abstract} We present enumerations of a class of toroidal graphs which give rise to semi-equivelar maps. There are eleven different types of semi-equivelar maps on the torus. These are of the types $\{3^{6}\}$, $\{4^{4}\}$, $\{6^{3}\}$, $\{3^{3}, 4^{2}\}$, $\{3^{2}, 4, 3, 4\}$, $\{3, 6, 3, 6\}$, $\{3^{4}, 6\}$, $\{4, 8^{2}\}$, $\{3, 12^{2}\}$, $\{4, 6, 12\}$, $\{3, 4, 6, 4\}$. We know the classification of the maps of types $\{3^{6}\}$, $\{4^{4}\}$, $\{6^{3}\}$ on the torus. In this article, we attempt to classify maps of types $\{3^{3}, 4^{2}\}$, $\{3^{2}, 4, 3, 4\}$, $\{3, 6, 3, 6\}$, $\{3^{4}, 6\}$, $\{4, 8^{2}\}$, $\{3, 12^{2}\}$, $\{4, 6, 12\}$, $\{3, 4, 6, 4\}$ on the torus.
\end{abstract}

{\small

{\bf AMS classification\,: 52B70, 05C30, 05C38}

{\bf Keywords\,:} Toroidal Graphs, Semi-Equivelar Maps, Cycles}



\section{Introduction}

A {\em map} $M$ is an embedding of a graph $G$ on a surface $S$ such that the closure of components of $S \setminus G$, called the faces of $M$, are closed $2$-cells, that is, each homeomorphic to $2$-disk. A map $M$ is said to be a polyhedral map (see Brehm and Schulte \cite{brehm:1997}) if the intersection of any two distinct faces is either empty, a common vertex, or a common edge. A {\em a-cycle} $C_{a}$ is a finite connected $2$-regular graph with $a$ vertices, and the {\em face sequence} of a vertex $v$ in a map is a finite sequence $(a^{p},b^{q},\cdots,m^{r})$ of powers of positive integers $a, b, \cdots, m \geq 3$ and $p, q, \cdots, r \geq 1$ in {\em cyclic} order such that through the vertex $v$, $p$ number of $C_{a}$ ($C_{a}$ denote the a-cycle), $q$ number of $C_{b}$, $\cdots$, $r$ number of $C_{m}$ are incident. A map $K$ is said to be {\em semi-equivelar} if face sequence of each vertex is same, see \cite{upa:2014}. 
Two maps of fixed type, on the torus, are called {\em isomorphic} if there exists a {\em homeomorphism} of the torus which sends vertices to vertices, edges to edges, faces to faces and preserves incidents. That is, if we consider two polyhedral complexes $K_{1}$ and $K_{2}$ then an isomorphism to be a map $f : K_{1}\rightarrow K_{2}$ such that $f|_{V(K_{1})} : V(K_{1}) \rightarrow V(K_{2})$ is a bijection and $f(\sigma)$ is a cell in $K_{2}$ if and only if $\sigma$ is a cell in $K_{1}$. There are eleven types $\{3^{6}\}$, $\{4^{4}\}$, $\{6^{3}\}$, $\{3^{3}, 4^{2}\}$, $\{3^{2}, 4, 3, 4\}$, $\{3, 6, 3, 6\}$, $\{3^{4}, 6\}$, $\{4, 8^{2}\}$, $\{3, 12^{2}\}$, $\{4, 6, 12\}$, $\{3, 4, 6, 4\}$ of semi-equivelar maps on the torus. In this article, we are interested to classify some of them up to isomorphism. It completes the classification of semi-equivelar maps on the torus. In this context, Altshuler \cite{alt:1972} has shown construction and enumeration of maps of types $\{3^{6}\}$ and $\{6^{3}\}$ on the torus. Kurth \cite{kurth:1986} has given an enumeration of semi-equivelar maps of types $\{3^{6}\}$, $\{4^{4}\}$, $\{6^{3}\}$ on the torus. Negami \cite{negami:1983} has studied uniqueness and faithfulness of embedding of a class of toroidal graphs. Brehm and K\"{u}hnel \cite{brehm:2008} have presented a classification of semi-equivelar maps of types $\{3^{6}\}$, $\{4^{4}\}$, $\{6^{3}\}$ on the torus. Tiwari and Upadhyay \cite{tiwari:2015} have classified semi-equivelar maps of types $\{3^{3}, 4^{2}\}$, $\{3^{2}, 4, 3, 4\}$, $\{3, 6, 3, 6\}$, $\{3^{4}, 6\}$, $\{4, 8^{2}\}$, $\{3, 12^{2}\}$, $\{4, 6, 12\}$, $\{3, 4, 6, 4\}$ on torus with up to twenty vertices. In this article, we devise a way of enumerating all semi-equivelar maps of types $\{3^{3}, 4^{2}\}$, $\{3^{2}, 4, 3, 4\}$, $\{3, 6, 3, 6\}$, $\{3^{4}, 6\}$, $\{4, 8^{2}\}$, $\{3, 12^{2}\}$, $\{4, 6, 12\}$, $\{3, 4, 6, 4\}$ on the torus and explicitly determine the maps with a small number of vertices. Therefore, we have the following theorem.

\begin{theo}\label{thm-main} The semi-equivelar maps with $n$ vertices of types $\{3^{3}, 4^{2}\}$, $\{3^{2}, 4, 3, 4\}$, $\{3, 6, 3, 6\}$, $\{3, 12^{2}\}$, $\{3^{4}, 6\}$, $\{4, 6, 12\}$, $\{3, 4, 6, 4\}$, $\{4, 8^{2}\}$ can be classified up to isomorphism on the torus. In Table \ref{table1}, \ref{table2}, \ref{table3}, \ref{table4}, \ref{table5}, \ref{table6}, \ref{table7}, \ref{table8}, we have given non-isomorphic objects for few vertices.
\end{theo}

More precisely, let $X$ = $\{3^{3}, 4^{2}\}$, $\{3^{2}, 4, 3, 4\}$, $\{3, 6, 3, 6\}$, $\{3, 12^{2}\}$, $\{3^{4}, 6\}$, $\{4, 6, 12\}$, $\{3, 4, 6, 4\}$ or $\{4, 8^{2}\}$ be a semi-equivelar type on the torus. We present an algorithmic approach of calculating different maps for the type $X$ on different number of vertices in the subsequent sections.
\section{Definitions} 

We now define some operations on graphs. Let $G_{1} = (V_{1}, E_{1})$ and $G_{2} = (V_{2}, E_{2})$ be two subgraphs of the same graph $G := (V, E)$. Then $union$ $G_{1}\cup G_{2}$ is a graph $G_3 = (V_3, E_3)$ where $V_3 = V_{1} \cup V_{2}$ and $E_3 = E_{1} \cup E_{2}$. Similarly, the $intersection$ $G_{1}\cap G_{2}$ is a graph $G_4 = (V_4, E_4)$ where $V_4 = V_{1}\cap V_{2}$ and $E_4 = E_{1}\cap E_{2}$. For more on graph theory see \cite{bondy:2008}.

We denote a cycle $u_{1} \mbox{-} u_{2}\mbox{-} \cdots \mbox{-} u_{k}\mbox{-} u_{1}$  by $C(u_{1},u_{2},\cdots,u_{k})$ 
and a path $w_{1} \mbox{-} w_{2} \mbox{-} \cdots \mbox{-} w_{x}$ by $P(w_{1}, w_{2}, \cdots, w_{x})$. Let $Q_{1} = P(u_{1}, \cdots, u_{k})$ be a path. We call a path $Q_{2} =
 P(v_{1}, \cdots, v_{r})$ to be a path extended from $Q_{1}$ if $V(Q_{1})\subset V(Q_{2}), E(Q_{1})\subset E(Q_{2})$, i.e., $Q_1$ is a subpath of $Q_2$. We say that the $Q_{2}$ is an {\em extended path} of the path $Q_1$.
 


We say that a cycle is {\em contractible} if it bounds a $2$-disk. If the cycle does not bound any $2$-disk on the torus then we say that the cycle is {\em non-contractible}. 


\section{Examples}\label{example}

\noindent\textbf{Example 2.1 :} 
Let $M$ be a semi-equivelar map of type $\{3^3, 4^2\}$ with $n$ vertices on the torus. The map $M$ has a $T(r, s, k)$ representation (defined later in Section \ref{33421}) for some $r, s, k \in \mathbb{N} \cup \{0\}$. Let the number of vertices $n = 14$. We get by Lemma \ref{lem3342-all-rep}, $n = rs = 14$ where $2 \mid s$. Hence, $s = 2$, $r = 7$ and $k = 2, 3, 4$ by Lemma \ref{lem3342-all-rep}. So, $T(r, s, k) = T(7, 2, 2),$  $T(7, 2, 3)$ and $T(7, 2, 4)$, see Figure 1, 2, 3 respectively. In $T(7, 2, 2)$, $C_{1,1} = C(u_1, u_2, \cdots, u_7)$ is a cycle of type $A_1$ (see definition of type $A_1$ in Section \ref{33421}), $C_{1,2} = C(u_1, u_8, u_3, u_{10}, u_5, u_{12}, u_7, u_{14}, u_2, u_9, u_4, u_{11}, u_6, u_{13})$ and $C_{1,3}= C(u_1, u_8, u_4, u_{11}, u_7, u_{14},$ $u_3, u_{10}, u_6, u_{13}, u_2, u_{9}, u_5, u_{12})$ are two cycles of type $A_2$ (see definition of type $A_2$ in Section \ref{33421}) and $C_{1,4} = C(u_3, u_{10}, u_5, u_{4})$ is a cycle of type $A_4$ (see definition of type $A_4$ in Section \ref{33421}). 
In $T(7, 2, 4)$, $C_{2,1} = C(v_1, v_2, \cdots, v_7)$ is of type $A_1$, $C_{2, 2 } = C(v_1, v_8, v_5, v_{12}, v_2, v_{9}, v_6, v_{13}, v_3,$ $ v_{10}, v_7, v_{14}, v_4, v_{11})$ and $C_{2, 3} = C(v_1, v_8, v_6, v_{13}, v_4, v_{11}, v_2,$ $v_{9}, v_7, v_{14}, v_5, v_{12}, v_3, v_{10})$ are of type $A_2$, and $C_{2,4} = C(v_5, v_{12}, v_3, v_{4})$ is of type $A_4$. 
In $T(7, 2, 3)$, $C_{3,1} = C(w_1, w_2, \cdots, w_7)$ is of type $A_1$, $C_{3, 2} = C(w_1, w_8, w_4, w_{11}, w_7, w_{14}, w_3, w_{10},$ $w_6, w_{13}, w_2, w_{9}, w_5, w_{12})$ and $C_{3,3} = C(w_1, w_8, w_5,$ $w_{12}, w_2, w_{9}, w_6, w_{13}, w_3, w_{10}, w_7, w_{14}, w_4,$ $w_{11})$ are of type $A_2$ and $C_{3,4} = C(w_4,$ $w_{5}, w_6, w_{7}, w_{11})$ is of type $A_4$.

In Section \ref{33421}, by Lemma \ref{lem3342-hom-length}, the cycles of type $A_1$ have same length and the cycles of type $A_2$ have at most two different lengths in $M$. So, $O_3 \not\cong O_1$ since length($C_{3, 4}$) $\neq$ length($C_{1, 4}$) and $O_3 \not\cong O_2$ since length($C_{3, 4}$) $\neq$ length($C_{2, 4}$). Thus, $O_3 \not\cong O_i$ for $i=1, 2$.  Now, length($C_{1, 1}$) = length($C_{2, 1}$), $\{$length($C_{1, 2}$), length($C_{1, 3})\}$ = $\{$length($C_{2, 2}$), length($C_{2, 3})\}$ and length($C_{1, 4}) = $length($C_{2, 4})$. We cut $T(7, 2, 4)$ along the path $P(v_5, v_{12}, v_3)$ and identify along the path $P(v_1, v_{8}, v_5)$. This gives a presentation of $T(7, 2, 4)$ in Figure 4. The Figure 4 has $T(7, 2, 2)$ representation. So, $O_1 \cong O_2$, (see the proof of Lemma \ref{lem3342-iso} for isomorphism map between $O_1$ and $O_2$). Thus, $O_1 \cong O_2$ and $O_3 \not\cong O_1, O_2$. Therefore, there are two semi-equivelar maps of type $\{3^{3}, 4^{2}\}$ with $14$ vertices on the torus up to isomorphism.

\vspace{5.9cm}

\begin{picture}(0,0)(5,-45)
\tiny
\setlength{\unitlength}{1.3mm}
\drawpolygon(0,0)(35,0)(35,10)(0,10)
\drawline[AHnb=0](0,5)(35,5)
\drawline[AHnb=0](5,0)(5,10)
\drawline[AHnb=0](10,0)(10,10)
\drawline[AHnb=0](15,0)(15,10)
\drawline[AHnb=0](20,0)(20,10)
\drawline[AHnb=0](25,0)(25,10)
\drawline[AHnb=0](30,0)(30,10)

\drawline[AHnb=0](0,5)(5,10)
\drawline[AHnb=0](5,5)(10,10)
\drawline[AHnb=0](10,5)(15,10)
\drawline[AHnb=0](15,5)(20,10)
\drawline[AHnb=0](20,5)(25,10)
\drawline[AHnb=0](25,5)(30,10)
\drawline[AHnb=0](30,5)(35,10)

\put(0,-1.1){\tiny ${u_{1}}$}
\put(5,-1.1){\tiny ${u_{2}}$}
\put(10,-1.1){\tiny ${u_{3}}$}
\put(15,-1.1){\tiny ${u_{4}}$}
\put(20,-1.1){\tiny ${u_{5}}$}
\put(25,-1.1){\tiny ${u_{6}}$}
\put(30,-1.1){\tiny ${u_{7}}$}
\put(35,-1.1){\tiny ${u_{1}}$}

\put(.5,3.8){\tiny${u_{8}}$}
\put(5.5,3.8){\tiny ${u_{9}}$}
\put(10.5,3.8){\tiny ${u_{10}}$}
\put(15.5,3.8){\tiny ${u_{11}}$}
\put(20.5,3.8){\tiny ${u_{12}}$}
\put(25.5,3.8){\tiny ${u_{13}}$}
\put(30.5,3.8){\tiny ${u_{14}}$}
\put(35.5,3.8){\tiny ${u_{8}}$}

\put(0.5,8.8){\tiny ${u_{3}}$}
\put(5.5,8.8){\tiny ${u_{4}}$}
\put(10.5,8.8){\tiny ${u_{5}}$}
\put(15.5,8.8){\tiny ${u_{6}}$}
\put(20.5,8.8){\tiny ${u_{7}}$}
\put(25.5,8.8){\tiny ${u_{1}}$}
\put(30.5,8.8){\tiny ${u_{2}}$}
\put(35.5,8.8){\tiny ${u_{3}}$}

\put(8,-4){\tiny  Figure 1 : {\boldmath T(7, 2, 2) : $O_1$}}
\end{picture}


\begin{picture}(0,0)(-46,-50)
\tiny
\setlength{\unitlength}{1.3mm}
\drawpolygon(0,0)(35,0)(35,10)(0,10)
\drawline[AHnb=0](0,5)(35,5)
\drawline[AHnb=0](5,0)(5,10)
\drawline[AHnb=0](10,0)(10,10)
\drawline[AHnb=0](15,0)(15,10)
\drawline[AHnb=0](20,0)(20,10)
\drawline[AHnb=0](25,0)(25,10)
\drawline[AHnb=0](30,0)(30,10)

\drawline[AHnb=0](0,5)(5,10)
\drawline[AHnb=0](5,5)(10,10)
\drawline[AHnb=0](10,5)(15,10)
\drawline[AHnb=0](15,5)(20,10)
\drawline[AHnb=0](20,5)(25,10)
\drawline[AHnb=0](25,5)(30,10)
\drawline[AHnb=0](30,5)(35,10)

\put(0,-1.1){\tiny ${v_{1}}$}
\put(5,-1.1){\tiny ${v_{2}}$}
\put(10,-1.1){\tiny ${v_{3}}$}
\put(15,-1.1){\tiny ${v_{4}}$}
\put(20,-1.1){\tiny ${v_{5}}$}
\put(25,-1.1){\tiny ${v_{6}}$}
\put(30,-1.1){\tiny ${v_{7}}$}
\put(35,-1.1){\tiny ${v_{1}}$}

\put(0.5,3.8){\tiny ${v_{8}}$}
\put(5.5,3.8){\tiny ${v_{9}}$}
\put(10.5,3.8){\tiny ${v_{10}}$}
\put(15.5,3.8){\tiny ${v_{11}}$}
\put(20.5,3.8){\tiny ${v_{12}}$}
\put(25.5,3.8){\tiny ${v_{13}}$}
\put(30.5,3.8){\tiny ${v_{14}}$}
\put(35.5,3.8){\tiny ${v_{8}}$}

\put(0.5,8.8){\tiny ${v_{5}}$}
\put(5.5,8.8){\tiny ${v_{6}}$}
\put(10.5,8.8){\tiny ${v_{7}}$}
\put(15.5,8.8){\tiny ${v_{1}}$}
\put(20.5,8.8){\tiny ${v_{2}}$}
\put(25.5,8.8){\tiny ${v_{3}}$}
\put(30.5,8.8){\tiny ${v_{4}}$}
\put(35.5,8.8){\tiny ${v_{5}}$}

\put(8,-4){\tiny  Figure 2 : {\boldmath T(7, 2, 4) : $O_2$}}

\end{picture}

\begin{picture}(0,0)(-97,-55)
\tiny
\setlength{\unitlength}{1.3mm}
\drawpolygon(0,0)(35,0)(35,10)(0,10)
\drawline[AHnb=0](0,5)(35,5)
\drawline[AHnb=0](5,0)(5,10)
\drawline[AHnb=0](10,0)(10,10)
\drawline[AHnb=0](15,0)(15,10)
\drawline[AHnb=0](20,0)(20,10)
\drawline[AHnb=0](25,0)(25,10)
\drawline[AHnb=0](30,0)(30,10)

\drawline[AHnb=0](0,5)(5,10)
\drawline[AHnb=0](5,5)(10,10)
\drawline[AHnb=0](10,5)(15,10)
\drawline[AHnb=0](15,5)(20,10)
\drawline[AHnb=0](20,5)(25,10)
\drawline[AHnb=0](25,5)(30,10)
\drawline[AHnb=0](30,5)(35,10)

\put(0,-1.1){\tiny ${w_{1}}$}
\put(5,-1.1){\tiny ${w_{2}}$}
\put(10,-1.1){\tiny ${w_{3}}$}
\put(15,-1.1){\tiny ${w_{4}}$}
\put(20,-1.1){\tiny ${w_{5}}$}
\put(25,-1.1){\tiny ${w_{6}}$}
\put(30,-1.1){\tiny ${w_{7}}$}
\put(35,-1.1){\tiny ${w_{1}}$}

\put(0.5,3.8){\tiny ${w_{8}}$}
\put(5.5,3.8){\tiny ${w_{9}}$}
\put(10.5,3.8){\tiny ${w_{10}}$}
\put(15.5,3.8){\tiny ${w_{11}}$}
\put(20.5,3.8){\tiny ${w_{12}}$}
\put(25.5,3.8){\tiny ${w_{13}}$}
\put(30.5,3.8){\tiny ${w_{14}}$}
\put(35.5,3.8){\tiny ${w_{8}}$}

\put(0.5,8.8){\tiny ${w_{4}}$}
\put(5.5,8.8){\tiny ${w_{5}}$}
\put(10.5,8.8){\tiny ${w_{6}}$}
\put(15.5,8.8){\tiny ${w_{7}}$}
\put(20.5,8.8){\tiny ${w_{1}}$}
\put(25.5,8.8){\tiny ${w_{2}}$}
\put(30.5,8.8){\tiny ${w_{3}}$}
\put(35.5,8.8){\tiny ${w_{4}}$}

\put(8,-4){\tiny  Figure 3 : {\boldmath T(7, 2, 3) : $O_3$}}

\end{picture}

\vspace{.8cm}

\begin{picture}(0,0)(0,-45)
\tiny
\setlength{\unitlength}{1.3mm}
\drawpolygon(0,0)(35,0)(35,10)(0,10)
\drawline[AHnb=0](0,5)(35,5)
\drawline[AHnb=0](5,0)(5,10)
\drawline[AHnb=0](10,0)(10,10)
\drawline[AHnb=0](15,0)(15,10)
\drawline[AHnb=0](20,0)(20,10)
\drawline[AHnb=0](25,0)(25,10)
\drawline[AHnb=0](30,0)(30,10)

\drawline[AHnb=0](0,5)(5,10)
\drawline[AHnb=0](5,5)(10,10)
\drawline[AHnb=0](10,5)(15,10)
\drawline[AHnb=0](15,5)(20,10)
\drawline[AHnb=0](20,5)(25,10)
\drawline[AHnb=0](25,5)(30,10)
\drawline[AHnb=0](30,5)(35,10)

\put(0,-1.1){\tiny ${v_{5}}$}
\put(5,-1.1){\tiny ${v_{4}}$}
\put(10,-1.1){\tiny ${v_{3}}$}
\put(15,-1.1){\tiny ${v_{2}}$}
\put(20,-1.1){\tiny ${v_{1}}$}
\put(25,-1.1){\tiny ${v_{7}}$}
\put(30,-1.1){\tiny ${v_{6}}$}
\put(35,-1.1){\tiny ${v_{5}}$}

\put(0.5,3.8){\tiny ${v_{12}}$}
\put(5.5,3.8){\tiny ${v_{11}}$}
\put(10.5,3.8){\tiny ${v_{10}}$}
\put(15.5,3.8){\tiny ${v_{9}}$}
\put(20.5,3.8){\tiny ${v_{8}}$}
\put(25.5,3.8){\tiny ${v_{14}}$}
\put(30.5,3.8){\tiny ${v_{13}}$}
\put(35.5,3.8){\tiny ${v_{12}}$}

\put(0.5,8.8){\tiny ${v_{3}}$}
\put(5.5,8.8){\tiny ${v_{2}}$}
\put(10.5,8.8){\tiny ${v_{1}}$}
\put(15.5,8.8){\tiny ${v_{7}}$}
\put(20.5,8.8){\tiny ${v_{6}}$}
\put(25.5,8.8){\tiny ${v_{5}}$}
\put(30.5,8.8){\tiny ${v_{4}}$}
\put(35.5,8.8){\tiny ${v_{3}}$}

\put(7,-4){\tiny  Figure 4 : {\boldmath T(7, 2, 2) : $O_2$}}

\end{picture}

\begin{picture}(0,0)(2,-27)
\tiny
\setlength{\unitlength}{1.3mm}
\drawpolygon(0,0)(40,0)(40,10)(0,10)
\drawline[AHnb=0](0,5)(40,5)
\drawline[AHnb=0](5,0)(5,10)
\drawline[AHnb=0](10,0)(10,10)
\drawline[AHnb=0](15,0)(15,10)
\drawline[AHnb=0](20,0)(20,10)
\drawline[AHnb=0](25,0)(25,10)
\drawline[AHnb=0](30,0)(30,10)
\drawline[AHnb=0](35,0)(35,10)

\drawline[AHnb=0](0,5)(5,10)
\drawline[AHnb=0](5,5)(10,0)
\drawline[AHnb=0](10,5)(15,10)
\drawline[AHnb=0](15,5)(20,0)
\drawline[AHnb=0](20,5)(25,10)
\drawline[AHnb=0](25,5)(30,0)
\drawline[AHnb=0](30,5)(35,10)
\drawline[AHnb=0](35,5)(40,0)

\put(0,-1.1){\tiny ${w_{2}}$}
\put(5,-1.1){\tiny ${w_{6}}$}
\put(10,-1.1){\tiny ${w_{10}}$}
\put(15,-1.1){\tiny ${w_{14}}$}
\put(20,-1.1){\tiny ${w_{4}}$}
\put(25,-1.1){\tiny ${w_{8}}$}
\put(30,-1.1){\tiny ${w_{12}}$}
\put(35,-1.1){\tiny ${w_{16}}$}
\put(40,-1.1){\tiny ${w_{2}}$}

\put(0.5,4){\tiny ${w_{1}}$}
\put(5.5,5.4){\tiny ${w_{5}}$}
\put(10.5,4){\tiny ${w_{9}}$}
\put(15.5,5.4){\tiny ${w_{13}}$}
\put(20.5,4){\tiny ${w_{3}}$}
\put(25.5,5.4){\tiny ${w_{7}}$}
\put(30.5,4){\tiny ${w_{11}}$}
\put(35.5,5.4){\tiny ${w_{15}}$}
\put(40.5,4){\tiny ${w_{1}}$}

\put(0.5,8.8){\tiny ${w_{4}}$}
\put(5.5,8.8){\tiny ${w_{8}}$}
\put(10.5,8.8){\tiny ${w_{12}}$}
\put(15.5,8.8){\tiny ${w_{16}}$}
\put(20.5,8.8){\tiny ${w_{2}}$}
\put(25.5,8.8){\tiny ${w_{6}}$}
\put(30.5,8.8){\tiny ${w_{10}}$}
\put(35.5,8.8){\tiny ${w_{14}}$}
\put(40.5,8.8){\tiny ${w_{4}}$}

\put(10,-4){\tiny  Figure 8 : {\boldmath T(8, 2, 4) : $O_4$}}
\end{picture}

\begin{picture}(0,0)(2,-8)
\tiny
\setlength{\unitlength}{1.3mm}
\drawpolygon(0,0)(40,0)(40,10)(0,10)
\drawline[AHnb=0](0,5)(40,5)
\drawline[AHnb=0](5,0)(5,10)
\drawline[AHnb=0](10,0)(10,10)
\drawline[AHnb=0](15,0)(15,10)
\drawline[AHnb=0](20,0)(20,10)
\drawline[AHnb=0](25,0)(25,10)
\drawline[AHnb=0](30,0)(30,10)
\drawline[AHnb=0](35,0)(35,10)

\drawline[AHnb=0](0,5)(5,10)
\drawline[AHnb=0](5,5)(10,0)
\drawline[AHnb=0](10,5)(15,10)
\drawline[AHnb=0](15,5)(20,0)
\drawline[AHnb=0](20,5)(25,10)
\drawline[AHnb=0](25,5)(30,0)
\drawline[AHnb=0](30,5)(35,10)
\drawline[AHnb=0](35,5)(40,0)

\put(0,-1.1){\tiny ${u_{1}}$}
\put(5,-1.1){\tiny ${u_{2}}$}
\put(10,-1.1){\tiny ${u_{3}}$}
\put(15,-1.1){\tiny ${u_{4}}$}
\put(20,-1.1){\tiny ${u_{5}}$}
\put(25,-1.1){\tiny ${u_{6}}$}
\put(30,-1.1){\tiny ${u_{7}}$}
\put(35,-1.1){\tiny ${u_{8}}$}
\put(40,-1.1){\tiny ${u_{1}}$}

\put(0.5,4){\tiny ${u_{9}}$}
\put(5.5,5.49){\tiny ${u_{10}}$}
\put(10.5,4){\tiny ${u_{11}}$}
\put(15.5,5.49){\tiny ${u_{12}}$}
\put(20.5,4){\tiny ${u_{13}}$}
\put(25.5,5.49){\tiny ${u_{14}}$}
\put(30.5,4){\tiny ${u_{15}}$}
\put(35.5,5.49){\tiny ${u_{16}}$}
\put(40.5,4){\tiny ${u_{9}}$}

\put(0.5,8.8){\tiny ${u_{5}}$}
\put(5.5,8.8){\tiny ${u_{6}}$}
\put(10.5,8.8){\tiny ${u_{7}}$}
\put(15.5,8.8){\tiny ${u_{8}}$}
\put(20.5,8.8){\tiny ${u_{1}}$}
\put(25.5,8.8){\tiny ${u_{2}}$}
\put(30.5,8.8){\tiny ${u_{3}}$}
\put(35.5,8.8){\tiny ${u_{4}}$}
\put(40.5,8.8){\tiny ${u_{5}}$}

\put(10,-4){Figure 9 : {\boldmath T(8, 2, 4) : $O_5$}}
\end{picture}

\begin{picture}(0,0)(-65,-47)
\tiny
\setlength{\unitlength}{1.3mm}
\drawpolygon(0,0)(20,0)(20,20)(0,20)
\drawline[AHnb=0](0,5)(20,5)
\drawline[AHnb=0](0,10)(20,10)
\drawline[AHnb=0](0,15)(20,15)

\drawline[AHnb=0](5,0)(5,20)
\drawline[AHnb=0](10,0)(10,20)
\drawline[AHnb=0](15,0)(15,20)

\drawline[AHnb=0](5,5)(10,0)
\drawline[AHnb=0](15,5)(20,0)

\drawline[AHnb=0](5,15)(10,10)
\drawline[AHnb=0](15,15)(20,10)

\drawline[AHnb=0](0,5)(5,10)
\drawline[AHnb=0](10,5)(15,10)

\drawline[AHnb=0](0,15)(5,20)
\drawline[AHnb=0](10,15)(15,20)

\put(0,-1.1){\tiny ${v_{1}}$}
\put(5,-1.1){\tiny ${v_{2}}$}
\put(10,-1.1){\tiny ${v_{3}}$}
\put(15,-1.1){\tiny ${v_{4}}$}
\put(20,-1.1){\tiny ${v_{1}}$}

\put(0.5,4){\tiny ${v_{5}}$}
\put(5.5,5.5){\tiny ${v_{6}}$}
\put(10.5,4){\tiny ${v_{7}}$}
\put(15.5,5.5){\tiny ${v_{8}}$}
\put(20.5,4){\tiny ${v_{5}}$}

\put(0.5,8.8){\tiny ${v_{9}}$}
\put(5.5,8.8){\tiny ${v_{10}}$}
\put(10.5,8.8){\tiny ${v_{11}}$}
\put(15.5,8.8){\tiny ${v_{12}}$}
\put(20.5,8.8){\tiny ${v_{9}}$}

\put(0.5,14){\tiny ${v_{13}}$}
\put(5.5,15.5){\tiny ${v_{14}}$}
\put(10.5,14){\tiny ${v_{15}}$}
\put(15.5,15.5){\tiny ${v_{16}}$}
\put(20.5,14){\tiny ${v_{13}}$}

\put(0.5,18.8){\tiny ${v_{1}}$}
\put(5.5,18.8){\tiny ${v_{2}}$}
\put(10.5,18.8){\tiny ${v_{3}}$}
\put(15.5,18.8){\tiny ${v_{4}}$}
\put(20.5,18.8){\tiny ${u_{1}}$}

\put(-1,-4){\tiny  Figure 6 : {\boldmath T(4, 4, 0) : $O_6$}}

\drawpolygon(25,0)(45,0)(45,20)(25,20)
\drawline[AHnb=0](25,5)(45,5)
\drawline[AHnb=0](25,10)(45,10)
\drawline[AHnb=0](25,15)(45,15)

\drawline[AHnb=0](30,0)(30,20)
\drawline[AHnb=0](35,0)(35,20)
\drawline[AHnb=0](40,0)(40,20)

\drawline[AHnb=0](30,5)(35,0)
\drawline[AHnb=0](40,5)(45,0)

\drawline[AHnb=0](30,15)(35,10)
\drawline[AHnb=0](40,15)(45,10)

\drawline[AHnb=0](25,5)(30,10)
\drawline[AHnb=0](35,5)(40,10)

\drawline[AHnb=0](25,15)(30,20)
\drawline[AHnb=0](35,15)(40,20)

\put(25,-1.1){\tiny ${w_{1}}$}
\put(30,-1.1){\tiny ${w_{2}}$}
\put(35,-1.1){\tiny ${w_{3}}$}
\put(40,-1.1){\tiny ${w_{4}}$}
\put(45,-1.1){\tiny ${w_{1}}$}

\put(25.5,4){\tiny ${w_{5}}$}
\put(30.5,5.5){\tiny ${w_{6}}$}
\put(35.5,4){\tiny ${w_{7}}$}
\put(40.5,5.5){\tiny ${w_{8}}$}
\put(45.5,4){\tiny ${w_{5}}$}

\put(25.5,8.8){\tiny ${w_{9}}$}
\put(30.5,8.8){\tiny ${w_{10}}$}
\put(35.5,8.8){\tiny ${w_{11}}$}
\put(40.5,8.8){\tiny ${w_{12}}$}
\put(45.5,8.8){\tiny ${w_{9}}$}

\put(25.5,14){\tiny ${w_{13}}$}
\put(30.5,15.5){\tiny ${w_{14}}$}
\put(35.5,14){\tiny ${w_{15}}$}
\put(40.5,15.5){\tiny ${w_{16}}$}
\put(45.5,14){\tiny ${w_{13}}$}

\put(25.5,18.8){\tiny ${w_{3}}$}
\put(30.5,18.8){\tiny ${w_{4}}$}
\put(35.5,18.8){\tiny ${w_{1}}$}
\put(40.5,18.8){\tiny ${w_{2}}$}
\put(45.5,18.8){\tiny ${w_{3}}$}

\put(25,-4){\tiny  Figure 7 : {\boldmath T(4, 4, 2) : $O_7$}}
\end{picture}

\begin{picture}(0,0)(-65,-17)
\tiny
\setlength{\unitlength}{1.25mm}

\drawline[AHnb=0](0,0)(40,0)
\drawline[AHnb=0](5,10)(45,10)
\drawline[AHnb=0](10,20)(50,20)

\drawline[AHnb=0](0,0)(10,20)
\drawline[AHnb=0](10,0)(20,20)
\drawline[AHnb=0](20,0)(30,20)
\drawline[AHnb=0](30,0)(40,20)
\drawline[AHnb=0](40,0)(50,20)

\drawline[AHnb=0](2.5,5)(5,0)
\drawline[AHnb=0](7.5,15)(15,0)
\drawline[AHnb=0](15,20)(25,0)
\drawline[AHnb=0](25,20)(35,0)
\drawline[AHnb=0](35,20)(42.5,5)
\drawline[AHnb=0](45,20)(47.5,15)

\put(0,-1.3){\tiny ${v_{1}}$}
\put(5,-1.3){\tiny ${v_{2}}$}
\put(10,-1.3){\tiny ${v_{3}}$}
\put(15,-1.3){\tiny ${v_{4}}$}
\put(20,-1.3){\tiny ${v_{5}}$}
\put(25,-1.3){\tiny ${v_{6}}$}
\put(30,-1.3){\tiny ${v_{7}}$}
\put(35,-1.3){\tiny ${v_{8}}$}
\put(40,-1.3){\tiny ${v_{1}}$}

\put(3,5){\tiny ${v_{1, 2, 16, 9}}$}
\put(13,5){\tiny ${v_{3, 4, 10, 11}}$}
\put(23,5){\tiny ${v_{5, 6, 12, 13}}$}
\put(33,5){\tiny ${v_{7, 8, 14, 15}}$}
\put(43,5){\tiny ${v_{1, 2, 16, 9}}$}

\put(5,9){\tiny ${v_{9}}$}
\put(10.5,9){\tiny ${v_{10}}$}
\put(15,9){\tiny ${v_{11}}$}
\put(20.5,9){\tiny ${v_{12}}$}
\put(25,9){\tiny ${v_{13}}$}
\put(30.5,9){\tiny ${v_{14}}$}
\put(35,9){\tiny ${v_{15}}$}
\put(40.5,9){\tiny ${v_{16}}$}
\put(45,9){\tiny ${v_{9}}$}

\put(8,15){\tiny ${v_{9, 10, 6, 7}}$}
\put(18,15){\tiny ${v_{11, 12, 8, 1}}$}
\put(28,15){\tiny ${v_{13, 14, 2, 3}}$}
\put(38,15){\tiny ${v_{15, 16, 4, 5}}$}
\put(48,15){\tiny ${v_{9, 10, 6, 7}}$}

\put(10,20.7){\tiny ${v_{7}}$}
\put(15,20.7){\tiny ${v_{8}}$}
\put(20,20.7){\tiny ${v_{1}}$}
\put(25,20.7){\tiny ${v_{2}}$}
\put(30,20.7){\tiny ${v_{3}}$}
\put(35,20.7){\tiny ${v_{4}}$}
\put(40,20.7){\tiny ${v_{5}}$}
\put(45,20.7){\tiny ${v_{6}}$}
\put(50,20.7){\tiny ${v_{7}}$}

\put(15,-5) {\tiny  Figure 10 : {\boldmath T(8, 2, 6)}}

\end{picture}

\begin{picture}(0,0)(-6,20)
\tiny
\setlength{\unitlength}{1.3mm}
\drawpolygon(-5,5)(10,5)(10,25)(-5,25)
\drawpolygon(15,5)(25,5)(25,25)(15,25)


\drawline[AHnb=0](0,5)(0,25)
\drawline[AHnb=0](5,5)(5,25)
\drawline[AHnb=0](10,5)(10,25)
\drawline[AHnb=0](15,5)(15,25)
\drawline[AHnb=0](20,5)(20,25)

\drawline[AHnb=0](-5,10)(10,10)
\drawline[AHnb=0](15,10)(25,10)
\drawline[AHnb=0](-5,15)(10,15)
\drawline[AHnb=0](15,10)(25,10)
\drawline[AHnb=0](-5,20)(10,20)
\drawline[AHnb=0](-5,10)(10,10)

\drawline[AHnb=0](-5,20)(0,25)
\drawline[AHnb=0](0,20)(5,15)
\drawline[AHnb=0](5,20)(10,25)
\drawline[AHnb=0](15,20)(20,25)
\drawline[AHnb=0](20,20)(25,15)
\drawline[AHnb=0](25,20)(25,25)
\drawline[AHnb=0](15,20)(25,20)

\drawline[AHnb=0](-5,10)(0,15)
\drawline[AHnb=0](0,10)(5,5)
\drawline[AHnb=0](5,10)(10,15)
\drawline[AHnb=0](15,10)(20,15)
\drawline[AHnb=0](20,10)(25,5)
\drawline[AHnb=0](25,10)(25,5)
\drawline[AHnb=0](15,15)(25,15)

\drawline [AHnb=0](10,5)(15,5)
\drawline [AHnb=0](10,10)(15,10)
\drawline [AHnb=0](10,15)(15,15)
\drawline [AHnb=0](10,20)(15,20)
\drawline [AHnb=0](10,25)(15,25)
\drawline [AHnb=0](10,10)(15,5)
\drawline [AHnb=0](10,20)(15,15)

\drawline[linecolor=black,linewidth=.1,AHnb=0](-5,25)(5,25)
\drawline[linecolor=black,linewidth=.1,AHnb=0](5,25)(5,5)
\drawline[linecolor=black,linewidth=.1,AHnb=0](5,25)(25,25)

\put(-5,4){\tiny ${v_{1}}$}
\put(0,4){\tiny ${v_{2}}$}
\put(5,4){\tiny ${v_{3}}$}
\put(10,4){\tiny ${v_{4}}$}
\put(15,4){\tiny ${v_{5}}$}
\put(20,4){\tiny ${v_{6}}$}
\put(25,4){\tiny ${v_{1}}$}

\put(-4.5,9){\tiny ${w_{1}}$}
\put(0.5,10.4){\tiny ${w_{2}}$}
\put(5.5,9){\tiny ${w_{3}}$}
\put(10.5,10.4){\tiny ${w_{4}}$}
\put(15.5,9){\tiny ${w_{5}}$}
\put(20.5,10.4){\tiny ${w_{6}}$}
\put(25.5,9.1){\tiny ${w_{1}}$}

\put(-4.5,14){\tiny ${x_{1}}$}
\put(0.5,14){\tiny ${x_{2}}$}
\put(5.5,14){\tiny ${x_{3}}$}
\put(10.5,14){\tiny ${x_{4}}$}
\put(15.5,14){\tiny ${x_{5}}$}
\put(20.5,14){\tiny ${x_{6}}$}
\put(25.5,14.1){\tiny ${x_{1}}$}

\put(-4.5,19){\tiny ${z_{1}}$}
\put(0.5,20.4){\tiny ${z_{2}}$}
\put(5.5,19){\tiny ${z_{3}}$}
\put(10.5,20.4){\tiny ${z_{4}}$}
\put(15.5,19){\tiny ${z_{5}}$}
\put(20.5,20.4){\tiny ${z_{6}}$}
\put(25.5,19){\tiny ${z_{1}}$}

\put(-5,25.5){\tiny ${v_{3}}$}
\put(0,25.5){\tiny ${v_{4}}$}
\put(5,25.5){\tiny ${v_{5}}$}
\put(10,25.5){\tiny ${v_{6}}$}
\put(15,25.5){\tiny ${v_{1}}$}
\put(20,25.5){\tiny ${v_{2}}$}
\put(25,25.5){\tiny ${v_{3}}$}

\put(4,0){\tiny  Figure 5 : {\boldmath T(6, 4, 2)}}
\end{picture}

\begin{picture}(0,0)(-42,16)
\tiny
\setlength{\unitlength}{1.3mm}
\drawpolygon(5,5)(20,5)(20,25)(5,25)
\drawpolygon(25,5)(40,5)(40,25)(25,25)


\drawline[AHnb=0](10,5)(10,25)
\drawline[AHnb=0](15,5)(15,25)
\drawline[AHnb=0](20,5)(20,25)
\drawline[AHnb=0](25,5)(25,25)
\drawline[AHnb=0](30,5)(30,25)
\drawline[AHnb=0](35,5)(35,25)

\drawline[AHnb=0](5,10)(20,10)
\drawline[AHnb=0](25,10)(40,10)
\drawline[AHnb=0](5,15)(20,15)
\drawline[AHnb=0](25,10)(40,10)
\drawline[AHnb=0](5,20)(20,20)
\drawline[AHnb=0](5,10)(20,10)

\drawline[AHnb=0](5,20)(10,25)
\drawline[AHnb=0](10,20)(15,25)
\drawline[AHnb=0](15,20)(20,25)
\drawline[AHnb=0](25,20)(30,25)
\drawline[AHnb=0](30,20)(35,25)
\drawline[AHnb=0](35,20)(40,25)
\drawline[AHnb=0](25,20)(40,20)

\drawline[AHnb=0](5,10)(10,15)
\drawline[AHnb=0](10,10)(15,15)
\drawline[AHnb=0](15,10)(20,15)
\drawline[AHnb=0](25,10)(30,15)
\drawline[AHnb=0](30,10)(35,15)
\drawline[AHnb=0](35,10)(40,15)
\drawline[AHnb=0](25,15)(40,15)

\drawline [AHnb=0](20,5)(25,5)
\drawline [AHnb=0](20,10)(25,10)
\drawline [AHnb=0](20,15)(25,15)
\drawline [AHnb=0](20,20)(25,20)
\drawline [AHnb=0](20,25)(25,25)
\drawline [AHnb=0](20,10)(25,15)
\drawline [AHnb=0](20,20)(25,25)

\drawline[linecolor=black,linewidth=.25,AHnb=0](20,5)(20,25)
\drawline[linecolor=black,linewidth=.25,AHnb=0](5,25)(20,25)
\drawline[linecolor=black,linewidth=.25,AHnb=0](30,25)(40,25)
\drawline[linecolor=black,linewidth=.2,AHnb=0](30,25)(25,20)
\drawline[linecolor=black,linewidth=.25,AHnb=0](25,20)(25,15)
\drawline[linecolor=black,linewidth=.2,AHnb=0](25,15)(20,10)

\put(5,3.5){\tiny ${v_{1}}$}
\put(10,3.5){\tiny ${v_{2}}$}
\put(15,3.5){\tiny ${v_{3}}$}
\put(20,3.5){\tiny ${v_{4}}$}
\put(25,3.5){\tiny ${v_{5}}$}
\put(30,3.5){\tiny ${v_{6}}$}
\put(35,3.5){\tiny ${v_{7}}$}
\put(38,3.5){\tiny ${v_{1}}$}

\put(5.5,8.8){\tiny ${w_{1}}$}
\put(10.5,8.8){\tiny ${w_{2}}$}
\put(15.5,8.8){\tiny ${w_{3}}$}
\put(20.5,8.8){\tiny ${w_{4}}$}
\put(25.5,8.8){\tiny ${w_{5}}$}
\put(30.5,8.8){\tiny ${w_{6}}$}
\put(35.5,8.8){\tiny ${w_{7}}$}
\put(37.8,8.8){\tiny ${w_{1}}$}

\put(5.5,15.5){\tiny ${x_{1}}$}
\put(7.7,15.5){\tiny ${x_{2}}$}
\put(12.7,15.5){\tiny ${x_{3}}$}
\put(17.7,15.5){\tiny ${x_{4}}$}
\put(22.7,15.5){\tiny ${x_{5}}$}
\put(27.7,15.5){\tiny ${x_{6}}$}
\put(32.7,15.5){\tiny ${x_{7}}$}
\put(38,15.5){\tiny ${x_{1}}$}

\put(5.5,18.8){\tiny ${z_{1}}$}
\put(10.5,18.8){\tiny ${z_{2}}$}
\put(15.5,18.8){\tiny ${z_{3}}$}
\put(20.5,18.8){\tiny ${z_{4}}$}
\put(25.5,18.8){\tiny ${z_{5}}$}
\put(30.5,18.8){\tiny ${z_{6}}$}
\put(35.5,18.8){\tiny ${z_{7}}$}
\put(38,18.8){\tiny ${z_{1}}$}

\put(5,25.5){\tiny ${v_{4}}$}
\put(9.5,25.5){\tiny ${v_{5}}$}
\put(14.5,25.5){\tiny ${v_{6}}$}
\put(19.5,25.5){\tiny ${v_{7}}$}
\put(24.5,25.5){\tiny ${v_{1}}$}
\put(29.5,25.5){\tiny ${v_{2}}$}
\put(34.5,25.5){\tiny ${v_{3}}$}
\put(38.5,25.5){\tiny ${v_{4}}$}

\put(14,0){\tiny  Figure 11 : $T(7, 4, 3)$}
\end{picture}

\begin{picture}(0,0)(-110,11)
\tiny
\setlength{\unitlength}{1.3mm}
\drawpolygon(-5,5)(10,5)(10,25)(-5,25)
\drawpolygon(15,5)(25,5)(25,25)(15,25)


\drawline[AHnb=0](0,5)(0,25)
\drawline[AHnb=0](5,5)(5,25)
\drawline[AHnb=0](10,5)(10,25)
\drawline[AHnb=0](15,5)(15,25)
\drawline[AHnb=0](20,5)(20,25)

\drawline[AHnb=0](-5,10)(10,10)
\drawline[AHnb=0](15,10)(25,10)
\drawline[AHnb=0](-5,15)(10,15)
\drawline[AHnb=0](15,10)(25,10)
\drawline[AHnb=0](-5,20)(10,20)
\drawline[AHnb=0](-5,10)(10,10)

\drawline[AHnb=0](0,20)(-5,25)
\drawline[AHnb=0](5,20)(0,15)
\drawline[AHnb=0](10,20)(5,25)
\drawline[AHnb=0](20,20)(15,25)
\drawline[AHnb=0](25,20)(20,15)
\drawline[AHnb=0](25,20)(25,25)
\drawline[AHnb=0](15,20)(25,20)
\drawline[AHnb=0](10,15)(15,20)

\drawline[AHnb=0](0,10)(-5,15)
\drawline[AHnb=0](5,10)(0,5)
\drawline[AHnb=0](10,10)(5,15)
\drawline[AHnb=0](20,10)(15,15)
\drawline[AHnb=0](25,10)(20,5)
\drawline[AHnb=0](25,10)(25,15)
\drawline[AHnb=0](15,10)(25,10)
\drawline[AHnb=0](10,5)(15,10)

\drawline[AHnb=0](25,10)(25,5)
\drawline[AHnb=0](15,15)(25,15)

\drawline [AHnb=0](10,5)(15,5)
\drawline [AHnb=0](10,10)(15,10)
\drawline [AHnb=0](10,15)(15,15)
\drawline [AHnb=0](10,20)(15,20)
\drawline [AHnb=0](10,25)(15,25)

\put(-5,3.8){\tiny ${w_{1}}$}
\put(0,3.8){\tiny ${w_{2}}$}
\put(5,3.8){\tiny ${w_{3}}$}
\put(10,3.8){\tiny ${w_{4}}$}
\put(15,3.8){\tiny ${w_{5}}$}
\put(20,3.8){\tiny ${w_{6}}$}
\put(23,3.8){\tiny ${w_{1}}$}

\put(-4.5,9){\tiny ${v_{1}}$}
\put(0.5,9){\tiny ${v_{2}}$}
\put(5.5,9){\tiny ${v_{3}}$}
\put(10.5,9){\tiny ${v_{4}}$}
\put(15.5,9){\tiny ${v_{5}}$}
\put(20.5,9){\tiny ${v_{6}}$}
\put(23,10.5){\tiny ${v_{1}}$}

\put(-4.8,15.5){\tiny ${z_{1}}$}
\put(0.5,13.5){\tiny ${z_{2}}$}
\put(5.2,15.5){\tiny ${z_{3}}$}
\put(10.5,13.5){\tiny ${z_{4}}$}
\put(15.2,15.5){\tiny ${z_{5}}$}
\put(20.5,13.5){\tiny ${z_{6}}$}
\put(23,13.5){\tiny ${z_{1}}$}

\put(-4.5,18.8){\tiny ${x_{1}}$}
\put(0.5,18.8){\tiny ${x_{2}}$}
\put(5.5,18.8){\tiny ${x_{3}}$}
\put(10.5,18.8){\tiny ${x_{4}}$}
\put(15.5,18.8){\tiny ${x_{5}}$}
\put(20.5,18.8){\tiny ${x_{6}}$}
\put(23,20.5){\tiny ${x_{1}}$}

\put(-5.5,25.5){\tiny ${w_{3}}$}
\put(-.5,25.5){\tiny ${w_{4}}$}
\put(4.5,25.5){\tiny ${w_{5}}$}
\put(9.5,25.5){\tiny ${w_{6}}$}
\put(14.5,25.5){\tiny ${w_{1}}$}
\put(19.5,25.5){\tiny ${w_{2}}$}
\put(22.5,25.5){\tiny ${w_{3}}$}

\put(4,0){\tiny Figure 12 : $R$}

\end{picture}

\vspace{1.2cm}

\noindent\textbf{Example 2.2 : }
Let $M$ be a semi-equivelar map of type $\{3^2,4, 3, 4\}$ with $16$ vertices on the torus. Similarly as above,  by Lemma \ref{lem32434-rep}, there are three representations of $M$, namely, $T(8, 2, 4),$  $T(4, 4, 0),$ and $T(4, 4, 2)$, see Figure 9, 6, 7 respectively. In $T(8, 2, 4), C_{1,1} = C(u_1, u_2, \cdots, u_8)$ and $C_{1, 2} = C(u_1, u_9, u_5, u_{13})$ are two cycles of type $B_1$ (see definition of type $B_1$ in Section \ref{324341}). 
 In $T(4, 4, 0), C_{2, 1} = C(v_1, v_2, v_3, v_4)$ and $C_{2, 2} = C(v_1, v_5, v_9, v_{13})$ are two cycles of type $B_1$. 
 In $T(4, 4, 2)$, $C_{3, 1}= C(w_1, w_2, w_3, w_4)$ and $C_{3, 2} = C(w_1, w_5, w_9, w_{13},$ $w_3, w_7, w_{11}, w_{15})$ are two cycles of type $B_1$.
In Section \ref{324341}, the cycles of type $B_1$ have at most two different lengths. So, $O_5 \not\cong  O_6$ since $\{$length($C_{1, 1})$, length($C_{1, 2}$)$\} \neq \{$length($C_{2, 1})$, length($C_{2, 2})\}$. Now, $\{$length($C_{1, 1})$, length($C_{1, 2}) \} = \{$length($C_{3, 1})$, length( $C_{3, 2})\}$. We identify boundaries of $O_7$ and cut along the cycle $C_{3, 2} = C(w_1, w_5, w_9,$ $w_{13}, w_3, w_7, w_{11}, w_{15})$ and next along $C_{3, 1}$. Thus, we get a $T(8, 2, 4)$ representation in Figure 8. So, $O_5 \cong O_7$ (see the proof of Lemma \ref{lem32434-iso} for isomorphism map). Therefore, we have two maps of type $\{3^{2}, 4, 3, 4\}$ with $16$ vertices on the torus up to isomorphism.

\section{Maps of type $\{3^{3}, 4^{2}\}$}\label{33421}

Let $M$ be a map of type $\{3^{3}, 4^{2}\}$ on the torus. Through each vertex in $M$ there are three distinct types of paths as follows.

\begin{defn} \label{defn3342-1} Let $P_1 :=P(\cdots, u_{i-1}, u_{i}, u_{i+1}, \cdots)$ be a path in edge graph of $M$. We say $P_1$ of type $A_1$ if all the triangles incident with an inner (degree two in $P_1$) vertex $u_i$ lie on one side and all quadrangles incident with $u_i$ lie on the other side of the subpath $P' = P(u_{i-1}, u_{i}, u_{i+1})$ (as in Figure 13) at $u_{i}$. Since the link of vertex $u_i$ is a cycle and the path $P'$ is a cord of cycle $lk(u_i)$, so, the path $P'$ divides the region into two parts. If $u_{t}$ is a boundary vertex (degree one in $P_1$) of $P_1$ then there is an extended path of $P_1$ where $u_{t}$ is an inner vertex.
\end{defn}

Observe that the link of a vertex in $M$ contains some vertices which are adjacent to the vertex. So, to identify the non adjacent vertices in the link we use bold letters. That is, if a description of link, say $lk(w)$ contains any bold letter $\textit{\textbf{a}}$ then it indicates non adjacent vertex to $w$. For example, in $lk(u_i)$ vertices $\textit{\textbf{a}}$ and $\textit{\textbf{c}}$ are non adjacent to $u_i$. In this article, we consider permutation of vertices in $lk(u_i)$ of a vertex $u_i$ counter clockwise locally at $u_i$. 

\vspace{-.3cm}

\begin{picture}(0,0)(-15,35)

\tiny
\setlength{\unitlength}{1.5mm}
\drawpolygon(5,5)(20,5)(20,10)(5,10)
\drawline[AHnb=0](12.5,5)(12.5,10)
\drawpolygon(12.5,10)(8.5, 15)(5,10)
\drawpolygon(12.5,10)(8.5, 15)(17.5,15)
\drawpolygon(17.5,15)(12.5, 10)(20,10)

\drawline[AHnb=0,dash={1.0 1.0 1.0 1.0}{0.0}](2,10)(5,10)
\drawline[AHnb=0,dash={1.0 1.0 1.0 1.0}{0.0}](20,10)(22.5,10)

\drawline[AHnb=0,dash={1.0 1.0 1.0 1.0}{0.0}](12.5,5)(12.5,2.5)
\drawline[AHnb=0,dash={1.0 1.0 1.0 1.0}{0.0}](12.5,5)(15,2.5)

\drawline[AHnb=0,dash={1.0 1.0 1.0 1.0}{0.0}](17.5,15)(20,17.5)
\drawline[AHnb=0,dash={1.0 1.0 1.0 1.0}{0.0}](8.5,15)(6,18)

\drawline[linecolor=black,linewidth=.2,AHnb=0](5,10)(20,10)
\drawline[linecolor=black,linewidth=.2,AHnb=0](12.5,5)(12.5,10)
\drawline[linecolor=black,linewidth=.2,AHnb=0](12.5,10)(8.5, 15)
\drawline[linecolor=black,linewidth=.2,AHnb=0](12.5,10)(17.5,15)

\put(5.5,8.8){\tiny ${u_{i-1}}$}
\put(13,8.8){\tiny  ${u_{i}}$}
\put(16.5,8.8){\tiny ${u_{i+1}}$}

\put(3.5,5){\tiny ${a}$}
\put(11,3.5){\tiny ${b}$}
\put(21,5){\tiny ${c}$}

\put(8.5,15.5){\tiny ${e}$}
\put(16.7,15.5){\tiny ${f}$}

\put(1,11){\tiny ${A_1}$}
\put(20,18){\tiny ${A_2}$}
\put(5,18.5){\tiny ${A_3}$}

\put(6,.5){\tiny  Figure 13 : $lk(u_i)$}


\end{picture}
\vspace{1.3cm}

\begin{picture}(0,0)(0,50)
\tiny
\setlength{\unitlength}{1.5mm}
\drawpolygon(5,5)(20,5)(20,15)(5,15)
\drawpolygon(25,5)(40,5)(40,15)(25,15)


\drawline[AHnb=0](10,5)(10,15)
\drawline[AHnb=0](15,5)(15,15)
\drawline[AHnb=0](30,5)(30,15)
\drawline[AHnb=0](35,5)(35,15)
\drawline[AHnb=0](5,10)(20,10)
\drawline[AHnb=0](25,10)(40,10)

\drawline[AHnb=0](5,10)(10,15)
\drawline[AHnb=0](10,10)(15,15)
\drawline[AHnb=0](15,10)(20,15)
\drawline[AHnb=0](25,10)(30,15)
\drawline[AHnb=0](30,10)(35,15)
\drawline[AHnb=0](35,10)(40,15)

\put(21,10){$\cdots$}
\put(21,5){$\cdots$}
\put(21,15){$\cdots$}

\put(5,3.5){\tiny ${v_{1}}$}
\put(10,3.5){\tiny ${v_{2}}$}
\put(15,3.5){\tiny ${v_{3}}$}
\put(20,3.5){\tiny ${v_{4}}$}
\put(25,3.5){\tiny ${v_{r-2}}$}
\put(30,3.5){\tiny ${v_{r-1}}$}
\put(35,3.5){\tiny ${v_{r}}$}
\put(40,3.5){\tiny ${v_{1}}$}

\put(5.5,8.7){\tiny ${w_{1}}$}
\put(10.5,8.7){\tiny ${w_{2}}$}
\put(15.5,8.7){\tiny ${w_{3}}$}
\put(20.5,8.7){\tiny ${w_{4}}$}
\put(25.5,8.7){\tiny ${w_{r-2}}$}
\put(30.5,8.7){\tiny ${w_{r-1}}$}
\put(35.5,8.7){\tiny ${w_{r}}$}
\put(40.5,8.7){\tiny ${w_{1}}$}

\put(5,15.5){\tiny ${x_{1}}$}
\put(10,15.5){\tiny ${x_{2}}$}
\put(15,15.5){\tiny ${x_{3}}$}
\put(20,15.5){\tiny ${x_{4}}$}
\put(25,15.5){\tiny ${x_{r-2}}$}
\put(30,15.5){\tiny ${x_{r-1}}$}
\put(35,15.5){\tiny ${x_{r}}$}
\put(40,15.5){\tiny ${x_{1}}$}

\put(18,0){\tiny  Figure 14 : Cylinder}

\end{picture}

\begin{picture}(0,0)(-20,90)
\tiny
\setlength{\unitlength}{1.5mm}
\drawpolygon(5,5)(20,5)(20,25)(5,25)
\drawpolygon(25,5)(40,5)(40,25)(25,25)


\drawline[AHnb=0](10,5)(10,25)
\drawline[AHnb=0](15,5)(15,25)
\drawline[AHnb=0](20,5)(20,25)
\drawline[AHnb=0](25,5)(25,25)
\drawline[AHnb=0](30,5)(30,25)
\drawline[AHnb=0](35,5)(35,25)
\drawline[AHnb=0](40,5)(40,25)
\drawline[AHnb=0](45,5)(45,25)
\drawline[AHnb=0](50,5)(50,25)
\drawline[AHnb=0](55,5)(55,25)
\drawline[AHnb=0](60,5)(60,25)

\put(42,5){\ldots}
\put(42,10){\ldots}
\put(42,15){\ldots}
\put(42,20){\ldots}
\put(42,25){\ldots}

\drawline[AHnb=0](5,10)(20,10)
\drawline[AHnb=0](25,10)(40,10)
\drawline[AHnb=0](5,15)(20,15)
\drawline[AHnb=0](25,10)(40,10)
\drawline[AHnb=0](5,20)(20,20)
\drawline[AHnb=0](5,10)(20,10)
\drawline[AHnb=0](45,10)(60,10)
\drawline[AHnb=0](45,20)(60,20)
\drawline[AHnb=0](45,15)(60,15)
\drawline[AHnb=0](45,5)(60,5)
\drawline[AHnb=0](45,25)(60,25)

\drawline[AHnb=0](5,20)(10,25)
\drawline[AHnb=0](10,20)(15,25)
\drawline[AHnb=0](15,20)(20,25)
\drawline[AHnb=0](25,20)(30,25)
\drawline[AHnb=0](30,20)(35,25)
\drawline[AHnb=0](35,20)(40,25)
\drawline[AHnb=0](25,20)(40,20)

\drawline[AHnb=0](5,10)(10,15)
\drawline[AHnb=0](10,10)(15,15)
\drawline[AHnb=0](15,10)(20,15)
\drawline[AHnb=0](25,10)(30,15)
\drawline[AHnb=0](30,10)(35,15)
\drawline[AHnb=0](35,10)(40,15)
\drawline[AHnb=0](25,15)(40,15)

\put(21,10){\ldots}
\put(21,5){\ldots}
\put(21,15){\ldots}
\put(21,20){\ldots}
\put(21,25){\ldots}

\drawline[AHnb=0](45,10)(50,15)
\drawline[AHnb=0](50,10)(55,15)
\drawline[AHnb=0](55,10)(60,15)
\drawline[AHnb=0](45,20)(50,25)
\drawline[AHnb=0](50,20)(55,25)
\drawline[AHnb=0](55,20)(60,25)

\put(5,4){\tiny ${v_{1}}$}
\put(10,4){\tiny ${v_{2}}$}
\put(15,4){\tiny ${v_{3}}$}
\put(20,4){\tiny ${v_{4}}$}
\put(25,4){\tiny ${v_{k}}$}
\put(30,4){\tiny ${v_{k+1}}$}
\put(35,4){\tiny ${v_{k+2}}$}
\put(40,4){\tiny ${v_{k+3}}$}
\put(45,4){\tiny ${v_{r-2}}$}
\put(50,4){\tiny ${v_{r-1}}$}
\put(55,4){\tiny ${v_{r}}$}
\put(60,4){\tiny ${v_{1}}$}

\put(5.5,9){\tiny ${w_{1}}$}
\put(10.5,9){\tiny ${w_{2}}$}
\put(15.5,9){\tiny ${w_{3}}$}
\put(20.5,9){\tiny ${w_{4}}$}
\put(25.5,9){\tiny ${w_{k}}$}
\put(30.5,9){\tiny ${w_{k+1}}$}
\put(35.5,9){\tiny ${w_{k+2}}$}
\put(40.5,9){\tiny ${w_{k+3}}$}
\put(45.5,9){\tiny ${w_{r-2}}$}
\put(50.5,9){\tiny ${w_{r-1}}$}
\put(55.5,9){\tiny ${w_{r}}$}
\put(60.5,9){\tiny ${w_{1}}$}

\put(5.5,14){\tiny ${x_{1}}$}
\put(10.5,14){\tiny ${x_{2}}$}
\put(15.5,14){\tiny ${x_{3}}$}
\put(20.5,14){\tiny ${x_{4}}$}
\put(25.5,14){\tiny ${x_{k}}$}
\put(30.5,14){\tiny ${x_{k+1}}$}
\put(35.5,14){\tiny ${x_{k+2}}$}
\put(40.5,14){\tiny ${x_{k+3}}$}
\put(45.5,14){\tiny ${x_{r-2}}$}
\put(50.5,14){\tiny ${x_{r-1}}$}
\put(55.5,14){\tiny ${x_{r}}$}
\put(60.5,14){\tiny ${x_{1}}$}

\put(5.5,19){\tiny ${z_{1}}$}
\put(10.5,19){\tiny ${z_{2}}$}
\put(15.5,19){\tiny ${z_{3}}$}
\put(20.5,19){\tiny ${z_{4}}$}
\put(25.5,19){\tiny ${z_{k}}$}
\put(30.5,19){\tiny ${z_{k+1}}$}
\put(35.5,19){\tiny ${z_{k+2}}$}
\put(40.5,19){\tiny ${z_{k+3}}$}
\put(45.5,19){\tiny ${z_{r-2}}$}
\put(50.5,19){\tiny ${z_{r-1}}$}
\put(55.5,19){\tiny ${z_{r}}$}
\put(60.5,19){\tiny ${z_{1}}$}

\put(5.5,24){\tiny ${v_{k+1}}$}
\put(10.5,24){\tiny ${v_{k+2}}$}
\put(15.5,24){\tiny ${v_{k+3}}$}
\put(20.5,24){\tiny ${v_{k+4}}$}
\put(25.5,24){\tiny ${v_{n}}$}
\put(30.5,24){\tiny ${v_{1}}$}
\put(35.5,24){\tiny ${v_{2}}$}
\put(40.5,24){\tiny ${v_{3}}$}
\put(45.5,24){\tiny ${v_{k-2}}$}
\put(50.5,24){\tiny ${v_{k-1}}$}
\put(55.5,24){\tiny ${v_{k}}$}
\put(60.5,24){\tiny ${v_{k+1}}$}

\put(25,0){\tiny  Figure 15 : $T(r,4,k)$}

\end{picture}

\begin{picture}(0,0)(-60,45)
\tiny
\setlength{\unitlength}{1.5mm}

\drawline[AHnb=0](5,5)(50,5)

\drawline[AHnb=0](10,5)(10,10)
\drawline[AHnb=0](15,5)(15,10)
\drawline[AHnb=0](20,5)(20,10)
\drawline[AHnb=0](25,5)(25,10)
\drawline[AHnb=0](30,5)(30,10)
\drawline[AHnb=0](35,5)(35,10)
\drawline[AHnb=0](40,5)(40,10)
\drawline[AHnb=0](45,5)(45,10)

\drawline[AHnb=0](7.5,15)(7.5,20)
\drawline[AHnb=0](12.5,15)(12.5,20)
\drawline[AHnb=0](17.5,15)(17.5,20)
\drawline[AHnb=0](22.5,15)(22.5,20)
\drawline[AHnb=0](27.5,15)(27.5,20)
\drawline[AHnb=0](32.5,15)(32.5,20)
\drawline[AHnb=0](37.5,15)(37.5,20)
\drawline[AHnb=0](42.5,15)(42.5,20)
\drawline[AHnb=0](47.5,15)(47.5,20)

\drawline[AHnb=0](10,10)(7.5,15)
\drawline[AHnb=0](10,10)(12.5,15)
\drawline[AHnb=0](15,10)(12.5,15)
\drawline[AHnb=0](15,10)(17.5,15)
\drawline[AHnb=0](20,10)(17.5,15)
\drawline[AHnb=0](20,10)(22.5,15)
\drawline[AHnb=0](25,10)(22.5,15)
\drawline[AHnb=0](25,10)(27.5,15)
\drawline[AHnb=0](30,10)(27.5,15)
\drawline[AHnb=0](30,10)(32.5,15)
\drawline[AHnb=0](35,10)(32.5,15)
\drawline[AHnb=0](35,10)(37.5,15)
\drawline[AHnb=0](40,10)(37.5,15)
\drawline[AHnb=0](40,10)(42.5,15)
\drawline[AHnb=0](45,10)(42.5,15)
\drawline[AHnb=0](45,10)(47.5,15)

\drawline[AHnb=0](10,25)(10,30)
\drawline[AHnb=0](15,25)(15,30)
\drawline[AHnb=0](20,25)(20,30)
\drawline[AHnb=0](25,25)(25,30)
\drawline[AHnb=0](30,25)(30,30)
\drawline[AHnb=0](35,25)(35,30)
\drawline[AHnb=0](40,25)(40,30)
\drawline[AHnb=0](45,25)(45,30)

\drawline[AHnb=0](7.5,20)(10,25)
\drawline[AHnb=0](12.5,20)(10,25)
\drawline[AHnb=0](12.5,20)(15,25)
\drawline[AHnb=0](17.5,20)(15,25)
\drawline[AHnb=0](17.5,20)(20,25)
\drawline[AHnb=0](22.5,20)(20,25)
\drawline[AHnb=0](22.5,20)(25,25)
\drawline[AHnb=0](27.5,20)(25,25)
\drawline[AHnb=0](27.5,20)(30,25)
\drawline[AHnb=0](32.5,20)(30,25)
\drawline[AHnb=0](32.5,20)(35,25)
\drawline[AHnb=0](37.5,20)(35,25)
\drawline[AHnb=0](37.5,20)(40,25)
\drawline[AHnb=0](42.5,20)(40,25)
\drawline[AHnb=0](42.5,20)(45,25)
\drawline[AHnb=0](47.5,20)(45,25)

\drawline[AHnb=0](7.5,35)(7.5,40)
\drawline[AHnb=0](12.5,35)(12.5,40)
\drawline[AHnb=0](17.5,35)(17.5,40)
\drawline[AHnb=0](22.5,35)(22.5,40)
\drawline[AHnb=0](27.5,35)(27.5,40)
\drawline[AHnb=0](32.5,35)(32.5,40)
\drawline[AHnb=0](37.5,35)(37.5,40)
\drawline[AHnb=0](42.5,35)(42.5,40)

\drawline[AHnb=0](10,5)(10,10)
\drawline[AHnb=0](15,5)(15,10)
\drawline[AHnb=0](20,5)(20,10)
\drawline[AHnb=0](25,5)(25,10)
\drawline[AHnb=0](30,5)(30,10)
\drawline[AHnb=0](35,5)(35,10)
\drawline[AHnb=0](40,5)(40,10)
\drawline[AHnb=0](45,5)(45,10)

\drawline[AHnb=0](10,30)(7.5,35)
\drawline[AHnb=0](10,30)(12.5,35)
\drawline[AHnb=0](15,30)(12.5,35)
\drawline[AHnb=0](15,30)(17.5,35)
\drawline[AHnb=0](20,30)(17.5,35)
\drawline[AHnb=0](20,30)(22.5,35)
\drawline[AHnb=0](25,30)(22.5,35)
\drawline[AHnb=0](25,30)(27.5,35)
\drawline[AHnb=0](30,30)(27.5,35)
\drawline[AHnb=0](30,30)(32.5,35)
\drawline[AHnb=0](35,30)(32.5,35)
\drawline[AHnb=0](35,30)(37.5,35)
\drawline[AHnb=0](40,30)(37.5,35)
\drawline[AHnb=0](40,30)(42.5,35)
\drawline[AHnb=0](45,30)(42.5,35)
\drawline[AHnb=0](45,30)(47.5,35)

\drawline[AHnb=0](5,10)(50,10)
\drawline[AHnb=0](5,15)(50,15)
\drawline[AHnb=0](5,20)(50,20)
\drawline[AHnb=0](5,25)(50,25)
\drawline[AHnb=0](5,30)(50,30)
\drawline[AHnb=0](5,35)(50,35)
\drawline[AHnb=0](5,40)(50,40)
\drawline[AHnb=0](5,45)(50,45)
\drawline[AHnb=0](47.5,35)(47.5,40)

\drawline[AHnb=0](7.5,40)(10,45)
\drawline[AHnb=0](12.5,40)(10,45)
\drawline[AHnb=0](12.5,40)(15,45)
\drawline[AHnb=0](17.5,40)(15,45)
\drawline[AHnb=0](17.5,40)(20,45)
\drawline[AHnb=0](22.5,40)(20,45)
\drawline[AHnb=0](22.5,40)(25,45)
\drawline[AHnb=0](27.5,40)(25,45)
\drawline[AHnb=0](27.5,40)(30,45)
\drawline[AHnb=0](32.5,40)(30,45)
\drawline[AHnb=0](32.5,40)(35,45)
\drawline[AHnb=0](37.5,40)(35,45)
\drawline[AHnb=0](37.5,40)(40,45)
\drawline[AHnb=0](42.5,40)(40,45)
\drawline[AHnb=0](42.5,40)(45,45)
\drawline[AHnb=0](47.5,40)(45,45)

\drawline[linecolor=black,linewidth=.2,AHnb=0](25,10)(22.5,15)
\drawline[linecolor=black,linewidth=.3,AHnb=0](22.5,15)(22.5,20)
\drawline[linecolor=black,linewidth=.2,AHnb=0](22.5,20)(20,25)
\drawline[linecolor=black,linewidth=.3,AHnb=0](20,25)(20,30)
\drawline[linecolor=black,linewidth=.2,AHnb=0](20,30)(17.5,35)
\drawline[linecolor=black,linewidth=.3,AHnb=0](17.5,35)(17.5,40)
\drawline[linecolor=black,linewidth=.2,AHnb=0](17.5,40)(15,45)

\drawline[linecolor=black,linewidth=.2,AHnb=0](15,45)(15,46)
\drawline[linecolor=black,linewidth=.2,AHnb=0](20,45)(20,46)
\drawline[linecolor=black,linewidth=.2,AHnb=0](35,45)(35,46)
\drawline[linecolor=black,linewidth=.2,AHnb=0](30,45)(30,46)

\drawline[AHnb=0](10,45)(10,46)
\drawline[AHnb=0](25,45)(25,46)
\drawline[AHnb=0](40,45)(40,46)
\drawline[AHnb=0](45,45)(45,46)

\drawline[AHnb=0](9.3,4)(10,5) \drawline[AHnb=0](10.7,4)(10,5)
\drawline[AHnb=0](14.3,4)(15,5) \drawline[AHnb=0](15.7,4)(15,5)
\drawline[AHnb=0](19.3,4)(20,5) \drawline[AHnb=0](20.7,4)(20,5)
\drawline[AHnb=0](24.3,4)(25,5) \drawline[AHnb=0](25.7,4)(25,5)
\drawline[AHnb=0](29.3,4)(30,5) \drawline[AHnb=0](30.7,4)(30,5)
\drawline[AHnb=0](34.3,4)(35,5) \drawline[AHnb=0](35.7,4)(35,5)
\drawline[AHnb=0](39.3,4)(40,5) \drawline[AHnb=0](40.7,4)(40,5)
\drawline[AHnb=0](44.3,4)(45,5) \drawline[AHnb=0](45.7,4)(45,5)

\put(25.2,5.7){\tiny ${u_{i-1}}$}
\put(25.2,9){\tiny ${u_{i}}$}
\put(20,14){\tiny ${u_{r}}$}
\put(18.5,19){\tiny ${u_{r-1}}$}
\put(15.8,24){\tiny ${u_{r-2}}$}
\put(16,29){\tiny ${u_{r-3}}$}
\put(13.2,34){\tiny ${u_{r-4}}$}
\put(13.2,39){\tiny ${u_{r-5}}$}
\put(11,44){\tiny ${u_{r-6}}$}

\drawline[linecolor=black,linewidth=.2,AHnb=0](25,10)(27.5,15)
\drawline[linecolor=black,linewidth=.3,AHnb=0](27.5,15)(27.5,20)
\drawline[linecolor=black,linewidth=.2,AHnb=0](27.5,20)(30,25)
\drawline[linecolor=black,linewidth=.3,AHnb=0](30,25)(30,30)
\drawline[linecolor=black,linewidth=.2,AHnb=0](30,30)(32.5,35)
\drawline[linecolor=black,linewidth=.3,AHnb=0](32.5,35)(32.5,40)
\drawline[linecolor=black,linewidth=.2,AHnb=0](32.5,40)(35,45)

\put(28,14){\tiny ${u_{i+1}}$}
\put(28,19){\tiny ${u_{i+2}}$}
\put(30.8,24){\tiny ${u_{i+3}}$}
\put(30.3,28.9){\tiny ${u_{i+4}}$}
\put(33.3,33.9){\tiny ${u_{i+5}}$}
\put(32.8,38.9){\tiny ${u_{i+6}}$}
\put(35.8,43.9){\tiny ${u_{i+7}}$}

\put(25.4,24){\tiny ${w_{i-1}}$}
\put(25.2,29){\tiny ${w_{i}}$}
\put(27.7,34){\tiny ${w_{i+1}}$}
\put(27.7,39){\tiny ${w_{i+2}}$}
\put(30.3,44){\tiny ${w_{i+3}}$}

\drawline[linecolor=black,linewidth=.2,AHnb=0](25,30)(22.5,35)
\drawline[linecolor=black,linewidth=.3,AHnb=0](22.5,35)(22.5,40)
\drawline[linecolor=black,linewidth=.2,AHnb=0](22.5,40)(20,45)
\drawline[linecolor=black,linewidth=.3,AHnb=0](25,30)(27.5,35)
\drawline[linecolor=black,linewidth=.3,AHnb=0](27.5,35)(27.5,40)
\drawline[linecolor=black,linewidth=.3,AHnb=0](27.5,40)(30,45)

\put(18.5,34){\tiny ${w_{r-8}}$}
\put(18,39){\tiny ${w_{r-9}}$}
\put(15.5,44){\tiny ${w_{r-10}}$}


\put(24,1.4){\tiny   Figure 16}
\end{picture}

\vspace{8.3cm}

\begin{defn}\label{defn3342-2}
Let $P_2 :=P( \cdots, v_{i-1}, v_{i}, v_{i+1}, \cdots)$ be a path in edge graph of $M$ for which $v_{i}, v_{i+1}$ are two consecutive inner vertices of $P_2$ or an extended path of $P_2$. We say $P_2$ of type $A_2$ if $lk(v_{i})=C(\textbf{a}, v_{i-1}, \textbf{b}, c, v_{i+1}, d, e)$ implies $lk(v_{i+1})=C(\textbf{a}_0, v_{i+2}, \textbf{b}_0, d, v_{i}, c, p)$ and $lk(v_{i})=C(\textbf{x},v_{i+1},\textbf{z}, l, v_{i-1}, k, m)$ implies $lk(v_{i+1})=C(\textbf{l},v_{i}, \textbf{m}, x, v_{i+2}, g, z)$. At least one of the former two conditions must occur for each vertex. 
\end{defn}

\begin{defn}\label{defn3342-3}
Let $P_3:=P(\cdots, w_{i-1}, w_{i}, w_{i+1}, \cdots)$ be a path in edge graph of $M$ for which $w_{i}, w_{i+1}$ are two inner vertices of $P_3$ or an extended path of $P_3$. We say $P_3$ of type $A_3$ if $lk(w_{i})=C(\textbf{a}, w_{i-1}, \textbf{b}, c, d, w_{i+1}, e)$ implies $lk(w_{i+1})=C(\textbf{a}_1, w_{i+2}, \textbf{b}_1, p, e, w_{i}, d)$ and $lk(w_{i})=C(\textbf{a}_2, w_{i+1}, \textbf{b}_2, p, e, w_{i-1}, d)$ implies $lk(w_{i+1})=C(\textbf{p}, w_{i}, \textbf{d}, a_2, z_1, w_{i+2}, b_2)$.
\end{defn}

Let $Q$ be a maximal path (path of maximal length) of type $A_{t}$ for a fixed $t \in \{1, 2, 3\}$. We show that there is an edge $e$ in $M$ such that $Q \cup e$ is a cycle of type $A_t$. So,

\begin{lemma}\label{lem3342-cycle}  If $P(u_{1}, \cdots, u_{r})$ is a maximal path of type $A_1$, $A_2$ or $A_3$ in $M$ then there is an edge $u_{r}u_{1}$ in $M$ such that $C(u_{1}, u_{2}, \cdots, u_{r})$ is a cycle.
\end{lemma}

\begin{proof} Let $Q = P(u_{1}, \cdots, u_{r})$ be of type $A_1$ and $lk(u_{r})=C(\textit{\textbf{x}}, y, \textit{\textbf{z}}, w, v, u, u_{r-1})$.

If $w = u_{1}$ then $Q = C(u_{1},$ $u_{2},$ $\cdots,$ $u_{r})$ is a cycle. 
If $w \neq u_{1}$. Then, either $w = u_i$ for some $2\le i \le r$ or $w \not= u_i$ for all $2 \le i \le r$. 
Suppose $w = u_i$ for some $2\le i \le r$. Observe that $L = P(u_{i-1}, u_i, u_{i+1}) \subset Q$ and $L' = P(u_{r}, w, x)$ are two paths of type $A_1$ through $u_i$. By Definition \ref{defn3342-1}, through each vertex in $M$ we have only one path of this particular type $A_1$. So, $L = L'$. This implies that $u_r = u_{i-1}$ or $u_r = u_{i+1}$. This is a contraction since, by assumption, $Q$ is a path and $u_i \ne u_j$ for all $i \ne j, 1 \le i, j \le r$. Therefore, $w \ne u_i$ for all $2\le i \le r$. So, if $w \not= u_i$ for all $1 \le i \le r$ then by the Definition \ref{defn3342-1}, $u_{r}$ is an inner vertex in the extended path of $Q$. Thus, we get a path namely $Q_1$ which is extended from $Q$. Hence, length$(Q)$ $<$ length$(Q_1)$. This is a contradiction as $Q$ is maximal. Therefore, $w = u_{1}$ and  $Q \cup u_ru_1$ = $C(u_{1}, u_{2},  \cdots, u_{r})$ is a cycle.

Let $W = P(u_{1}, u_{2}, \cdots, u_{r})$ be of type $A_{2}$. We follow similar argument as in $Theorem~1$ \cite{alt:1972}. Let \textit{$lk(u_{r-1})=( \textbf{v}, u_{r}, \textbf{z}, p, q,   u_{r-2}, x)$} and \textit{$lk(u_{r})=(\textbf{p},u_{r-1}, \textbf{x}, v, e, u_{r+1}, z)$}. 

If $u_{r+1} = u_{1}$ then $C(u_{1}, u_{2},  \cdots,  u_{r})$ is a cycle. 
If $u_{r+1} \neq u_{1}$. Then, either $u_{r+1} = u_i$ for some $2\le i \le r$ or $u_{r+1} \not= u_i$ for all $2 \le i \le r$. Suppose $u_{r+1} = u_i$ for some $2\le i \le r$. Then, $u_{r+1} = u_i$ defines a cycle $R = C(u_i, u_{i+1},  \cdots,  u_r)$. Now by assumption, $W$ is a path. That is, $u_i \ne u_j$ for all $1 \le i, j \le r$ and $i \ne j$. By Definition \ref{defn3342-2}, through each vertex in $M$ we have exactly two paths of type $A_2$. Hence, we have either \textit{$lk(u_i) = C(\textbf{a}, u_{i-1}, \textbf{b}, c, u_{i+1}, d, e)$} or \textit{$lk(u_i) = C(\textbf{a}, u_{i-1}, \textbf{b}, c, z, u_{i+1}, e)$}. If \textit{$lk(u_i) = C(\textbf{a}, u_{i-1}, \textbf{b}, c, z, u_{i+1}, e)$} then $u_{i+1} = u_r$. But, $u_i \ne u_j$ for all $1 \le i, j \le r$ and $i \ne j$. Hence, \textit{$lk(u_i) = C(\textbf{a}, u_{i-1}, \textbf{b}, c, u_{i+1}, d, e)$}.  Thus, from the cycles $lk(u_r)$ and $lk(u_i)$, $z = u_{i+1}$, $p = u_{i+2}$, $d = u_r$ and $u_iu_{i+1}u_r$ is a triangle (see Figure 16). Consider cycle $R$ and faces incident to it. So, these faces $u_{r-3}w_iw_{r-8}$, $u_{r-3}w_{r-8}u_{r-4}$, $[u_{r-4}, w_{r-8}, w_{r-9}, u_{r-5}]$, $u_{r-5}w_{r-9}w_{r-10}$, $u_{r-5}w_{r-10}u_{r-6}$, $\cdots$, $w_{i+3}u_{i+7}u_{i+6}$, $w_{i+3}u_{i+6}w_{i+2}$, $[w_{i+2}, u_{i+6}, u_{i+5}, w_{i+1}]$, $w_{i+1}u_{i+5}u_{i+4}$, $w_{i+1}u_{i+4}w_i$ define a new cycle $R' = C(w_i, w_{i+1},$ $\cdots$, $w_{r-9}, w_{r-8})$ (see Figure 16). Observe that, $R'$ is a cycle of same type as $R$ since the faces $[u_{r-2}, u_{r-3}, w_i, w_{i-1}]$ and $[u_{i+3}, u_{i+4}, w_i, w_{i-1}]$ have a common edge $w_{i-1}w_i$, and length($R'$) $<$ length($R$). Similarly, we consider cycle $R'$ and repeat the process as above. Thus, we get a sequence of cycles of same type as $R$. But, in this sequence, the length of cycles is gradually decreasing. So, after finite number of steps cycle of type as $R$ may not exist since the map is finite. Therefore, $u_{r+1} \ne u_i$ for all $2\le i \le r$. By Definition \ref{defn3342-2}, \textit{$lk(u_{r})=(\textbf{z}, u_{r-1}, \textbf{x}, v, u_{r+1}, w)$} implies \textit{$lk(u_{r+1})= C(\textbf{y}, u_{r}, \textbf{x}, a, w, b, c)$} for some vertices {\em b, w}. Hence, we define a new path $L := P(u_1,$ $\cdots$, $u_r) \cup P(u_{r}, u_{r+1})$ which is of type $A_{2}$. So, we have a path $Q$ with length($Q$) $>$ length($P$). This gives a contradiction as $P$ is maximal. Therefore, $u_{r+1} = u_{1}$, that is, $C(u_{1}, u_{2},$ $\cdots$, $u_{r})$ is cycle of type $A_{2}$. 

We use similar argument as above for maximal path of type $A_{3}$. Similarly, we get an edge which defines a cycle of type $A_3$.  This completes the proof.
\end{proof}

So, every maximal path of type $A_1$, $A_2$ or $A_3$ is a cycle. In this article, we use the terminology cycle in place of maximal path since it is a cycle. Let $C_1$ and $C_2$ be two cycles of type $A_{t}$ for a fixed $t \in \{1, 2, 3\}$. Then, 

\begin{lemma}\label{lem3342-hom} (a) If $C_{1}, C_{2}$ are two cycles of same type $A_{1}$ such that $C_{1} \cap C_{2} \neq \emptyset$ then $C_{1} = C_{2}$. (b) If $C'_{1}, C'_{2}$ are two cycles of same type $A_t$ for a fixed $t \in \{2, 3\}$ such that $E(C'_{1}) \cap E(C'_{2}) \neq \emptyset$ then $C'_{1} = C'_{2}$.
\end{lemma}

\begin{proof} Let $C_{1} :=C(u_{1, 1}, u_{1, 2},\cdots,u_{1, r})$ and $C_{2} :=C(u_{2, 1}, u_{2, 2},\cdots, u_{2, s})$ be two cycles of type $A_{1}$. If $C_{1} \cap C_{2} \neq \emptyset$ then $V(C_{1} \cap C_{2}) \not= \emptyset$. Let $w \in V(C_{1} \cap C_{2})$. The cycles $C_{1}$ and $C_{2}$ are both well defined at the common vertex $w$. Let $lk(w) = C(\textit{\textbf{w}}_{1}, w_{2},  \textit{\textbf{w}}_{3}, w_{4}, w_{5}, w_{6}, w_{7})$. By Definition \ref{defn3342-1}, $w_{4}, w_{7} \in V(C_{1} \cap C_{2})$. So, $P(w_{4}, w, w_{7})$ is part of $C_{t}$ for $t \in \{1, 2\}$. Let $w = u_{1,t_{1}} = u_{2,t_{2}}$. Then $w_{4} = u_{1,t_{1}-1} = u_{2,t_{2}-1}$ and $w_{7} = u_{1,t_{1}+1} = u_{2,t_{2}+1}$ for some $t_{1} \in \{1, \cdots, r\}$ and $t_{2}\in \{1, \cdots, s\}$. We can argue for $w_4$ and $w_7$ as we did for $w$ to get two vertices, $u_{1,t_{1}-2} = u_{2,t_{2}-2}$ and $u_{1,t_{1}+2} = u_{2,t_{2}+2}$. This process stops after finite number of steps as $r$ and $s$ both are finite. Let $ r < s $. Then $u_{1, 1} = u_{2, I+1}, u_{1, 2} = u_{2, I+2},$ $\cdots$, $u_{1, r} = u_{2, I+r}$ and $u_{1, 1} = u_{2, I+r+1}$ for some $I \in \{1,\cdots,s\}$. Hence $u_{1, 1} = u_{2, I+1} = u_{2, I+r+1}$. This implies that $I+1 = I+r+1$ and the cycle $C_{2}$ contains a cycle of length $r$. This gives $r = s$. Hence $C_{1} = C_{2}$. 

Let $C'_{1}, C'_{2}$ be two cycles of same type $A_2$ and $E(C'_{1}) \cap E(C'_{2}) \neq \emptyset$. Let $uv \in E(C'_{1}) \cap E(C'_{2})$. We proceed with the vertex $u$ of edge $uv$ and in a similar way as we did for the cycles of type $A_1$. (This argument we have also used in Lemma \ref{lem32434-1}.) Thus, we get $C'_1 = C'_2$. Similarly we argue for the case of cycles of type $A_{3}$ to show that $C'_1 = C'_2$. This completes the proof.
\end{proof}

 Now, we show that the cycle of type $A_{t}$ for each $t \in \{1, 2, 3\}$ is non-contractible.


\begin{lemma}\label{lem3342-non-con}  If a cycle $C$ is of type $A_{t}$ for some $t \in \{1, 2, 3\}$ in $M$ then $C$ is non-contractible.
\end{lemma}

\begin{proof} Let $C$ be a cycle of type $A_{t}$ for a fixed $t \in \{1, 2, 3\}$ in $M$. We claim that the cycle $C$ is non-contractible. Let the cycle $C$ be of type $A_{1}$. Suppose, $C$ is contractible. Let $D_{C}$ be a $2$-disk bounded by the cycle $C$. Let $f_0, f_1$ and $f_2$ denote the number of vertices, edges and faces of $D_{C}$ respectively. Let there be $n$ internal vertices and $m$ boundary vertices. So, $f_0 = n+m$, $f_1 = (5n + 3m)/2$ and $f_2 = n+(n+m)/2$ if quadrangles are incident with $C$, and $f_0 = n+m$, $f_1 = (5n + 4m)/2$ and $f_2 = 3n/2+m$ if triangles are incident with $C$ in $D_{C}$. In both the cases, $f_0 - f_1 + f_2 = 0$. This is not possible since the Euler characteristic of the 2-disk $D_{C}$ is 1.
Therefore, $C$ is non-contractible. 

We argue with the similar argument for the cycles of types $A_{2}$ and $A_{3}$. Suppose a cycle $W$ of type $A_2$ is contractible. Let $D_{W}$ be a $2$-disk which is bounded by the cycle $W$. Similarly as above, calculate $f_0, f_1, f_2$ and we get $f_0 - f_1 + f_2 = 0$, a contradiction. 

We can argue similar way for cycle of type $A_3$. This completes the proof.
\end{proof}

Let $C$ be a cycle of type $A_{t}$ for a fixed $t \in \{1, 2, 3\}$ in $M$. Let $S$ be a set of faces which are incident at $u$ for all $u \in V(C)$. The geometric carrier $|S|$ is a cylinder since $C$ is non-contractible. Let $S_{C} := |S|$. Observe that a {\em cylinder (or an annulus)} in $M$ is a sub-complex of $M$ with two boundary cycles. If the boundary cycles of a cylinder are same, it is a torus and we say that the cylinder has identical boundary components. Clearly, the $S_C$ is a cylinder and has two boundary cycles. Let $\partial S_{C} = \{C_{1}, C_{2}\}$. Then, length($C$) = length($C_{1}$) = length($C_{2}$) by Lemma \ref{lem3342-length}.

\begin{lemma} \label{lem3342-length} If $C$ is a cycle of type $A_{i}$ for a fixed $i \in \{1, 2, 3\}$ such that $S_{C}$ is a cylinder and $\partial S_{C} = \{C_{1}, C_{2}\}$ then length($C$) = length($C_{1}$) = length($C_{2}$).
\end{lemma}

\begin{proof} Let $C$ be a cycle of type $A_{1}$. Let $F_{1},$ $F_{2}, \cdots,F_{r}$ be a sequence of faces in order which are incident with $C$ and lie on one side of $C$. These faces $F_{1},$ $F_{2}, \cdots, F_{r}$ are also incident with $C_{t}$ and lie on one side of $C_{t}$ for a fixed $t \in \{1, 2\}$. Without loss of generality we assume that $C_t = C_1$. Let $\widehat{F}_{1},$ $\widehat{F}_{2}, \cdots,\widehat{F}_{r}$ denote the sequence of faces incident with $C$ that lie on the other side of $C$. Let $\widehat{F}_i$ be a $\widehat{d}_i$-gon. Then, we get $T_1 := \{\widehat{d}_1, \widehat{d}_2, \cdots, \widehat{d}_r\}$ the sequence of face-types corresponding to the sequence  $\widehat{F}_{1},$ $\widehat{F}_{2}, \cdots,\widehat{F}_{r}$. Again, let $W_{1},$ $W_{2}, \cdots,W_{r}$ denote the sequence of faces incident with $C_1$ that lie on the other side of $C_1$. Similarly, let $T_2 := \{d_1, d_2, \cdots, d_r\}$ for $d_i$-gon $W_i$, $i =1, 2, \cdots, r$. Since $F_{1},$ $F_{2}, \cdots, F_{r}$ is a sequence of faces that lie on one side of both $C$ and $C_1$, so, there exists a $j$ such that $\widehat{d}_1 = d_{j},$ $\widehat{d}_2= d_{j+1}, \cdots, \widehat{d}_{k-j+1} = d_{k},$ $\widehat{d}_{k-j+2} = d_{1}, \cdots, \widehat{d}_k = d_{j-1}$. So, the cycle $C_{1}$ is of type $A_{1}$. Similarly we argue as above for $C_2$ and we get that the cycle $C_{2}$ is of type $A_{1}$. For example in Figure 15, cycle $C = C(x_{1},\cdots,x_{r})$, $\partial S_{C}=\{C_1(w_{1},\cdots,w_{r}), C_2(z_{1},\cdots, z_{r})\}$, and $C,$ $C_1$ and $C_2$ are cycles of same type $A_1$. 

Similarly we repeat above argument for the other two types $A_{j}$ for $j \in \{2, 3\}$ and we get the similar results. So, the boundary cycles of $S_C$ for a cycle $C$ of type $A_i$, $i \in \{1, 2, 3\}$ are also of type $A_i$.

Suppose length($C$) $\neq$ length($C_{1}$) $\neq$ length($C_{2}$). Let $C := C(u_{1},$ $u_{2},\cdots, u_{r})$, $C_{1} :=C(v_{1},$ $v_{2},$ $\cdots, v_{s})$ and $C_{2} :=C(w_{1}, w_{2},\cdots w_{l})$ denote three cycles of type $A_{1}$. The link $lk$($u_1$) contains the vertices $v_1$ and $w_1$. Let  $P(v_{1}, u_{1}, w_{1})$  be a path of type either $A_{2}$ or $A_{3}$ through $u_1$.  Without loss of generality, we assume $r < s$ and $r \not= l$ since $r \not= s \not= l$. Now, the path $P(v_{1}, u_{1}, w_{1})$ is a shortest path between $v_{1}$ and $w_{1}$ via $u_{1}$. So it follows that the path $P(v_{i}, u_{i}, w_{i})$ is also a shortest path of type either $A_{2}$ or $A_{3}$ between $v_{i}$ and $w_{i}$ via $u_{i}$. Since, by assumption $r < s$, so, the link $lk$($u_r$) contains the vertices $v_r,$ $v_{r+j}$ for $j(>0)$ and $w_{r}$. This gives that the link of $u_{r}$ is different from the link of $u_{r-1}$. This is a contradiction as $M$ is a semi-equivelar map. Therefore, $r = s = l$, that is, length($C$) = length($C_{1}$) = length($C_{2}$). 

Similarly we argue for the cycles of type $A_{j}$ for $j \in \{2, 3\}$. This completes the proof.
\end{proof}

Let $C_1$ and $C_2$ be two cycles of same type in a semi-equivelar map $M$ on the torus. We denote $S_{C_1, C_2}$ a cylinder if the boundary components are $C_1$ and $C_2$. We say that the cycle $C_1$ is {\em homologous} to $C_2$ if $S_{C_1, C_2}$ exists. In particular, if $C_1 = C_2$ then consider $S_{C_1, C_2} = C_1 (= C_2)$, and hence, $C_1$ is homologous to $C_2$.

So by above lemma, the cycle $C$ and the boundary cycles of $S_C$ are homologous. Let $C_{1},$ $C_{2}, \cdots, C_{m}$ be a list of cycles which are homologous to $C$ in $M$. Then, the cycles have same length. That is,

\begin{lemma} \label{lem3342-all-length} If $C_{1}, C_{2}, \cdots, C_{m}$ are homologous cycles of type $A_{t}$ for a fixed $t \in \{1, 2, 3\}$ then length($C_{i}$) = length($C_{j}$) for $1 \le i, j \le m$.
\end{lemma}

\begin{proof} Let $C_{i}$ be a cycle of type $A_{1}$. Then, we have a cylinder $S_{C_{i}}$. Let $C_{t_1}, C_{t_2}, \cdots, C_{t_m}$ denote a sequence of cycles $\{C_{1}, C_{2},  \cdots,  C_{m}\}$ such that $\partial S_{C_{t_j}} = \{C_{t_{j-1}}, C_{t_{j+1}}\}$ for $2 \le j \le (m-1)$. So, by Lemma \ref{lem3342-length}, length($C_{t_1}$) =  length($C_{t_j}$) for $j \in \{1,$ $2, \cdots, m\}$. Thus, length($C_{i}$) =  length($C_{j}$) for $1 \le i, j \le m$. Similarly we argue for the cycles of type $A_{j}$ for $j \in \{2, 3\}$. This completes the proof.
\end{proof}

\textbf{A (r, s, k)-representation in $M$ :} Let $v \in V(M)$. By Lemma \ref{lem3342-cycle}, there are three cycles of types $A_1, A_2, A_2$ through $v$. Let $L_{1}, L_{2}, L_{3}$ be three cycles through the vertex $v$ where the cycle $L_i$ is of type $A_{i}$, $i = 1, 2, 3$. Let $L_1 := C(a_1, a_2, \cdots, a_r)$. We cut map $M$ along the cycle $L_{1}$. We get a cylinder which is bounded by identical cycle $L_1$. We denote it by $N_1$. We call that a cycle is a {\em horizontal cycle} if the cycle is $L_1$ or homologous to $L_1$. Similarly, we say that a cycle is a {\em vertical cycle} if the cycle is $L_i$ or homologous to $L_i$ for $i \in \{2, 3\}$. Observe that the horizontal and vertical cycles are non-contractible by Lemma \ref{lem3342-non-con}. Again, we say that a path is a {\em vertical path} if the path is part of a vertical cycle. We consider $N_1$ and make another cut in $N_1$ starting through the vertex $v$ along a path $Q \subset L_{3}$ until reaching $L_{1}$ again for the first time where the starting adjacent face to the horizontal cycle $L_{1}$. 
(For example in Figure 15, $v = v_{1}$.) Assume that along $P := P(w_1 (=a_1), w_2, \cdots, w_m) \subset L_3$ we took the second cut in $N_1$. Thus, we get a planar representation which is denoted by $N_2$. 

\smallskip

\noindent $Claim~1.$  The representation $N_2$ is connected.

\smallskip

Observe that the $N_1$ is connected as $L_1$ is a non-contractible cycle. Suppose $N_2$ is disconnected. This implies that there exists a $2$-disk namely $D_{P_1\cup Q_1}$ which is bounded by a cycle  $P_1 \cup Q_1 = P(u_j, \cdots, u_i) \cup P(a_t \cdots, a_s)$ where $P_1 \subset L_3$, $Q_1 \subset L_1$, $u_i = a_t$ and $u_j = a_s$. We consider faces which are incident with $P_1$ and $Q_1$ in $D_{P_1 \cup Q_1}$. (In this article, $\Box$ represents quadrangular face and $\triangle$ represents triangular face.) Observe that if the quadrangular faces are incident with $Q_1$ and $\Box_i, \triangle_{i, 1}, \triangle_{i, 2}, \Box_{i+1}, \triangle_{i+1, 1}, \triangle_{i+1, 2}, \cdots, \triangle_{j-1, 1}, \triangle_{j-1, 2}, \Box_j$ are incident with $P_1$ in $D_{P_1 \cup Q_1}$ then as in Lemma \ref{lem3342-non-con}, we calculate number of vertices $f_0$, edges $f_1$ and faces $f_2$ of $D_{P_1\cup Q_1}$. So, we get $f_0 - f_1 + f_2 = 0$. Similarly, if the triangular faces are incident with $Q_1$ then also we calculate $f_0, f_1$ and $f_2$ of $D_{P_1\cup Q_1}$. Similarly, we get $f_0 - f_1 + f_2 = 0$. Which is a contradiction in both the case as $D_{P_1\cup Q_1}$ is $2$-disk. So, $N_2$ is connected. 

Observe that $N_2$ is planer and bounded by $L_1$ and $Q$. Let $s$ denote the number of cycles which are homologous to $L_{1}$ along $P$ in $N_2$.
(For example in Figure 15, we took second cut along the path $P(v_{1},$ $w_{1},$ $x_{1},$ $z_{1},$ $v_{k+1})$ which is part of $L_{3}$ and $s = 4$.)
Now, in $N_2$, length($L_{1}$) = $r$ and number of horizontal cycles along $P$ is $s$. So, we denote $N_2$ by {\em $(r, s)$-representation}. 

Observe that $N_2$ is bounded by cycle $L_1$, path $Q$, cycle $L_1'$ and path $Q'$ where $L_1 = L_1'$ and $Q = Q'$.  We say that the cycle $L_1$ is a {\em lower (base) horizontal cycle} and $L_1'$ is an {\em upper horizontal cycle} in $(r, s)$-representation. Without loss, we may assume that the incident faces of $L_1$ are quadrangles. (For example in Figure 15, $C(v_1, \cdots, v_r)$ is a lower horizontal cycle and $C(v_{k+1}, \cdots, v_k)$ is an upper horizontal cycle.) So, we get identification of vertical sides of $N_2$ in the natural manner but the identification of the horizontal sides needs some shifting so that a vertex in the lower(base) side is identified with a vertex in the upper side. 
Let $L_1' = C(a_{k+1}(=w_m), \cdots, a_k)$. Then, $a_{k+1}$ is the starting vertex in $L_1'$. In $(r, s)$-representation, let $k$ $=:$ length$(P(a_1, \cdots ,a_{k+1}))$ if $P(a_1, \cdots ,a_{k+1})$ is part of $L_1$. Thus, we get a new {\em $(r, s, k)$-representation} of the $(r, s)$-representation. We call boundaries of $(r, s, k)$-representation are the cycles and paths along which we took the cuts to construct $(r, s, k)$-representation of $M$. (For example in Figure 15, the vertex $v_{k+1}$ is the starting vertex of the upper horizontal cycle $C(v_{k+1},v_{k+2},\cdots,v_{k})$ and $k = $ length($P(v_{1}, v_{2}, \cdots, v_{_{k+1}}$)).) 
\hfill$\Box$

\smallskip

By the above construction, we see that every map of type $\{3^{3}, 4^{2}\}$ has a $(r, s, k)$-representation. 
We use $T(r, s, k)$ to represent $(r, s, k)$-representation. Therefore, we have the following lemma.

\begin{lemma} \label{lem3342-rep}  The map of type $\{3^{3}, 4^{2}\}$ on the torus has a $T(r, s, k)$ representation.
\end{lemma}

Let $T(r, s, k)$ be a representation of $M$. It has two identical upper and lower horizontal cycles of type $A_1$, namely, $C_{lh} := C(u_{1}, u_{2}, \cdots, u_{r})$ and $C_{uh} := C(u_{k+1}, u_{k+2}, \cdots, u_{k})$ in $T(r, s, k)$ respectively. (For example in Figure 11, $C_{lh} = C(v_{1}, v_{2},$ $ \cdots, v_{7})$ and $C_{uh} = C(v_{4}, v_{5}, \cdots, v_{3})$.) We define a new cycle in $T(r, s, k)$ using $C_{lh}$ and $C_{uh}$. The vertex $u_{k+1} \in V(C_{lh})$ is the starting vertex of $C_{uh}$ in $T(r, s, k)$. In $T(r, s, k)$, we define two paths $Q_{2} = P(u_{k+1}, \cdots, u_{k_{1}})$ of type $A_{2}$ and $Q_{3} = P(u_{k+1}, \cdots, u_{k_{2}})$ of type $A_{3}$ through $u_{k+1}$ in $T(r, s, k)$ where $u_{k_{1}}, u_{k_{2}} \in V(C_{uh})$. Clearly, the paths $Q_2$ and $Q_3$ are not homologous to horizontal cycles, that is, these are not part of cycles of type $A_1$. We define two edge disjoint paths $Q_{2}' = P(u_{k_{1}}, \cdots, u_{k+1})$ and $Q_{3}' = P(u_{k_{2}}, \cdots, u_{k+1})$ in $C_{uh}$ where $Q_{2}' \cup Q_{3}' :=P(u_{k_{1}}, \cdots, u_{k+1}, \cdots, u_{k_{2}}) \subset C_{uh}$ is a path in $T(r, s, k)$. Let $C_{4, 1} := Q_{2}'\cup Q_2 := C(u_{k+1}, \cdots, u_{k_{1}}, \cdots, u_{k+1})$ and $C_{4, 2} := Q_{3}'\cup Q_3 :=C(u_{k+1}, \cdots, u_{k_{2}}, \cdots, u_{k+1})$. So, we define a new cycle $C_4$ using lengths of $C_{4, 1}$ and $C_{4, 2}$ as follows: 
\begin{equation}
C_{4} := \left\{
\begin{array}{rl}
C_{4, 1}, &  \mbox{if length}(C_{4, 1}) \leq \mbox{length}(C_{4, 2}),\\
C_{4, 2}, &  \mbox{if length}(C_{4, 1}) > \mbox{length}(C_{4, 2}).
\end{array}
\right.
\end{equation}
It follows from the definition of $C_4$ that length($C_{4}$) := min$\{$length($Q_{2}'$) + length($Q_{2}$), length($Q_{3}'$) + length($Q_{3}$ )$\}$ = min$\{k+s, (r-\frac{s}{2}-k)(mod~r)+s\}$. We say that the cycle $C_{4}$ is of type $A_{4}$ in $T(r, s, k)$. (In this section, we use $(r+\frac{s}{2}-k)$ in place of $(r+\frac{s}{2}-k)(mod~r)$.) (For example in Figure 11, $C_{4, 1} = Q_{2}'\cup Q_{2} = P(v_{4}, v_{3}, v_{2}) \cup P(v_{2}, z_{5}, x_{5}, w_{4}, v_{4})$ and $C_{4, 2} = Q_{3}'\cup Q_{3} = P(v_{4}, v_{5}, v_{6}, v_{7}) \cup P(v_{7}, z_{4}, x_{4}, w_{4}, v_{4})$.) So, we have the cycles of four types $A_1, A_2, A_3$ and $A_4$ in $T(r, s, k)$. 

We show that the cycles of type $A_{1}$ have same length and the cycles of type $A_{2}$ have at most two different lengths in $M$. So, 

\begin{lemma}\label{lem3342-hom-length}  In $M$, the cycles of type $A_{1}$ have unique length and the cycles of type $A_{2}$ have at most two different lengths.
\end{lemma}

\begin{proof} Let $C_1$ be a cycle of type $A_1$ in $M$. By the preceding argument of this section, the geometric carrier $S_{C_1}$ of the faces which are incident with $C_1$ is a cylinder and $\partial S_{C_1} := \{C_2, C_0\}$ where the cycle $C_2$ is homologous to $C_1$ and length($C_1$) = length($C_2$). Similarly as above, the $S_{C_2}$ is a cylinder and which is bounded by two homologous cycles $C_1$ and another, say $C_3$ of type $A_1$. Again, we consider the cycle $C_3$ and continue with above process. In this process, let $C_i$ denote a cycle at $i^{th}$ step such that $\partial S_{C_{i}} = \{C_{i-1}, C_{i+1}\}$ and length($C_{i-1}$) = length($C_i$) = length$(C_{i+1})$. Since $M$ consists of finite number of vertices, it follows that this process stops after, say $t+1$ number of steps when the cycle $C_0$ appears in this process. Thus, we get $C_1$, $C_2$, $\cdots$, $C_t$ cycles of type $A_1$ which are homologous to $C_1$ and cover all the vertices of $M$ as the vertices of $S_{C_i}$ are the vertices of $C_{i-1} \cup C_i \cup C_{i+1}$ for $1 \leq i \leq t$.  It is clear from the definition that there is only one cycle of type $A_1$ through each vertex in $M$. Therefore, the cycles $C_1$, $C_2$, $\cdots$, $C_t$ are the only cycles of type $A_1$ in $M$. Since these cycles are homologous to each other, it follows that length($C_1$) = length($C_i$) for $i \in \{1, \cdots, t\}$ by Lemma \ref{lem3342-all-length}. This implies that the cycles of type $A_{1}$ have unique length in $M$. 

Let $L_1, L_2, L_3$ be three cycles through a vertex of types $A_1$, $A_2, A_3$ respectively in $M$. We repeat above process and consider $L_2$ in place of $C_1$. Similarly, here we get a sequence of cycles, namely, $R_1 (=L_2), R_2, \cdots, R_k$ of type $A_2$ which are homologous to each other. Since $R_i$ and $R_j$ are homologous to each other for $1 \le i, j \le k$, it follows that $l_1 = $ length$(L_2)$ = length $(R_{i})$ for $1 \leq i \leq k$ by Lemma \ref{lem3342-all-length}. Similarly, again we consider the cycle $L_3$ and repeat above argument. Let $l_2 =$ length($L_3$). Since the cycles $L_2$ and $L_3$ are mirror image of each other, it follows that they define same type of cycles. So, the map $M$ contains the cycles of type $A_{2}$ of lengths $l_1$ and $l_2$. Therefore, the cycles of type $A_{2}$ have at most two different lengths in $M$. This completes the proof.
\end{proof}

We define admissible relations among $r$, $s$, $k$ of $T(r, s, k)$ such that $T(r, s, k)$ represents a map after identifying their boundaries.


\begin{lemma}\label{lem3342-all-rep} The maps of type $\{3^3, 4^2\}$ of the form $T(r, s, k)$ exist if and only if the following holds : (i) $s \geq  2$ even, (ii) $rs \geq  10$, (iii) $2 \le k \le r-3$ if $s = 2 $ $\&$ $ 0 \leq k \leq r-1$ if $s \ge 4$.
\end{lemma}

\begin{proof} Let $T(r, s, k)$ be a representation of $M$. In $T(r, s, k)$, the $s$ denote the number of horizontal cycles of type $A_{1}$. By the preceding argument of this section and Lemma \ref{lem3342-hom-length}, the cycles of type $A_1$ are homologous to  each other, cover all the vertex of $M$ and have length $r$. So, the number of vertex $n$ of $M$ = length of the cycle of type $A_1 \times$ the number of cycles of type $A_1$ = $rs$. 

Let $C$ be a cycle of type $A_{1}$. By the definition of $A_{1}$, the triangles incident with $C$ lie on one side and $4$-gons lie on the other side of $C$. If $s = 1$ then $T(r, 1, k)$ contains one horizontal cycle namely $C$. Since the incident faces of $C$ are either triangles or $4$-gons, it implies that the faces of $M$ are either only $3$-gons or $4$-gons. This is a contradiction as $M$ consists of both types of faces. So, $s \ge 2$ for all $r$. If $s$ is not an even integer and  $C$ is the base horizontal cycle in $T(r, s, k)$ then the incident faces of a vertex in $C$ are all $3$-gons or $4$-gons after identification of the boundaries of $T(r, s, k)$. This is a contradiction as both $3$- and $4$-gons are incident at each vertex of $M$. So, $s$ is even.

If $s = 2$ and $r < 5$ then the representation $T(4, 2, k)$ has two horizontal cycles. If $C(u_{1}, u_{2}, u_{3}, u_{4})$ and $C(u_{5}, u_{6}, u_{7}, u_{8})$ are two horizontal cycles in $T(4, 2, k)$ then the link $lk(u_{6})$ is not a cycle. So, $r \neq 4$. Similarly one can see that $r \neq 1, 2, 3$. Thus, $r \geq 5$. 
If $s \ge 4$ and $r < 3$ then one can see as above that some vertex has link which is not a cycle. So, by combining above all cases, $r \geq  3$ and $rs \geq  10$.

If $s = 2$ and $k \in \{1, \cdots, r-1\} \setminus \{2, \cdots, r-3\}$ then we proceed with as above and we get some vertex whose link is not a cycle. Thus, $s = 2$ implies $ k \in \{2, \cdots, r-3\}$. Similarly we repeat the above argument for $s \geq 4$ and we get that $k \in \{1, \cdots, r-1\}$ if $s \geq 4$. This completes the proof.
\end{proof}

Let $M_{1}$ and $M_{2}$ be two maps of type $\{3^{3}, 4^{2}\}$ with same number of vertices on the torus and $T_i := T(r_{i}, s_{i}, k_{i})$, $i \in \{1, 2\}$ denote $M_{i}$. If $a_{i, 1}$ = length of the cycle of type $A_1$, $a_{i, 2}$ = length of the cycle of type $A_2$, $a_{i, 3}$ = length of the cycle of type $A_3$ and $a_{i, 4}$ = length of the cycle of type $A_4$ in $T_i$ then we say that $T(r_i, s_i, k_i)$ has cycle-type $(a_{i, 1}, a_{i, 2}, a_{i, 3}, a_{i, 4})$ if $a_{i, 2} \le a_{i, 3}$ or $(a_{i, 1}, a_{i, 3}, a_{i, 2}, a_{i, 4})$ if $a_{i, 3} < a_{i, 2}$.
Now, we show the following isomorphism lemma. 

\begin{lemma}\label{lem3342-iso}
The map $M_{1} \cong M_{2}$ if and only if they have same cycle-type.
\end{lemma}

\begin{proof} We first assume that the maps $M_{1}$ and $M_{2}$ have same cycle-type. This gives that $a_{1, 1}$ = $a_{2, 1}$, $\{a_{1, 2}, a_{1, 3}\}$ = $\{a_{2, 2}, a_{2, 3}\}$ and $a_{1, 4} = a_{2, 4}$. The maps $M_i$ for $i \in \{1, 2\}$ have a $T_i = T(r_i, s_i, k_i)$ representation.

\smallskip

\noindent $Claim.$  $T_1 \cong T_2$.

\smallskip

The $T_1$ has $s_1$ number of horizontal cycles of type $A_1$, namely, $C(1, 0) := C(u_{0,0}, u_{0,1},$ $ \cdots, u_{0, r_{1}-1})$, $C(1, 1) := C(u_{1, 0}, u_{1, 1}, \cdots, u_{1, r_{1}-1})$, $\cdots$, $C(1, s_{1}-1) := C(u_{s_{1}-1, 0}, u_{s_{1}-1, 1},$ $ \cdots, u_{s_{1}-1, r_{1}-1})$ in order. Similarly, the $T_2$ has $s_2$ number of horizontal cycles of type $A_1$, namely, $C(2, 0) := C(v_{0, 0}, v_{0, 1}, \cdots, v_{0, r_{2}-1})$, $C(2, 1) :=C(v_{1, 0}, v_{1, 1}, \cdots$, $v_{1, r_{2}-1}, v_{1, 0})$, $\cdots$, $C(2, s_{2}-1) := C(v_{s_{2}-1, 0},  v_{s_{2}-1, 1}, \cdots, v_{s_{2}-1, r_{2}-1})$ in order. Now we have the following cases. 

\textbf{Case 1 :} If $(r_{1}, s_1, k_1) = (r_{2}, s_{2}, k_{2})$ then $r_{1} = r_{2}, s_{1} = s_{2}, k_{1} = k_{2}$. We define a map $f_{1} : V(T(r_{1}, s_{1}, k_{1})) \rightarrow V(T(r_{2}, s_{2}, k_{2}))$ such that $f_1(u_{t,i})$ = $v_{t,i}$ for $0 \leq t \leq s-1$ and $0\leq i\leq r-1$. Observe that $lk(u_{t,i}) = C(\textit{\textbf{u}}_{t-1, i-1}, u_{t-1, i}, \textit{\textbf{u}}_{t-1, i+1}, u_{t, i+1}, u_{t+1, i+1},$ $u_{t+1, i}, u_{t, i-1})$ is the link of the vertex $u_{t,i}$ in $T(r_{i}, s_{i}, k_{i})$. By $f_1$, $f_1(lk(u_{t,i})) = C(f_1(\textit{\textbf{u}}_{t-1, i-1}),$ $f_1(u_{t-1, i}), f_1(\textit{\textbf{u}}_{t-1, i+1}),  f_1(u_{t, i+1}),$ $f_1(u_{t+1, i+1}),$ $f_1(u_{t+1, i}),$ $f_1(u_{t, i-1}))  =  C(\textit{\textbf{v}}_{t-1, i-1}, v_{t-1, i},$ $\textit{\textbf{v}}_{t-1, i+1}, v_{t, i+1},$ $v_{t+1, i+1}, v_{t+1, i}, v_{t, i-1})$. So, $f_1(lk(u_{t,i})) = lk(v_{t, i})$. This implies that the map $f_1$ sends vertices to vertices, edges to edges, faces to faces and also, preserves incidents. Therefore, the map $f_{1}$ defines an isomorphism map between $T(r_{1}, s_{1}, k_{1})$ and $T(r_{2}, s_{2}, k_{2})$. Thus, $T_{1} \cong T_{2}$ by $f_1$.       

\textbf{Case 2 :} If $r_{1} \neq r_{2}$ then length$(C_{1, 1}) \neq$ length$(C_{2, 1})$. This implies that $a_{1, 1} \neq a_{2, 1}$, a contradiction since $a_{1, 1} = a_{2, 1}$. So, $r_{1} = r_{2}$.

\textbf{Case 3 :} If $s_{1} \neq s_{2}$ then $n_{1} = r_{1}s_{1} \neq r_{1} s_{2} = n_{2}$ as $r_{1} = r_{2}$ by Case 2. This is a contradiction since $n_1= n_2$. So, $s_{1} = s_{2}$.       

\textbf{Case 4 :} Suppose $k_{1} \neq k_{2}$. By assumption, $a_{1, 4} = a_{2, 4}$, length($C_{1, 4}$) = length($C_{2, 4}$). This implies that min$\{k_{1}+s_{1}, r_{1}+\frac{s_{1}}{2}-k_{1}\}$ = min$\{k_{2}+s_{2}, r_{2}+\frac{s_{2}}{2}-k_{2}\}$. It follows that $k_{1}+s_{1} \neq k_{2}+s_{2}$ since $k_{1} \neq k_{2}$ and $s_{1} = s_{2}$. This gives that $k_{1} + s_{1} = r_{2}+\frac{s_{2}}{2}-k_{2} = r_{1}+\frac{s_{1}}{2}-k_{2}$ as $r_1=r_2$ and $s_1 = s_2$. That is, $k_{2} = r_{1} - k_{1}+\frac{s_{1}}{2}-s_{1} = r_{1} - k_{1}-\frac{s_{1}}{2}$. In this case, identify $T_2$ along vertical identical boundary $P(v_{0,0}, v_{1,0}, \cdots, v_{s_{2}-1, 0}, v_{0, k_{1}})$ of $T_2$ and then cut along the path $Q = P(v_{0,0}, v_{1,0}, v_{2,1}, \cdots, v_{s_{2}-1, \frac{s_{2}}{2}-1}, v_{0, \frac{s_{2}}{2}+k_{2}})$ of type $A_{2}$ through vertex $v_{0, 0}$. We get a new $(r, s, k)$-representation of $M_2$ and we denote it by $R$. This process defines the map $f_{2} : V(T(r_{2}, s_{2}, k_{2}))\rightarrow V(R)$ such that $f_{2}(v_{t, i}) = v_{t,(r_{2}-i+[\frac{t}{2}])(mod~r_{2})}$ for $0\leq t \leq s_{2}-1$ and $0\leq i\leq r_{2}-1$. In $R$, the base horizontal cycle is $C'(2, 0) := C(v_{0, 0}, v_{0, r_{2}-1}, \cdots, v_{0, 1})$, upper horizontal cycle is $C(v_{0,k_{2}+\frac{s_{2}}{2}}, v_{0,k_{2}+\frac{s_{2}}{2}-1}, \cdots, v_{0,k_{2}+\frac{s_{2}}{2}+1})$ and the length of the path $P(v_{0, 0}$,$v_{0, r_{2}-1}, \cdots, v_{0,k_{2}+\frac{s_{2}}{2}})$ in $C'(2, 0)$ is $r_{2}-\frac{s_{2}}{2}-k_{2}$. In this process, we are not changing the both length of the horizontal cycles and number of horizontal cycles which are homologous to the cycle $C'(2, 0)$. So, we get $R = T(r_{2}, s_{2}, r_{2}-k_{2}-\frac{s_{2}}{2})$. Now $r_{2}-k_{2}-\frac{s_{2}}{2} = r_{2} -(r_{1} - k_{1} - \frac{s_{1}}{2})-\frac{s_{2}}{2} = k_{1}$ since $r_{1} = r_{2}, s_{1} = s_{2}$ and $k_2=r_1-k_1-\frac{s_1}{2}$. Thus, by $f_1$, $T(r_{2}, s_{2}, r_{2}-k_{2}-\frac{s_{2}}{2}) \cong T(r_{1}, s_{1}, k_{1})$. So, $T_{1} \cong T_{2}$. So, by Cases 1, 2, 3, 4, claim follows. 

Therefore, by $f_1$, $M_{1} \cong M_{2}$.

Conversely, let $M_{1} \cong M_{2}$. Then, there is an isomorphism map $f : V(M_{1}) \rightarrow V(M_{2})$. Let $C_{1, j}$ be cycle of type $A_{j}$ for $j = 1 , 2, 3 ,4 $ in $M_{1}$. By $f$, consider $C_{2, j} := f(C_{1, j})$ for $j = 1, 2, 3, 4$. So, length$(C_{1, j})$ = length$(f(C_{1, j}))$ = length$(C_{2, j})$ for $j = 1, 2, 3, 4$ since $f$ is an isomorphism map. Hence, $M_1$ and $M_2$ have same cycle-type.
\end{proof}

Thus, we state the following corollary.

\begin{cor} \label{cor3342} 
$T(r_{1}, s_{1}, k_{1}) \not\cong T(r_{2}, s_{2}, k_{2}) ~ \forall ~ r_{1} \neq r_{2}$, $T(r_{1}, s_{1}, k_{1}) \not\cong T(r_{2}, s_{2}, k_{_{2}})$ $~ \forall ~ s_{1} \neq s_{2}$, $T(r_{1}, s_{1}, k_{1}) \not\cong T(r_{1}, s_{1}, k_{2})$ if $s_{1} = 2$ and $k_{2}\in \{2, 3,\cdots,r_{1}-3\} \setminus \{k_{1}, r_{1}-k_{1}-1\} $, $T(r_{1}, s_{1}, k_{1}) \not\cong T(r_{1}, s_{1}, k_{2})$ if $s_{1} \geq  4$ and $k_{2}\in \{0, 1, \cdots,r_{1}-1\} \setminus \{k_{1}, r_{1}-k_{1}-\frac{s_{1}}{2}\}$, $T(r_{1}, s_{1}, k_{1}) \cong T(r_{1}, s_{1}, r_{1}-k_{1}-1)$ if $s_{1} = 2$ and $r_{1} \geq 5$, and $T(r_{1}, s_{1}, k_{1}) \cong T(r_{1}, s_{1}, r_{1}-\frac{s_{1}}{2}-k_{1})$ if $s_{1} \geq 4$ and $r_{1} \geq 3$.
\end{cor}

\begin{proof} If $r_{1} \neq r_{2}$ then it follows that $a_{1, 1} \neq a_{2, 1}$. This implies that $T(r_{1}, s_{1}, k_{1}) \not\cong T(r_{2}, s_{2}, k_{2})$ by Lemma \ref{lem3342-iso}. So, $T(r_{1}, s_{1}, k_{1}) \not\cong T(r_{2}, s_{2}, k_{2}) ~\forall~ r_{1} \neq r_{2}$. Again, $s_{1} \neq s_{2}$ implies $r_{1} \neq r_{2}$ since $r_1s_1=r_2s_2$. This implies that $T(r_{1}, s_{1}, k_{1}) \not\cong T(r_{2}, s_{2}, k_{_{2}}) ~\forall~ s_{1} \neq s_{2}$. If $k_{1} \neq k_{2}$, $r_{1} = r_{2}$ and $s_{1} = s_{2}$ then by the argument in the proof of Lemma \ref{lem3342-iso}, $T(r_{1}, s_{1}, k_{1}) \cong T(r_{1}, s_{1}, k_{2})$ if and only if $k_{2} = r_{1} - k_{1}-\frac{s_{1}}{2}$. So, $T(r_{1}, s_{1}, k_{1}) \cong T(r_{1}, s_{1}, k_{2})$ if $s_1 = 2$ and $k_{2} \neq r_{1} - k_{1}- 1$. Thus, $T(r_{1}, s_{1}, k_{1}) \not\cong T(r_{1}, s_{1}, k_{2})$ if $s_{1} = 2$ and $k_{2}\in \{2,3,\cdots,r_{1}-3\} \setminus \{k_{1}, r_{1}-k_{1}-1\}$, and $T(r_{1}, 2, k_{1}) \cong T(r_{1}, 2, r_{1}-k_{1}-1)$ if $r_{1} \geq 5$ (by Lemma \ref{lem3342-all-rep}). Again, $T(r_{1}, s_{1}, k_{1}) \cong T(r_{1}, s_{1}, k_{2})$ if $s_1 \geq 4$ and $k_{2} \neq r_{1} - \frac{s_1}{2 }- k_{1}$. So, $T(r_{1}, s_{1}, k_{1}) \not\cong T(r_{1}, s_{1}, k_{2})$ if $s_{1} \geq  4$ and $k_{2}\in \{0, 1, \cdots,r_{1}-1\} \setminus \{k_{1}, r_{1}-k_{1}-\frac{s_{1}}{2}\}$, and $T(r_{1}, s_{1}, k_{1}) \cong T(r_{1}, s_{1}, r_{1}-\frac{s_{1}}{2}-k_{1})$ if $s_{1} \geq 4$ and $r_{1} \geq 3$ (by Lemma \ref{lem3342-all-rep}). This completes the proof.
\end{proof}

\begin{table}
\tiny 
\caption{Maps of type $\{3^{3}, 4^{2}\}$}
\centering 
\begin{tabular}{c c c c} 
\hline\hline 
$n$&Equivalence classes& Length of cycles & $i(n)$\\ [0.3ex] 
\hline 
10 & T(5, 2, 2) & (5, $\{10,10\}$, 4) & 1(10)\\
\hline
12& T(6, 2, 2), T(6, 2, 3) & (6, $\{6, 4\}$, 4) & 3(12) \\
 & T(3, 4, 0), T(3, 4, 1) & (3, $\{4, 12\}$, 4)&\\

 & T(3, 4, 2) & (3, $\{12, 12\}$, 6)&\\
\hline
14 & T(7, 2, 2), T(7, 2, 4) & (7, $\{14, 14\}$, 4) & 2(14)\\

 & T(7, 2, 3) & (7, \{14, 14\}, 5)&\\
\hline
16& T(8, 2, 2), T(8, 2, 5) & (8, $\{8, 16\}$, 4) & 5(16) \\

 & T(8, 2, 3), T(8, 2, 4) & (8, $\{16, 4\}$, 5)&\\

 & T(4, 4, 0), T(4, 4, 2) & (4, $\{4, 8\}$, 4)&\\

 & T(4, 4, 1) & (4, $\{16, 16\}$, 5)&\\

 & T(4, 4, 3) & (4, $\{16, 16\}$, 7)&\\
\hline
18& T(9, 2, 2), T(9, 2, 6) & (9, $\{18, 6\}$, 4) & 5(18)\\

 & T(9, 2, 3), T(9, 2, 5) & (9, $\{6, 18\}$, 5)&\\

 & T(9, 2, 4) & (9, $\{18, 18\}$, 6)&\\

 & T(3, 6, 0) & (3, $\{6, 6\}$, 6)&\\

 & T(3, 6, 1), T(3, 6, 2) & (3, $\{18, 18\}$, 7)&\\
\hline
20& T(10, 2, 2), T(10, 2, 7) & (10, $\{10, 20\}$, 4) & 6(20)\\

 & T(10, 2, 3), T(10, 2, 6) & (10, $\{20, 10\}$, 5)&\\

 & T(10, 2, 4), T(10, 2, 5) & (10,$\{10, 4\}$, 6)&\\

 & T(5, 4, 0), T(5, 4, 3) & (5, $\{4, 20\}$, 4)&\\

 & T(5, 4, 1), T(5, 4, 2) & (5, $\{20, 20\}$, 5)&\\

 & T(5, 4, 4) & (5, $\{20, 20\}$, 8)&\\
\hline
22& T(11, 2, 2), T(11, 2, 8) & (11, $\{22, 22\}$, 4) & 4(22)\\

 & T(11, 2, 3), T(11, 2, 7) & (11, $\{22, 22\}$, 5)&\\

 & T(11, 2, 4), T(11, 2, 6) & (11, $\{22, 22\}$, 6)&\\

 & T(11, 2, 5) & (11, $\{22, 22\}$, 7)&\\
[1ex] 
\hline\hline 
\end{tabular}
\label{table1} 
\end{table}

We calculate all possible $T(r, s, k)$ representations on $n$ vertices by Lemma \ref{lem3342-all-rep}. Then, we calculate lengths of the cycles of type $A_{i}$ for $i \in \{1, 2, 3, 4\}$. Next, we classify all $T(r, s, k)$ representation by Lemma \ref{lem3342-iso} up to isomorphism. So, by the Lemmas \ref{lem3342-all-rep}, \ref{lem3342-iso}, maps of type $\{3^{3}, 4^{2}\}$ can be classified up to isomorphism. We repeat this same argument in the Sections \ref{324341}, \ref{36361}, \ref{31221}, \ref{3461}, \ref{46121}, \ref{34641}, \ref{4881}. We have done the above calculations for the vertices $n \le 22$. We have listed the obtained objects in the form of their $(r, s, k)$-representation in Table \ref{table1}. In Table \ref{table1}, we have used $n$ to denote the number of vertices of a map. We put $T(r_{1}, s_{1}, k_{1})$ and $T(r_{2}, s_{2}, k_{2})$ in a single equivalence class if $T(r_{1}, s_{1}, k_{1})$ and $T(r_{2}, s_{2}, k_{2})$ are isomorphic. We have used $(a_{1}, \{a_{2}, a_{3}\}, a_{4})$ to denote a permutation of lengths of cycles where $a_{j}$ $=$ length$(C_{1, j})$ for $j \in \{1, 2, 3, 4\}$ and $\{a_{2}, a_{3}\}$ denotes a set of lengths of the cycles $C_{1, 2}$ and $C_{1, 3}$ of type $A_{2}$. We have also used $i(n)$ where $i$ denote the number of non-isomorphic objects of type $\{3^{3}, 4^{2}\}$ on $n$ vertices up to isomorphism. The above notations are also used in Tables \ref{table2}, \ref{table3}, \ref{table4}, \ref{table5}, \ref{table6}, \ref{table7}, \ref{table8}.

\section{Maps of type $\{3^{2}, 4, 3, 4\}$}\label{324341}

Let $M$ be a map of type $\{3^{2}, 4, 3, 4\}$ on the torus. Through each vertex in $M$ there is a path as follows.

\begin{defn}\label{defn32434} Let $P(\cdots, u_{i-1}, u_{i}, u_{i+1}, \cdots)$ be a path in edge graph of $M$. We say the path $P$ of type $B_1$ if $lk(u_{i})=C(\textbf{a},$ $u_{i+1}, b,$ $c, \textbf{d}, u_{i-1}, e)$ implies $lk(u_{i-1})=C(\textbf{f}, g,$ $e, u_{i}, \textbf{c},$ $d, u_{i-2})$ and $lk(u_{i+1})=C(\textbf{e},$ $a, k, u_{i+2}, \textbf{l},$ $b, u_{i})$, and  $lk(u_{i})=C(\textbf{e},$ $h, k, u_{i+1}, \textbf{l},$ $b, u_{i-1})$ implies $lk(u_{i-1})=C(\textbf{h}, u_{i},$ $b, c, \textbf{d},$ $u_{i-2}, e)$ and $lk(u_{i+1})=C(\textbf{s},$ $u_{i+2}, t, l, $\textbf{b},$ u_{i}, k)$.
\end{defn}

In Figure 17, lk$(u_i) = C(\textit{\textbf{a}}, b, c, u_{i+1}, \textit{\textbf{f}}, e, u_{i-1})$ and path $P(u_{i-1}, u_i, u_{i+1})$ is part of a path of type $B_1$. Let $P$ be a maximal path of type $B_{1}$. Then, by the next Lemma \ref{lem32434-cycle}, it defines a cycle.

\vspace{-1.2cm}

\begin{picture}(0,0)(-80,40)
\tiny
\setlength{\unitlength}{1.5mm}
\drawpolygon(5,5)(20,5)(20,10)(5,10)
\drawline[AHnb=0](12.5,5)(12.5,10)
\drawpolygon(12.5,10)(20,10)(20,15)(12.5,15)
\drawline[AHnb=0](12.5,10)(20,5)
\drawline[AHnb=0](5,10)(12.5,15)

\drawline[linecolor=black,linewidth=0.2,AHnb=0](5,10)(20,10)
\drawline[linecolor=black,linewidth=0.22,AHnb=0](12.5,5)(12.5,15)

\drawline[AHnb=0,dash={1.0 1.0 1.0 1.0}{0.0}](20,10)(24,10)
\drawline[AHnb=0,dash={1.0 1.0 1.0 1.0}{0.0}](5,10)(1,10)

\drawline[AHnb=0,dash={1.0 1.0 1.0 1.0}{0.0}](12.5,5)(12.5,2.5)
\drawline[AHnb=0,dash={1.0 1.0 1.0 1.0}{0.0}](12.5,15)(12.5,18)

\put(2,11){\tiny ${B_1}$}
\put(12.5,19){\tiny ${B_1}$}

\put(5.5,8.8){\tiny ${u_{i-1}}$}
\put(13,11){\tiny ${u_{i}}$}
\put(16.3,8.8){\tiny ${u_{i+1}}$}

\put(3.5,5){\tiny ${a}$}
\put(13,3.5){\tiny ${b}$}
\put(20.5,5){\tiny ${c}$}

\put(13,15.5){\tiny ${e}$}
\put(20,15.5){\tiny ${f}$}

\put(6,1){\tiny  Figure 17 : $lk(u_i)$}


\end{picture}

\begin{picture}(0,0)(3,30)
\tiny
\setlength{\unitlength}{1.5mm}
\drawpolygon(5,5)(20,5)(20,15)(5,15)
\drawpolygon(25,5)(40,5)(40,15)(25,15)


\drawline[AHnb=0](10,5)(10,15)
\drawline[AHnb=0](15,5)(15,15)
\drawline[AHnb=0](30,5)(30,15)
\drawline[AHnb=0](35,5)(35,15)
\drawline[AHnb=0](5,10)(20,10)
\drawline[AHnb=0](25,10)(40,10)

\drawline[AHnb=0](5,10)(10,15)
\drawline[AHnb=0](10,10)(15,5)
\drawline[AHnb=0](15,10)(20,15)
\drawline[AHnb=0](25,10)(30,5)
\drawline[AHnb=0](30,10)(35,15)
\drawline[AHnb=0](35,10)(40,5)

\put(21,10){$\cdots$}
\put(21,5){$\cdots$}
\put(21,15){$\cdots$}

\put(5,3.8){\tiny ${v_{1}}$}
\put(10,3.8){\tiny ${v_{2}}$}
\put(15,3.8){\tiny ${v_{3}}$}
\put(20,3.8){\tiny ${v_{4}}$}
\put(25,3.8){\tiny ${v_{r-2}}$}
\put(30,3.8){\tiny ${v_{r-1}}$}
\put(35,3.8){\tiny ${v_{r}}$}
\put(40,3.8){\tiny ${v_{1}}$}

\put(5.5,9){\tiny ${w_{1}}$}
\put(10.5,10.5){\tiny ${w_{2}}$}
\put(15.5,9){\tiny ${w_{3}}$}
\put(20.5,9){\tiny ${w_{4}}$}
\put(25.5,10.5){\tiny ${w_{r-2}}$}
\put(30.5,9){\tiny ${w_{r-1}}$}
\put(35.5,10.5){\tiny ${w_{r}}$}
\put(40.5,9){\tiny ${w_{1}}$}

\put(5.5,13.8){\tiny ${x_{1}}$}
\put(10.5,13.8){\tiny ${x_{2}}$}
\put(15.5,13.8){\tiny ${x_{3}}$}
\put(20.5,13.8){\tiny ${x_{4}}$}
\put(25.5,13.8){\tiny ${x_{r-2}}$}
\put(30.5,13.8){\tiny ${x_{r-1}}$}
\put(35.5,13.8){\tiny ${x_{r}}$}
\put(40.5,13.8){\tiny ${x_{1}}$}

\put(16,1){\tiny  Figure 18 : Cylinder}

\end{picture}

\begin{picture}(0,0)(12,95)
\tiny
\setlength{\unitlength}{1.5mm}

\drawline[AHnb=0](9,20)(46,20)
\drawline[AHnb=0](9,25)(46,25)
\drawline[AHnb=0](9,30)(46,30)
\drawline[AHnb=0](9,35)(46,35)
\drawline[AHnb=0](9,40)(46,40)
\drawline[AHnb=0](9,45)(46,45)

\drawline[AHnb=0](10,19)(10,46)
\drawline[AHnb=0](15,19)(15,46)
\drawline[AHnb=0](20,19)(20,46)
\drawline[AHnb=0](25,19)(25,46)
\drawline[AHnb=0](30,19)(30,46)
\drawline[AHnb=0](35,19)(35,46)
\drawline[AHnb=0](40,19)(40,46)
\drawline[AHnb=0](45,19)(45,46)

\drawline[AHnb=0](10,20)(15,25)
\drawline[AHnb=0](20,20)(25,25)
\drawline[AHnb=0](30,20)(35,25)
\drawline[AHnb=0](40,20)(45,25)

\drawline[AHnb=0](20,25)(15,30)
\drawline[AHnb=0](30,25)(25,30)
\drawline[AHnb=0](40,25)(35,30)

\drawline[AHnb=0](10,30)(15,35)
\drawline[AHnb=0](20,30)(25,35)
\drawline[AHnb=0](30,30)(35,35)
\drawline[AHnb=0](40,30)(45,35)

\drawline[AHnb=0](20,35)(15,40)
\drawline[AHnb=0](30,35)(25,40)
\drawline[AHnb=0](40,35)(35,40)

\drawline[AHnb=0](10,40)(15,45)
\drawline[AHnb=0](20,40)(25,45)
\drawline[AHnb=0](30,40)(35,45)
\drawline[AHnb=0](40,40)(45,45)

\drawline[AHnb=0](10,45)(9,46)\drawline[AHnb=0](20,45)(19,46)
\drawline[AHnb=0](30,45)(29,46)\drawline[AHnb=0](40,45)(39,46)
\drawline[AHnb=0](10,35)(9,36)\drawline[AHnb=0](10,25)(9,26)

\drawline[AHnb=0](15,20)(16,19)\drawline[AHnb=0](25,20)(26,19)
\drawline[AHnb=0](35,20)(36,19)\drawline[AHnb=0](45,20)(46,19)
\drawline[AHnb=0](45,30)(46,29)\drawline[AHnb=0](45,40)(46,39)

\put(12,24){\tiny ${u_{i}}$}
\put(16.2,24){\tiny ${u_{i+1}}$}
\put(20.2,24){\tiny ${u_{i+2}}$}
\put(26.2,24){\tiny ${u_{i+3}}$}
\put(30.2,24){\tiny ${u_{i+4}}$}
\put(36.2,24){\tiny ${u_{i+5}}$}
\put(41,24){\tiny ${u_{i+6}}$}

\put(13.2,29){\tiny ${u_{r}}$}
\put(15.5,34){\tiny ${u_{r-1}}$}
\put(10.8,39){\tiny ${u_{r-2}}$}
\put(15.5,44){\tiny ${u_{r-3}}$}

\drawline[linecolor=black,linewidth=.2,AHnb=0](15,25)(46,25)
\drawline[linecolor=black,linewidth=.2,AHnb=0](15,25)(15,46)

\drawline[linecolor=black,linewidth=.2,AHnb=0](20,30)(46,30)
\drawline[linecolor=black,linewidth=.2,AHnb=0](20,30)(20,46)

\put(18,29){\tiny ${w_{i}}$}
\put(21,29){\tiny ${w_{i+1}}$}
\put(26,29){\tiny ${w_{i+2}}$}
\put(31,29){\tiny ${w_{i+3}}$}
\put(36,29){\tiny ${w_{i+4}}$}
\put(40.8,29){\tiny ${w_{i+5}}$}

\put(20.2,34){\tiny ${w_{r-2}}$}
\put(20.2,39){\tiny ${w_{r-3}}$}
\put(20.2,44){\tiny ${w_{r-4}}$}

\put(24,16){\tiny   Figure 20}

\end{picture}

\vspace{4cm}

\begin{picture}(0,0)(-55,25)
\tiny
\setlength{\unitlength}{1.5mm}
\drawpolygon(5,5)(20,5)(20,25)(5,25)
\drawpolygon(25,5)(40,5)(40,25)(25,25)
\drawpolygon(45,5)(55,5)(55,25)(45,25)


\drawline[AHnb=0](10,5)(10,25)
\drawline[AHnb=0](15,5)(15,25)
\drawline[AHnb=0](20,5)(20,25)
\drawline[AHnb=0](25,5)(25,25)
\drawline[AHnb=0](30,5)(30,25)
\drawline[AHnb=0](35,5)(35,25)
\drawline[AHnb=0](50,5)(50,25)
\drawline[AHnb=0](55,5)(55,25)

\drawline[AHnb=0](5,10)(20,10)
\drawline[AHnb=0](25,10)(40,10)
\drawline[AHnb=0](5,15)(20,15)
\drawline[AHnb=0](25,10)(40,10)
\drawline[AHnb=0](5,20)(20,20)
\drawline[AHnb=0](5,10)(20,10)
\drawline[AHnb=0](45,10)(55,10)
\drawline[AHnb=0](45,15)(55,15)
\drawline[AHnb=0](45,20)(55,20)

\drawline[AHnb=0](5,20)(10,25)
\drawline[AHnb=0](10,20)(15,15)
\drawline[AHnb=0](15,20)(20,25)
\drawline[AHnb=0](25,20)(30,15)
\drawline[AHnb=0](30,20)(35,25)
\drawline[AHnb=0](35,20)(40,15)
\drawline[AHnb=0](25,20)(40,20)
\drawline[AHnb=0](50,10)(55,5)
\drawline[AHnb=0](45,10)(50,15)
\drawline[AHnb=0](50,20)(55,15)
\drawline[AHnb=0](45,20)(50,25)

\drawline[AHnb=0](5,10)(10,15)
\drawline[AHnb=0](10,10)(15,5)
\drawline[AHnb=0](15,10)(20,15)
\drawline[AHnb=0](25,10)(30,5)
\drawline[AHnb=0](30,10)(35,15)
\drawline[AHnb=0](35,10)(40,5)
\drawline[AHnb=0](25,15)(40,15)

\put(21,10){$\cdots$}
\put(21,5){$\cdots$}
\put(21,15){$\cdots$}
\put(21,20){$\cdots$}
\put(21,25){$\cdots$}

\put(42,10){$\cdots$}
\put(42,5){$\cdots$}
\put(42,15){$\cdots$}
\put(42,20){$\cdots$}
\put(42,25){$\cdots$}

\put(5,4){\tiny ${v_{1}}$}
\put(10,4){\tiny ${v_{2}}$}
\put(15,4){\tiny ${v_{3}}$}
\put(20,4){\tiny ${v_{4}}$}
\put(25,4){\tiny ${v_{2k}}$}
\put(30,4){\tiny ${v_{2k+1}}$}
\put(35,4){\tiny ${v_{2k+2}}{}$}
\put(40,4){\tiny ${v_{2k+3}}$}
\put(45,4){\tiny ${v_{r-1}}$}
\put(50,4){\tiny ${v_r}{}$}
\put(55,4){\tiny ${v_{1}}$}

\put(5.5,9){\tiny ${w_{1}}$}
\put(10.5,10.5){\tiny ${w_{2}}$}
\put(15.5,9){\tiny ${w_{3}}$}
\put(20.5,9){\tiny ${w_{4}}$}
\put(25.5,10.5){\tiny ${w_{2k}}$}
\put(30.5,9){\tiny ${w_{2k+1}}$}
\put(35.5,10.5){\tiny ${w_{2k+2}}$}
\put(40.5,9){\tiny ${w_{2k+3}}$}
\put(45.5,9){\tiny ${w_{r-1}}$}
\put(50.5,10.5){\tiny ${w_{r}}$}
\put(55.5,9){\tiny ${w_{1}}$}

\put(5.5,14){\tiny ${x_{1}}$}
\put(10.5,14){\tiny ${x_{2}}$}
\put(15.5,14){\tiny ${x_{3}}$}
\put(20.5,14){\tiny ${x_{4}}$}
\put(25.5,14){\tiny ${x_{2k}}$}
\put(30.5,14){\tiny ${x_{2k+1}}$}
\put(35.5,14){\tiny ${x_{2k+2}}$}
\put(40.5,14){\tiny ${x_{2k+3}}$}
\put(45.5,14){\tiny ${x_{r-1}}$}
\put(50.5,14){\tiny ${x_{r}}$}
\put(55.5,14){\tiny ${x_{1}}$}

\put(5.5,19){\tiny ${z_{1}}$}
\put(10.5,20.5){\tiny ${z_{2}}$}
\put(15.5,19){\tiny ${z_{3}}$}
\put(20.5,19){\tiny ${z_{4}}$}
\put(25.5,20.5){\tiny ${z_{2k}}$}
\put(30.5,19){\tiny ${z_{2k+1}}$}
\put(35.5,20.5){\tiny ${z_{2k+2}}$}
\put(40.5,19){\tiny ${z_{2k+3}}$}
\put(45.5,19){\tiny ${z_{r-1}}$}
\put(50.5,20.5){\tiny ${z_{r}}$}
\put(55.5,19){\tiny ${z_{1}}$}

\put(5.2,24){\tiny ${v_{2k+1}}$}
\put(10.5,24){\tiny ${v_{2k+2}}$}
\put(15.2,24){\tiny ${v_{2k+3}}$}
\put(20.5,24){\tiny ${v_{2k+4}}$}
\put(25.5,24){\tiny ${v_{r}}$}
\put(30.2,24){\tiny ${v_{1}}$}
\put(35.5,24){\tiny ${v_{2}}$}
\put(40.5,24){\tiny ${v_{3}}$}
\put(45.2,24){\tiny ${v_{2k-1}}$}
\put(50.5,24){\tiny ${v_{2k}}$}
\put(55.5,24){\tiny ${v_{2k+1}}$}

\put(21,0) {\tiny  Figure 19 : $T(r, 4, 2k)$}

\end{picture}

\vspace{2.5 cm}

\begin{lemma}\label{lem32434-cycle}  If $P$ is a maximal path of type $B_{1}$ in $M$ then there exists an edge $e$ such that $P \cup e$ is a cycle.
\end{lemma}

\begin{proof} Let $P(u_{1}, u_{2}, \cdots, u_{r})$ be a maximal path of type $B_{1}$ and $lk(u_{r}) = C( u_{r-1}, \textit{\textbf{a}}, b, c, d, \textit{\textbf{f}},$ $e)$. If $d = u_{1}$ then $C(u_{1}, u_{2}, \cdots, u_{r})$ is a cycle. Suppose $d \neq u_{1}$ and $d  = u_i$ for some $2 \le i \le r$. Then, it defines a cycle $L = C(u_i, u_{i+1}, \cdots, u_r)$. By the similar argument as in Lemma \ref{lem3342-cycle} and by Definition \ref{defn32434}, either $f = u_{i+1}$, $d = u_i$, $c = u_{i-1}$ or $c = u_{i+1}$, $d = u_i$, $f = u_{i-1}$. In both the cases, by considering faces incident with the cycle, we get a new cycle $C(w_i, w_{i+1}, \cdots, w_{r-2})$ (see Figure 20) of same type as $L$ with lesser length. By induction, it is impossible similarly as in Lemma \ref{lem3342-cycle}. Therefore, $d \ne u_i$ for $2 \le i \le r$. So, we get a path $Q$ which is extended from $P$ with length($P$) $<$ length($Q$). This is a contradiction as $P$ is maximal. Therefore, $d = u_{1}$ and the path $P$ defines the cycle $C(u_{1},$ $u_{2},$ $\cdots,$ $u_{r})$. So, every maximal path of type $B_1$ is a cycle.
\end{proof}

In Figure 19, the path $P(v_{1}, v_{2}, \cdots, v_{r})$ is of type $B_1$ and the cycle $C(u_{1}, u_{2},\cdots, u_{r})$ is of type $B_1$. Let $C_{1}$ and $C_{2}$ be two cycles of type $B_{1}$. We claim that 

\begin{lemma}\label{lem32434-1} If $C_{1}$ and $C_{2}$ are two cycles of type $B_{1}$ and $E(C_{1}) \cap E(C_{2})\not= \emptyset$ then $C_{1} = C_{2}$.
\end{lemma}

\begin{proof} Let $C_{1}:=C(u_{1, 1}, u_{1, 2}, \cdots, u_{1, r})$ and $C_{2}:=C(u_{2, 1}, u_{2, 2}, \cdots,  u_{2, s})$ and $E(C_{1}) \cap E(C_{2}) \not= \emptyset$. Then, there is an edge $e \in E(C_{1} \cap C_{2})$. Let $e = yx$. The cycles $C_{1},$ $C_{2}$ are both well defined at the vertices $y$ and $x$. Let $lk(x) = C(\textit{\textbf{a}},$ $b,$ $c,$ $w,$ $\textit{\textbf{d}},$ $e,$ $y)$. By Definition \ref{defn32434}, $w \in V(C_{1}\cap C_2)$. So, the path $P(y, x, w)$ is part of both $C_{1}$ and $C_{2}$. This implies that $y = u_{1, t_{1}-1} = u_{2, t_{2}-1}$, $x = u_{1, t_{1}} = u_{2, t_{2}}$ and $w = u_{1, t_{1}+1} = u_{2, t_{2}+1}$ for some $t_{1}\in \{1,$ $\cdots$, $r\}$ and $t_{2} \in \{1,$ $\cdots$, $s\}$. Again, we argue similarly for the edge $xw$ as we did for the edge $e$. We continue with above process. This process stops after, say $r$ number of steps. Let $t_2 > t_1$ and $t_{2}-t_{1} = m$ for some $m$. Then, by this process, we get $u_{1, 1} = u_{2, m+1},$ $u_{1, 2} = u_{2, m+2},\cdots, u_{1, r} = u_{2, m+r}$ and $u_{1, 1} = u_{2, m+r+1}$. This implies that $m+1 = m+r+1$ and $r = s$ as $u_{1, m+r+1} = u_{2, m+1}$ and $C_{2}$ is a cycle. Hence, $C_{1}$ = $C_{2}$. Again, let $lk(x) = C(\textit{\textbf{a}},b, c, w, \textit{\textbf{d}}, y, e)$. Then, by Definition \ref{defn32434}, $ b \in V(C_{1} \cap C_2)$. Similarly again we repeat above argument and we get $C_{1}$ = $C_{2}$. Therefore, by combining above two cases, $E(C_{1}) \cap E(C_{2})\not= \emptyset$ implies $C_{1} = C_{2}$. This completes the proof.
\end{proof}

Let $C$ be a cycle of type $B_1$. Similarly we argue as in Lemma \ref{lem3342-non-con} for the cycles of type $B_1$ and so, the cycle $C$ is non-contractible. Let $S := \{F \in F(M) ~|~ V(C) \cap V(F) \neq \emptyset\}$. The cylinder $S_C = |S|$ has two boundary cycles which are either disjoint or identical by Lemma \ref{lem32434-2}. 

\begin{lemma}\label{lem32434-2} Let $C$ be a cycle of type $B_{1}$, $\partial S_{C} = \{C_{1}, C_{2}\}$. If $C_{1} \cap C_{2}\not= \emptyset$ then $C_{1} = C_{2}$.
\end{lemma}

\begin{proof} The cycle $C$ is of type $B_{1}$ and $\partial S_{C} = \{C_{1}, C_{2}\}$. We argue similarly as in Lemma \ref{lem3342-cycle} for the cycle $C$ and the cylinder $S_C$. So, the cycles $C,$ $C_1$ and $C_2$ are of same type $B_1$. Let $C_{1} \cap C_{2} \neq \emptyset$ and $u \in V(C_{1} \cap C_{2})$. Suppose $C_{1} \cap C_2$ does not contain any edge which is incident at $u$. By Definition \ref{defn32434}, the number of incident edges that lie on one side of the cycle $C_{i}$ is two and on the other side is one at each vertex of $C_i$. Let $d_i$ denote the number of incident edges which are incident at $u$ and does not belong to $E(C_i)$. Hence, $d_1 + d_2 = 3$. The vertex $u \in V(C_1)$ and $u \in V(C_2)$. Since $C_{1} \cap C_2$ does not contain any edge at $u$, so, the cycles $C_1$ and $C_2$ both contain two different edges which are incident at the vertex $u$. This implies that degree$(u) \geq (d_1 + d_2 + 4) = 7$. This is a contradiction as the degree of the vertex $u$ is five. Therefore, $C_{1} \cap C_{2}$ contains an edge at the vertex $u$. This implies that $C_{1} = C_{2}$ by Lemma \ref{lem32434-1}. Again, if $C_{1} \cap C_{2}$ contains an edge then by Lemma \ref{lem32434-1}, $C_{1} = C_{2}$. Therefore, boundary cycles of a cylinder are either identical or disjoint.
\end{proof}

We show that the cycles of type $B_{1}$ have at most two different lengths in the next Lemma \ref{lem32434-hom-length}.

\begin{lemma}\label{lem32434-hom-length} In $M$, the cycles of type $B_{1}$ have at most two different lengths.
\end{lemma}

\begin{proof} We proceed with as in the case of Lemma \ref{lem32434-2}. There are two cycles of type $B_1$ through each vertex of $M$ (by the definition of cycle of type $B_1$). Let $u \in V(M)$. Let $C_1$ and $C_1'$ denote two cycles through a vertex $u$. Consider cylinder $S_{C_1}$ which is defined by cycle $C_1$. Let $\partial S_{C_1} = \{C_2, C_0\}$. The cycles $C_1,$ $C_2,$ and $C_0$ are homologous to each other and length$(C_1)$ = length$(C_2)$ = length$(C_0)$ by the similar argument of Lemma \ref{lem3342-all-length}. Again, we proceed with above argument for the cycle $C_2$ in place of $C_1$ and continue. In this process, let $C_i$ denote a cycle at $i^{th}$ step where $\partial S_{C_i} = \{C_{i+1}, C_{i-1}\}$ and length$(C_{i-1})$ = length$(C_i)$ = length$(C_{i+1})$. Let after $k+1$ number of steps the process stops when the cycle $C_1$ appears. Thus, the cycles $C_i$, $C_j$ are homologous for every $1 \le i, j \le k$ where $\cup_{i=1}^{k}V(C_{i}) = V(M)$ and $l_1=$length$(C_1)$ = length$(C_i)$ for all $1 \leq i \leq k$. Again, we proceed with above process for $C_1'$ in place of $C_1$. Similarly, we get a sequence of homologous cycles, namely, $C_1'$, $C_2'$, $\cdots$, $C_{k_1}'$ such that $\cup_{i=1}^{k_1}V(C_{i}') = V(M)$ and $l_2=$length$(C_i')$ for all $1 \leq i \leq k_1$. So, $M$ contains cycles of type $B_{1}$ at most two different lengths $l_1$ and $l_2$. This completes the proof of lemma.
\end{proof}

Similarly as in Section \ref{33421}, observe that every map of type $\{3^{2}, 4, 3, 4\}$ on the torus has a $T(r, s, k)$ representation for some $r, s, k$. We define admissible relations among $r$, $s$, $k$ of $T(r, s, k)$ such that $T(r, s, k)$ represent a map after identifying their boundaries in next lemma. We omit the proof of next lemma as its argument similarly as in Lemma \ref{lem3342-all-rep}.


\begin{lemma}\label{lem32434-rep} The maps of type $\{3^{2}, 4, 3, 4\}$ of the form $T(r, s, k)$ exist if and only if the following holds : (i) $s \geq  2$ even, (ii) $2 \mid  r$, (iii) $rs \geq  16$, (iv) $k \in \{2t+4 ~:~ 0 \leq t \leq \frac{r-8}{2}\}$ if $s = 2 $ $\&$ $ k \in \{2t ~:~ 0 \leq t < \frac{r}{2}\}$ if $s \ge 4$.
\end{lemma}



Let $T_i = T(r_i, s_i, k_i)$, $i = \{1, 2\}$ denote $M_i$ of type $\{3^{2}, 4, 3, 4\}$ on the torus with $n_i$ vertices and $n_1 = n_2$. Let $C_{i, 1}$ and $C_{i, 2}$ be two non-homologous cycles of type $B_{1}$ in $M_i$ for $i = 1, 2$ and
$a_{i, j}$ = length$(C_{i, j})$. Then, 

\begin{lemma}\label{lem32434-iso}
The map $M_{1} \cong M_{2}$ if and only if $(a_{1, 1}, a_{1, 2})$ = $(a_{2, t_1}, a_{2, t_2})$ for $t_1 \neq t_2 \in \{1, 2\}$.
\end{lemma}

\begin{proof} We first assume that $(a_{1, 1}, a_{1, 2}) = (a_{2, t_1}, a_{2, t_2})$ where $t_1, t_2 \in \{1, 2\}$ and $t_1 \neq t_2$. This implies that $\{a_{1, 1}, a_{1, 2}\} = \{a_{2, 1}, a_{2, 2}\}$.

\smallskip

\noindent $Claim.$  $T_{1} \cong T_{2}$. 

\smallskip

From the definition of $(r_i, s_i, k_i)$-representation, the $T_i$ has $s_i$ number of horizontal cycles of type $B_1$, namely, $C(1, 0) := C(u_{0,0}, u_{0,1}, \cdots, u_{0, r_{1}-1})$, $C(1, 1) :=C(u_{1, 0},$ $u_{1, 1}, \cdots,$ $u_{1, r_{1}-1})$, $\cdots$, $C(1, s_{1}-1) :=C(u_{s_{1}-1, 0},$ $u_{s_{1}-1, 1}, \cdots,$ $u_{s_{1}-1, r_{1}-1})$ in $T_1$ and $C(2, 0) := C(v_{0, 0},$ $v_{0, 1}, \cdots,$ $ v_{0, r_{2}-1}), C(2, 1) := C(v_{1, 0},$ $v_{1, 1}, \cdots,$ $v_{1, r_{2}-1},$ $v_{1, 0})$, $\cdots$, $C(2, s_{2}-1) := C(v_{s_{2}-1, 0},$ $v_{s_{2}-1, 1},\cdots,$ $v_{s_{2}-1, r_{2}-1})$ in $T_2$. 

\textbf{Case 1 :} If $(r_{1}, s_{1}, k_1) = (r_{2}, s_{2}, k_{2})$ i.e. $r_{1} = r_{2}, s_{1} = s_{2}, k_{1} = k_{2}$ then similarly as in the proof of the Lemma \ref{lem3342-iso}, we define an isomorphism map $f_{1} : V(T(r_{1}, s_{1}, k_{1})) \rightarrow V(T(r_{2}, s_{2}, k_{2}))$ such that $f(u_{t,i}) = v_{t,i}$ for $0\leq t\leq s_1-1$ and $0\leq i\leq r_1-1$. So, $T_{1} \cong T_{2}$ by $f_{1}$.     

\textbf{Case 2 :} Suppose $r_{1} = r_{2}, s_{1} = s_{2}, k_{1} \neq k_{2}$. Since $r_1 = r_2$, it implies that the vertical cycles in $T_1$ and $T_2$ have same length. 

We define a cycle in $T_1$ as in equation (1) of Section \ref{33421}. Let $C_{lh}$ denote the base horizontal cycle and $C_{uh}$ denote the upper horizontal cycle in $T_2$. Let $Q$ be a path through $u_{k_1+1}$ of type $B_{1}$ and not homologous to $C_{lh}$. Again, let $Q'$ and $Q''$ denote two edge disjoint paths in $C_{uh}$ such that $C_{uh} = Q' \cup Q''$. Hence as in equation (1), we define a new cycle $C_3(1)$ using the above paths in $T_1$. Similarly, there is a cycle $C_3(2)$ as $C_3(1)$ in $T_2$. Since $\{a_{1, 1}, a_{1, 2}\} = \{a_{2, 1}, a_{2, 2}\}$ and $r_1 = r_2$, it follows that length($C_3(1)$) = length($C_3(2)$).  

Since length($C_3(1)$) = length($C_3(2)$), it implies that min$\{s_{1}+k_{1}, r_{1}+s_{1}-k_{1}\}$ = min$\{s_{2}+k_{2}, r_{2}+s_{2}-k_{2}\}$. It follows that $r_{1}+s_{1}-k_{1} = s_{2}+k_{2}$ since $k_{1} \neq k_{2}$. So, $k_{2} = r_{1}-k_{1}$ as $s_1 = s_2$. We proceed similarly as in the Lemma \ref{lem3342-iso}. In this process, identify $T(r_{2}, s_{2}, k_{2})$ along vertical boundary and cut along a path $Q := P(v_{0,i}, v_{1,i}, \cdots, v_{s_{2}-1, i}, v_{0, i+k_{2}})$ for some odd  $0 \le i \le r_1-1$. Thus, we get a new representation, say $R$ of $M_2$ with a map $f_{2} : V(T(r_{2}, s_{2}, k_{2}))\rightarrow V(R)$ such that $f_{2}(v_{t, i'}) = v_{t,(i+r_{2}-i')(mod~r_{2})}$ for $0\leq t\leq s_{2}-1$ and $0\leq i'\leq r_{2}-1$. Clearly, $f_2$ maps cycle $C(2, t) :=C(v_{t, 0}, v_{t, 1}, \cdots, v_{t, r_{2}-1})$ to cycle $C'(2, t) :=C(v_{t, i}, v_{t, i-1}, \cdots, v_{t, r_{2}-1}, v_{t, 0}, v_{t, 1}, \cdots, v_{t, i+1})$. Since the path $Q_1 := P(v_{0, i}, v_{0, i-1}, \cdots, v_{0, i+k_2}) \subset C'(2, 0) := C(v_{0, i}, v_{0, i-1}, \cdots, v_{0, i+1})$ and length$(Q_1) = i+r_{2}-k_{2}-i = r_{2}-k_{2}$, it follows that $R$ has $s_{2}$ number of horizontal cycles of length $r_{2}$ and the cycle of type $B_{2}$ has length $r_{2} - k_{2}$ as length$(Q_1) = r_2 - k_2$. In $R$, the faces incident with the base horizontal cycle $C' = C(v_{0, i}, v_{0, i-1}, \cdots, v_{0, r_{2}-1}, v_{0, 0}, v_{0, 1}, \cdots, v_{0, i+1})$ at  $v_{0, i-1}$ have two possibilities as follows. If one triangle and one $4$-gon incident at $v_{0, i-1}$ in $C'$ then $R$ has a $(r, s, k)$-representation. If two triangles and one $4$-gon incident at $v_{0, i-1}$ then $R$ does not follow the definition of $(r, s, k)$-representation. (For an example, the Figure 12 is an example of $R$ which does not follow the definition of $(r, s, k)$-representation since the number of incident triangles at $w_2$ is two.) In this case, if $C'(2, 0), C'(2, 1),$ $\cdots$, $C'(2, s_2-1)$ denote a sequence of horizontal cycles in $R$ then we identify $R$ along $C(v_{0, i}, v_{0, i-1}, \cdots, v_{0, i+1})$ and cut along $C(v_{1, i}, v_{1, i-1}, \cdots, v_{1, i+1})$. Thus, we get a new representation of $M_2$, say $R'$ where $C'(2, 1) := C(v_{1, i}, v_{1, i-1}, \cdots, v_{1, i+1})$ denote the base horizontal cycle. In this process, $C'(2, 1)\rightarrow C'(2, 0), C'(2, 0)\rightarrow C'(2, 1), C'(2, s_2-1) \rightarrow C'(2, 2),$ $\cdots$, $C'(2, 2) \rightarrow C'(2, s_2-1)$. This process defines a map $f_3 : R \rightarrow R' $ such that $f_3(C'(2, t)) = C'(2, 1-t (mod~s_{2}))$ for $0 \le t \le s_2-1$. In $R'$, $C'(2, 1), C'(2, 0), C'(2, s_2-1),$ $\cdots$, $C'(2, 2)$ denote the sequence of horizontal cycles in $R'$. (For an example of $R'$, the Figure 5 is a $T(6, 4, 2)$ representation which is defined from $R$ in Figure 12. In Figure 12, we cut $R$ along the cycle $C(v_{1}, v_{2}, \cdots, v_{6})$ and identify along $C(w_{1}, w_{2}, \cdots, w_{6})$. Hence, we get a representation $T(6, 4, 2)$ in Figure 5.) In the above process, we are redefining $R$ to a desire representation $R'$. In this process, the length of the horizontal cycles of type $B_{1}$ are remain unchanged as we are only changing the order of the horizontal cycles. So, $R'$ has a well defined $T(r_{2}, s_{2}, r_{2}-k_{2})$ representation. Thus, $T(r_{2}, s_{2}, r_{2}-k_{2}) = T(r_{1}, s_{1}, k_{1})$ since $r_{1} = r_{2}, s_{1} = s_{2}, k_{2}=r_{1}-k_{1}$. So, $M_{2}$ has a $T(r_{1}, s_{1}, k_{1})$ representation. Therefore, by $f_{1}$, $T_{1} \cong T_{2}$.

\textbf{Case 3 :} If $r_{1} \neq r_{2}$, it implies that $a_{1, 1} \neq a_{2, 1}$. So by assumption $\{a_{1, 1}, a_{1, 2}\}$ = $\{a_{2, 1}, a_{2, 2}\}$, we get that $a_{1, 1}$ = $a_{2, 2}$. In this case, we identify boundaries of $T(r_{2}, s_{2}, k_{2})$ and cut along the hole cycle $C(2, 2)$ in place of $C(2, 1)$. Then, take another cut along $C(2, 1)$ until we reaching $C(2, 2)$ again for the first time. Hence, we get $r_{1}$ = length($C(1, 1)$) = length($C(2, 2)$) = $r_{2}$. Thus, $r_{1}s_{1} = r_{2}s_{2}$ implies that $s_{1} = s_{2}$. Since $r_1 = r_2$ and $s_1 = s_2$, it implies that we are in Case 2. So similarly as in Case 2, we define maps $f_{1}$, $f_{2}$ and $f_{3}$. Thus, by $f_{1}$,  $f_{2}$ and $f_{3}$, $T_{1} \cong T_{2}$.       

\textbf{Case 4 :} If $s_{1} \neq s_{2}$ then it implies that $n_{1} = r_{1}s_{1} \neq r_{2}s_{2} = n_{2}$ if $r_1 = r_2$. Which is a contradiction as $n_1 = n_2$. If $r_{1} \neq r_{2}$ then we are in Case 3. So, by combining above two Cases 2 and 3, we get an isomorphism map $f_1$ if $s_1 \not= s_2$. Again, let $k_1 \neq k_2$. Here, we have the following cases. If $r_1 \neq r_2$ then we are in Case 3. Similarly, if $s_1 \neq s_2$ then we are in Case 4. If $r_1 = r_2$, $s_1 = s_2$ and $k_1 \neq k_2$ then we are in Case 2. Thus, $(a_{1, 1}, a_{1, 2})$ = $(a_{2, t_1}, a_{2, t_2})$ where $t_1, t_2 \in \{1, 2\}$ defines $T_{1} \cong T_{2}$. So, by Cases 1, 2, 3, 4, the claim follows. 

So, by $f_1$, $M_1 \cong M_2$.

Conversely, let $M_{1} \cong M_{2}$. Similarly as in Lemma \ref{lem3342-iso}, let $f : V(M_{1})\rightarrow V(M_{2})$ such that $C_{2, j} := f(C_{1, j})$ for $j \in \{1, 2\}$. So, $\{a_{1, 1}, a_{1, 2}\}$ = $\{a_{2, 1}, a_{2, 2}\}$. Thus, it implies that $(a_{1, 1}, a_{1, 2})$ = $(a_{2, t_1}, a_{2, t_2})$ where $t_1 \ne t_2 \in \{1, 2\}$. 
\end{proof}


\begin{table}
\tiny 
\caption{Maps of type $\{3^{2}, 4, 3, 4\}$}
\centering 
\begin{tabular}{c c c c} 
\hline\hline
$n$&Equivalence classes& Length of cycles & $i(n)$\\ [0.3ex] 
\hline 
16 & T(8, 2, 4), T(4, 4, 2) & (8, 4) & 2(16) \\
   & T(4, 4, 0) & (4, 4)& \\
\hline
20 & T(10, 2, 4), T(10, 2, 6) & (10, 10) & 1(20) \\
\hline
24 & T(12, 2, 4), T(12, 2, 8), & (12, 6) & 3(24) \\
   & T(6, 4, 2), T(6, 4, 4)& &\\
   & T(12, 2, 6), T(4, 6, 2) & (4, 12)&\\
   & T(6, 4 ,0), T(4, 6, 0) & (4, 6)&\\
   \hline
28 & T(14, 2, 4), T(14, 2, 10) & (14, 14) & 1(28)\\
   & T(14, 2, 6), T(14, 2, 8) & &\\
   \hline
32 & T(16, 2, 4), T(16, 2, 12), & (16, 8) & 5(32)\\
   &T(8, 4, 2), T(8, 4, 6)& &\\
   & T(16, 2, 6), T(16, 2, 10) & (16, 16)&\\
   & T(16, 2, 8), T(4, 8, 2) & (4, 16)&\\
   & T(8, 4, 4) & (8, 8)&\\
   & T(8, 4, 0), T(4, 8, 0) & (4, 8)&\\
[1ex] 
\hline\hline 
\end{tabular}
\label{table2} 
\end{table}

As in Section \ref{33421}, by Lemmas \ref{lem32434-rep}, \ref{lem32434-iso}, the maps of type $\{3^{2}, 4, 3, 4\}$ can be classified up to isomorphism. We have done calculation for the vertices $\le 32$ and listed the obtained objects in the form of their $T(r, s, k)$ representation in Table \ref{table2}.

\section{Maps of type $\{3, 6, 3, 6\}$}\label{36361}

Let $M$ be a semi-equivelar map of type $\{3, 6, 3, 6\}$ on the torus. We define path in $M$ as follows. Through each vertex in $M$ there are two paths of type $X_1$ as shown in Figure 20.

\begin{defn}\label{defn3636} Let $Q_1 := P(\cdots,u_{i-1},u_i,u_{i+1},\cdots)$ be a path in edge graph of $M$. Let $A(v)$ denote a set of incident edges through $v$ in $M$. We say the path $Q_1$ of type $X_1$ if $A(u_{i})\setminus E(Q_1)$ is a set of two edges where one edge lie on one side and remaining one lie on the other side of $P(u_{i-1}, u_{i}, u_{i+1})$.
\end{defn}

\vspace{-.7cm}

\begin{picture}(0,0)(-50,37)
\tiny
\setlength{\unitlength}{1.5 mm}
\drawpolygon(5,5)(10,5)(12.5,10)(10,15)(5,15)(2.5,10)
\drawpolygon(15,5)(20,5)(22.5,10)(20,15)(15,15)(12.5,10)

\drawline[AHnb=0](10,15)(20,15)
\drawline[AHnb=0](10,5)(20,5)

\drawline[linecolor=black,linewidth=0.2,AHnb=0](10,5)(15,15)
\drawline[linecolor=black,linewidth=0.2,AHnb=0](10,15)(15,5)

\drawline[AHnb=0,dash={1.0 1.0 1.0 1.0}{0.0}](10,5)(7.2,0)
\drawline[AHnb=0,dash={1.0 1.0 1.0 1.0}{0.0}](15,5)(17.8,0)

\drawline[AHnb=0,dash={1.0 1.0 1.0 1.0}{0.0}](10,15)(7.2,20)
\drawline[AHnb=0,dash={1.0 1.0 1.0 1.0}{0.0}](15,15)(17.8,20)

\put(18,21){\scriptsize${X_1}$}
\put(6,21){\scriptsize${X_1}$}

\put(4.8,3.5){\scriptsize${a}$}
\put(10.1,3.8){\scriptsize${u_{i-1}}$}
\put(14.4,3.3){\scriptsize${b}$}
\put(19.8,3.5){\scriptsize${c}$}
\put(23,10){\scriptsize${d}$}
\put(13.2,10){\scriptsize${u_i}$}
\put(1.5,10){\scriptsize${e}$}
\put(4.8,16){\scriptsize${f}$}
\put(9.8,16){\scriptsize${g}$}
\put(11.6,16){\scriptsize${u_{i+1}}$}
\put(19.8,16){\scriptsize${h}$}

\put(6,-4){\scriptsize Figure 20 : $lk(u_i)$}

\end{picture}

\vspace{4.5cm}

For example in Figure 20, $P(u_{i-1}, u_i, u_{i+1})$ and $P(b, u_i, g)$ are two paths through $u_i$ and both are part of paths of type $X_1$. Let $P(u_1, \cdots, u_r)$ be a maximal path of type $X_{1}$ in $M$. Now, consider vertex $u_r$, $lk(u_r)$ and argue similarly as in Lemma \ref{lem3342-cycle}. So, we get an edge $u_1u_r$ such that $P \cup \{u_1u_r\}$ is a cycle of type $X_1$. So, there are two cycles of type $X_{1}$ through a vertex. Let $L_{1}(v),$ $L_{2}(v)$ be two cycles of type $X_1$ through a vertex $v$. We proceed with similar argument as in Section \ref{33421} and we get a connected $T(r, s, k)$ representation of $M$. In this process, we cut $M$ along the cycles $L_{1}(v)$ and $L_{2}(v)$ where we take second cut along the cycle $L_{2}$ and the starting adjacent face to the base horizontal cycle $L_{1}$ is a $3$-gon. Thus, every map $M$ has a $T(r, s, k)$ representation. The Figure 10 is an example of $T(8, 2, 6)$ representation of a map with $24$ vertices on the torus.

Now, we show that map of type $\{3, 6, 3, 6\}$ contains three cycles of type $X_{1}$ up to homologous.

\begin{lemma} \label{lem3636-hom-length} The map $M$ contains at most three cycles of type $X_{1}$ of different lengths.
\end{lemma}

\begin{proof}
Let $\triangle(u, v, w)$ be a $3$-gon in $M$. The $\triangle(u, v, w)$ has three edges $e_{1} = uv$, $e_{2} = vw$ and $e_{3} = uw$. By the definition of cycle of type $X_{1}$, $M$ contains at least three cycles, say $C_{1}$, $C_{2}$ and $C_{3}$ where $C_{i}$ contains edge $e_{i}$ for $i \in \{1, 2, 3\}$ and $C_{i}$ does not contain $e_{j}$ for $j \neq i$. Since $C_{i}$ does not contain $e_{j}$ for $j \neq i$, so, cycles are not identical. Again since, cycles are not identical and $V(C_{i})\cap C(C_{j})$ is a vertex of $\triangle$ for $i \neq j$, so, cycles are not homologous. (In Figure 10, $v_1v_2v_{1, 2, 16, 9}$ denote a face and the cycles which are of type $X_1$ and contains the edges $v_1v_2$, $v_1v_{1, 2,16, 9}$ and $v_{1, 2,16, 9}v_2$ are namely, $L_1 = C(v_1, v_2, \cdots, v_8),$ $L_2 = C(v_1, v_{1, 2,16, 9}, v_9, v_{9, 10, 6, 7}, \cdots, v_{11, 12, 8, 1})$ and $L _3 = C(v_2, v_{1, 2,16, 9}, v_{16}, v_{15, 16, 4, 5}, \cdots, v_{13, 14, 2, 3})$ respectively. The cycles $L_1$, $L_2$ and $L_3$ in Figure 10 are not homologous to each other.) Let $\triangle_{1}$ be an $3$-gon in $T(r, s, k)$ and $\triangle \neq \triangle_1$. Observe that there is a cycle of type $X_{1}$ through an edge $\triangle_{1}$ and homologous $C_{i}$ for some $i$. That is, there is a cylinder which is bounded by two cycles of type $X_1$ and containing edges of $\triangle, \triangle_1$. This is true for any $3$-gon in $T(r, s, k)$. We proceed as in the case of Lemma \ref{lem3342-all-length} and thus, the homologous cycles of type $X_1$ have same length. So, there are three different cycles $C_1, C_2, C_3$ of type $X_1$ up to homologous in $M$. So, the map $M$ contains at most three cycles of type $X_{1}$ of different lengths. (In Figure 10, consider face $v_{11}v_{12}v_{11, 12, 8, 1}$ and cycles of type $X_1$ which are containing the edges $v_{11}v_{12}$, $v_{11}v_{11, 12, 8, 1}$ and $v_{12}v_{11, 12, 8, 1}$ are namely, $L_1' = C(v_9, v_{10}, \cdots, v_{16}),$ $L_2' = C(v_{11}, v_{11, 12, 8, 1}, v_1, v_{1, 2, 16, 9}, \cdots, v_{3, 4, 10, 11})$ and $L_3' = C(v_{12}, v_{11, 12, 8, 1}, v_{8}, v_{7, 8, 14, 15}, \cdots, v_{5, 6, 12, 13})$ respectively. The cycles which are containing the edges $v_1v_2$ and $v_{11}v_{12}$ are $L_1$ and $L_1'$ respectively, and the cycles $L_1$ and $L_1'$ are homologous. The cycles which are containing $v_1v_{1, 2,16, 9}$ and $v_{11}v_{11, 12, 8, 1}$ are $L_2$ and $L_2'$ respectively and $C_2 = C_2'$. Also, the cycles which are containing $v_{2}v_{1, 2, 16, 9}$ and $v_{12}v_{11, 12, 8, 1}$ are $L_3$ and $L_3'$ respectively and $L_3 = L_3'$.)
\end{proof} 

We define admissible relations among $r, s, k$ of $T(r, s, k)$ in the next lemma. In this lemma we omit some cases as thier similar cases are discussed in the previous sections.

\begin{lemma}\label{lem3636-1}
The maps of type $\{3, 6, 3, 6\}$ of the form $T(r, s, k)$ exist if and only if the following holds : (i) $s \geq  1$, (ii) $2 \mid r$, (iii) number of vertices of $T(r, s, k)$ = $\frac{3}{2}rs \geq  21$, (iv) $r \geq 14$ if $s = 1$, (v) $r \geq 8$ if $s = 2$, (vi) $r \geq 6$ if $s \geq 3$, (vii) $k \in \{2t+6 \colon 0 \leq t \leq \frac{r-10}{2}\}\setminus \{2(\frac{r-10}{4})+6\}$ if $s = 1$, $k \in \{2t+6 \colon 0 \leq t \leq \frac{r-8}{2}\}$ if $s = 2 $ $\&$ $ k \in \{2t \colon 0 \leq t < \frac{r}{2}\}$ if $s \geq 3$.
\end{lemma}

\begin{proof} Let $C_{0}(u_{0,0}, u_{0,1}, \cdots, u_{0, r-1})$, $C_{1}(u_{1, 0}, u_{1, 1}, \cdots, u_{1, r-1})$,$\cdots$, $C_{s-1}(u_{s-1, 0}, u_{s-1, 1}, \cdots,\linebreak u_{s-1, r-1})$ be horizontal cycles of type $X_{1}$ in $T(r, s, k)$. By the definition of $T(r, s, k)$, $T(r, s, k)$ contains $s$ number of horizontal cycles of type $X_1$. Observe that, the number of adjacent vertices which are lie on one side of a horizontal cycle and not belong to any horizontal cycles is $\frac{r}{2}$. So, the total number of vertices in $T(r, s, k)$ is $(r+\frac{r}{2})s$. This implies that $n = (r+\frac{r}{2})s = \frac{3}{2}rs$.

By Euler's formula, the number of $6$-gons in $T(r, s, k)$ is $2n/6$ and it is an integer. This implies that $6 \mid 2n$. So, $3 \mid \frac{3}{2} sr$ as $n = \frac{3}{2}rs$. Thus, $2 \mid r$ if $s = 1$. Again, if $s \ge 2$ and $2 \nmid r$, it gives that the link $lk(u_{1})$ is not type $\{3, 6, 3, 6\}$. Which is a contradiction. So, $2 \mid r ~\forall~ s \ge 1$.       

Let $s = 1$. If $r < 14$ then $r \in \{2, 4, 6, 8, 10, 12\}$. If $r = 2, 4, 6, 8, 10, 12$ then there is a vertex in $T(r, s, k)$ whose link is not a cycle. So, $r \geq 14$ if $s = 1$.
Similarly as above, we get that $r \geq 8$ if $s = 2$ and $r \geq 6$ if $s \geq 3$. So, $\frac{3}{2}rs \geq  21$.

If $s = 1$ and $k \in\{t \colon 0\leq t\leq r-1\}\setminus(\{2t+6 \colon 0 \leq t \leq \frac{r-10}{2}\}\setminus  \{2(\frac{r-10}{4})+6\})$ then similarly as above we get some vertex whose link is not a cycle. Similarly, we repeat same argument as above for other two cases. So, $k \in \{2t+6 \colon 0 \leq t \leq \frac{r-8}{2}\}$ if $s = 2$ and $k \in \{2t \colon 0 \leq t \leq \frac{r}{2}\}$ if $s \geq 3$.
\end{proof}

Let $M_{i}$, $i = 1, 2$ be maps of type $\{3, 6, 3, 6\}$ on $n_i$ number of vertices and $n_1 = n_2$. Let $C_{i, j}$, $j = 1, 2, 3$ denote cycles which are of type $X_{1}$ and non-homologous in $T_i = T(r_{i}, s_{i}, k_{i})$. 
Let $a_{i, j}$ = length($C_{i, j}$) for $i= 1, 2$ and $j = 1, 2, 3$.
Then,

\begin{lemma}\label{lem3636-iso} The map $M_{1} \cong M_{2}$ if and only if $(a_{1, 1}, a_{1, 2}, a_{1, 3})$ = $(a_{2, t_1}, a_{2, t_2}, a_{2, t_3})$ for $t_1 \neq t_2 \neq t_3 \in \{1, 2, 3\}$.
\end{lemma}

\begin{proof} We first assume that $(a_{1, 1}, a_{1, 2}, a_{1, 3}) = (a_{2, t_1}, a_{2, t_2}, a_{2, t_3})$ where $t_i \in \{1, 2, 3\}$ and $t_i \neq t_j$. This implies that $\{a_{1, 1}, a_{1, 2}, a_{1, 3}\} = \{a_{2, 1}, a_{2, 2}, a_{2, 3}\}$.

\smallskip

$Claim : T_{1} \cong T_{2}$.

\smallskip

\textbf{Case 1 :} If $(r_{1}, s_{1}, k_{1}) = (r_{2},s_{2},k_{2})$ then $T(r_{1}, s_{1}, k_{1}) = T(r_2, s_2, k_2) = T(r, s, k)$. Let $C(1, 0) := C(u_{0,0}, u_{0,1}, \cdots, u_{0,r-1})$, $C(1, 1) :=C(u_{1,0}, u_{1,1}, \cdots u_{1,r-1})$, $\cdots$, $C(1, s) :=C(u_{s,0}, u_{s,1}, \cdots, u_{s, r-1}\}$ denote sequence of horizontal cycles of type $X_{1}$ in $T_1$. Again, let $G_1(t, t+1) := \{w_{t,0},  w_{t,1}, \cdots,  w_{t,\frac{r-2}{2}})$ which is a set of vertices where $ w_{t,i}$ is adjacent with both $u_{t, 2i}$ and $u_{t+1, 2i}$ in $T_1$ and not belong to both $C(1, t)$ and $C(1, t+1)$ for $0\leq t \leq s$. (For example in Figure 10, $G_1(0, 1) = \{v_{1, 2, 16, 9}, v_{3, 4, 10, 11}, v_{5, 6, 12, 13}, v_{7, 8, 14, 15}\}$ where $ x_{0, 0} = v_{1, 2, 16, 9},  x_{0, 1} = v_{3, 4, 10, 11},   x_{0, 2} = v_{5, 6, 12, 13},  x_{0, 3} = v_{7, 8, 14, 15}$.) Similarly, let $C(2, 0) := C(v_{0, 0}, v_{0,1}, \cdots, v_{0,r-1})$, $C(2, 1) := C(v_{1,0}, v_{1,1}, \cdots v_{1,r-1})$, $\cdots, C(2, s) :=C(v_{s,0}, v_{s,1}, \cdots $ $ v_{s, r-1})$ denote sequence of horizontal cycles of type $X_{1}$ in $T_2$ and define $G_2(t, t+1) :=\{x_{t,0},  x_{t, 1}, \cdots, $ $  x_{t,\frac{r-2}{2}}\}$ which is a set of vertices where the vertex $x_{t, i}$ is adjacent with both $v_{t,2i}$ and $v_{t+1,2i}$ in $T_2$ for $0\leq t \leq s$. Now we define an isomorphism map $f : V(T(r_{1}, s_{1}, k_{1}))\rightarrow V(T(r_{2}, s_{2}, k_{2}))$ such that $f(u_{t, i}) = v_{t, i} ~\forall~ 0 \leq i\leq r-1$, $0 \leq t\leq s-1$ and $f( w_{t, i}) =  x_{t, i}$ for the vertices of $G_1(t, t+1)$ and $G_2(t, t+1)$ for all $0 \leq t\leq s-1$. By $f$, the link $lk(u_{t, i})$ maps to the link of $lk(v_{t, i})$ and $lk(w_{t, i})$ maps to $lk( x_{t, i})$. So, by $f$, $T_{1} \cong T_{2}$.

\textbf{Case 2 :} Let $(r_{1}, s_{1}, k_{1}) \neq (r_{2},s_{2},k_{2})$. If $r_1 \ne r_2$, we identify boundaries of $T(r_{2}, s_{2}, k_{2})$ and cut $M_2$ along cycle of length $r_1$ and followed by, make another cut along a cycle of type $X_1$ to get a $(r, s, k)$-representation. Thus we get a new $T(r_2', s_2', k_2')$ representation of $M_2$. It implies that $r_{1} = r_2'$ and $s_{1} = s_{2}'$ as $n_{1} = \frac{3}{2}r_{1}s_{1} = \frac{3}{2}r_{2}'s_{2}'=n_{2}$. By this process, we get a new representation $T(r_{1}, s_{1}, k_{3}')$ of $M_{2}$.
If $k_{1} = k_{3}'$ then $M_{1} \cong M_{2}$ by $f$ in Case 1. If $k_{1} \neq k_{3}'$ similarly as in Lemma \ref{lem32434-iso}, we make a cut along a path which is homologous to boundary path and identify along the boundary path. Thus, we  get an another representation $T(r_1, s_1, k_3'')$ of $M_2$ and $k_{1} = k_{3}''$. So, the $M_{2}$ has a $T(r_{1}, s_{1}, k_{1})$ representation since $k_{1} = k_{3}''$. Therefore, there exists $f$ and $T_{1} \cong T_{2}$ by $f$. This completes the Claim.

So, by $f$, $M_1 \cong M_2$.

Conversely, let $M_{1} \cong M_{2}$. We proceed as in the converse part of Lemma \ref{lem32434-iso} and we get $\{a_{1, 1}, a_{1, 2}, a_{1, 3}\}$ = $\{a_{2, 1}, a_{2, 2}, a_{2, 3}\}$. That is, $(a_{1, 1}, a_{1, 2}, a_{1, 3})$ = $(a_{2, t_1}, a_{2, t_2}, a_{2, t_3})$ for $t_1 \neq t_2 \neq t_3 \in \{1, 2, 3\}$.
\end{proof}

As in Section \ref{33421}, by Lemmas \ref{lem3636-1}, \ref{lem3636-iso}, the maps of type $\{3,6,3,6\}$ can be classified up to isomorphism on different number of vertices. We have done the calculation for the vertices $\le 30$. We have listed the obtained objects in the form of their $T(r, s, k)$ representation in Table \ref{table3}.

\begin{table}
\tiny 
\caption{Maps of type $\{3, 6, 3, 6\}$}
\centering 
\begin{tabular}{c c c c} 
\hline\hline 
$n$&Equivalence classes& Length of cycles & $i(n)$\\ [0.3ex] 
\hline 
21 & T(14, 1, 6), T(14, 1, 10) & (14, 14, 14) & 1(21) \\
\hline
24 & T(16, 1, 6), T(16, 1, 12) & (16, 16, 8) & 2(24) \\
&T(8, 2, 6)& &\\
 & T(16, 1, 8), T(16 , 1, 10) & (16, 16, 4) &\\
\hline
27 & T(18, 1, 6), T(18, 1, 8) &(18, 18, 6) & 2(27)\\
& T(18, 1, 12), T(18, 1, 14) & &\\
& T(6, 3, 2), T(6, 3, 4) & &\\
 & T(6, 3, 0) & (6, 6, 6) &\\
\hline
30 & T(20, 1, 6), T(20, 1, 8) & (20, 20, 10) & 2(30)\\
& T(20, 1, 14), T(20, 1, 16)& &\\
& T(10, 2, 2), T(10, 2, 6) & &\\
& T(10, 2, 8) & &\\

& T(20, 1, 10), T(20, 1, 12) & (20, 4, 10) &\\
&T(10, 2, 0), T(10, 2, 4)& &\\
[1ex] 
\hline\hline 
\end{tabular}
\label{table3} 
\end{table}

\section{Maps of type $\{3, 12^2\}$}\label{31221}

Let $M$ be a semi-equivelar map of type $\{3, 12^2\}$ on the torus. We define a fixed type of path $G_1$ in the edge graph of $M$ as shown in Figure 21. Let $Q(i) := P(u_{i}, u_{i+1}, u_{i+2}, u_{i+3}, u_{i+4})$ be a path in $M$ where \textit{$lk(u_{i})=C(u_{i-1}, \textbf{a}, \textbf{b}, \textbf{c}, \textbf{d}, \textbf{e}, \textbf{f}, \textbf{g}, \textbf{g}', \textit{\textbf{u}}_{i+2}, u_{i+1}, u, \textit{\textbf{t}}, \textit{\textbf{v}}, \textit{\textbf{w}}, \textit{\textbf{x}}, \textit{\textbf{y}}, \textit{\textbf{z}}, \textbf{a}', \textbf{u}_{i-3},$ $\textbf{u}_{i-2})$, $lk(u_{i+1})=C(u_i, \textit{\textbf{u}}_{i-1}, \textit{\textbf{a}}, \textbf{b}, \textbf{c}, \textbf{d}, \textbf{e}, \textbf{f}, \textbf{g}, \textbf{g}', u_{i+2}, \textit{\textbf{u}}_{i+3}, \textit{\textbf{u}}_{i+4}, \textbf{o}',$  $\textit{\textbf{o}}, \textbf{p}, \textbf{q}, \textit{\textbf{r}}, \textbf{s}, \textbf{t}, u)$, \newline $lk(u_{i+2})=C(u_{i+1}, \textit{\textbf{u}}_{i}, \textit{\textbf{u}}_{i-1}, \textbf{a}, \textbf{b}, \textbf{c}, \textbf{d}, \textbf{e}, \textbf{f}, \textbf{g}, g', u_{i+3}, \textit{\textbf{u}}_{i+4}, \textbf{o}',$  $\textit{\textbf{o}}, \textbf{p}, \textbf{q}, \textbf{r}, \textbf{s}$, $\textbf{t}, \textbf{u})$,
$lk(u_{i+3})=C(u_{i+2}, {g}', \textbf{g},$ $\textit{\textbf{h}}, \textbf{i}, \textbf{j}, \textbf{k}, \textbf{l}, \textbf{m}, \textbf{n}, \textit{\textbf{u}}_{i+6}, \textit{\textbf{u}}_{i+5}, u_{i+4}, \textbf{o}', \textit{\textbf{o}}, \textbf{p}, \textbf{q}, \textbf{r}, \textbf{s}, \textbf{t}, \textbf{u} \textit{\textbf{u}}_{i+1})$, and \newline
$lk(u_{i+4})=C(u_{i+3}, \textbf{g}', \textit{\textbf{g}}, \textbf{h}, \textbf{i}, \textbf{j},$ $ \textbf{k}, \textbf{l}, \textbf{m}, \textbf{n}, \textit{\textbf{u}}_{i+6}, u_{i+5}, o', \textbf{o}, \textit{\textbf{p}}, \textbf{q}, \textbf{r}, \textbf{s}, \textbf{t}, \textbf{u}, \textit{\textbf{u}}_{i+1}, \textit{\textbf{u}}_{i+2})$}.

\begin{defn}\label{defn3122-1} Let $R_1 := P(\cdots,v_{i-1},v_i,v_{i+1},\cdots)$ be a path in edge graph of $M$. We say $R_1$ of type $G_1$ if $L_1 :=P(v_{t}, v_{t+1}, v_{t+2}, v_{t+3}, v_{t+4})$ is a subpath of $R_1$ or is a part in an extended path of $R_1$ then either $L_1\mapsto Q(i)$ by $v_{j}\mapsto u_{j}$, $L_1 \mapsto Q(i+1)$ by $v_{j}\mapsto u_{j+1}$, $L_1 \mapsto Q(i+2)$ by $v_{j}\mapsto u_{j+2}$ or $L_1 \mapsto Q(i+3)$ by $v_{j}\mapsto u_{j+3}$ for $j \in \{t, t+1, t+2, t+3, t+4\}$.
\end{defn}

\begin{defn}\label{defn3122-2} Let $R_2 := P(\cdots,x_{i-1},x_i,x_{i+1},\cdots)$ be a path in edge graph of $M$. We say $R_2$ of type $G_1'$ if $L_2 : = P(x_{t}, x_{t+1}, x_{t+2}, x_{t+3}, x_{t+4})$ is a subpath of $R_2$ or is a part in the extended path of $R_2$ then either $L_2 \mapsto Q(i)$ by $x_{j}\mapsto u_{2t+4-j}$, $L_2 \mapsto Q(i+1)$ by $x_{j}\mapsto u_{2t+4-j}$, $L_2 \mapsto Q(i+2)$ by $x_{j}\mapsto u_{2t+4-j}$ or $L_1 \mapsto Q(i+3)$ by $v_{j}\mapsto u_{2t+4-j}$ for $j \in \{t, t+1, t+2, t+3, t+4\}$.
\end{defn}

\vspace{0.8 cm}
\smallskip
\begin{picture}(50,-40)(-40,2)
\tiny
\setlength{\unitlength}{3.4mm}
\drawpolygon(-0.5,0)(1,-1)(2,-1)(3.5,0)(4,1)(4,2)(3.3,3)(2,4)(1,4)(-0.5,3)(-1,2)(-1,1)
\drawpolygon(4.5,0)(6,-1)(7,-1)(8.5,0)(9,1)(9,2)(8.3,3)(7,4)(6,4)(4.5,3)(4,2)(4,1)
\drawpolygon(9.5,0)(11,-1)(12,-1)(13.5,0)(14,1)(14,2)(13.3,3)(12,4)(11,4)(9.5,3)(9,2)(9,1)

\drawline[linecolor=black,linewidth=0.1,AHnb=0](1,-1)(2,-1)
\drawline[linecolor=black,linewidth=0.1,AHnb=0](2,-1)(3.5,0)
\drawline[linecolor=black,linewidth=0.1,AHnb=0](6,-1)(7,-1)
\drawline[linecolor=black,linewidth=0.1,AHnb=0](7,-1)(8.5,0)
\drawline[linecolor=black,linewidth=0.1,AHnb=0](11,-1)(12,-1)
\drawline[linecolor=black,linewidth=0.1,AHnb=0](12,-1)(13.5,0)

\drawline[linecolor=black,linewidth=0.1,AHnb=0](-.5,0)(1,-1)
\drawline[linecolor=black,linewidth=0.1,AHnb=0](4.5,0)(6,-1)
\drawline[linecolor=black,linewidth=0.1,AHnb=0](4.5,0)(3.5,0)
\drawline[linecolor=black,linewidth=0.1,AHnb=0](9.5,0)(11,-1)
\drawline[linecolor=black,linewidth=0.1,AHnb=0](9.5,0)(8.5,0)
\drawline[linecolor=black,linewidth=0.1,AHnb=0](14.5,0)(16,-1)
\drawline[linecolor=black,linewidth=0.1,AHnb=0](14.5,0)(13.5,0)

\drawpolygon (-1,2)(-1.7,3)(-0.5,3)
\drawpolygon (4,2)(3.3,3)(4.5,3)
\drawpolygon (9,2)(8.3,3)(9.5,3)

\drawpolygon (1,-1)(2,-1)(1.5,-2)
\drawpolygon (1.5,-2)(2,-1)(3.5,0)(4.5,0)(6,-1)(6.5,-2)(6.5,-3)(5.9,-4)(4.7,-5)(3,-5)(1.8,-4)(1.5,-3)
\drawpolygon (6.5,-2)(7,-1)(8.5,0)(9.5,0)(11,-1)(11.5,-2)(11.5,-3)(10.9,-4)(9.7,-5)(8,-5)(6.8,-4)(6.5,-3)
\drawpolygon (11.5,-2)(12,-1)(13.5,0)(14.5,0)(16,-1)(16.5,-2)(16.5,-3)(15.9,-4)(14.7,-5)(13,-5)(11.8,-4)(11.5,-3)

\drawpolygon (13.5,0)(14.5,0)(14,1)

\drawpolygon (16.5,-3)(15.9,-4)(16.8,-4)
\drawpolygon (11.5,-3)(10.9,-4)(11.8,-4)
\drawpolygon (6.5,-3)(5.9,-4)(6.8,-4)

\drawline[AHnb=0,dash={.5 .5 .5 .5}{0.0}](-.5,0)(-2,0)
\drawline[AHnb=0,dash={.5 .5 .5 .5}{0.0}](16,-1)(18.5,-1)

\put(-3,0){$G_1$}

\put(2.8,0){\scriptsize c}
\put(4.5,0){\scriptsize d}
\put(5.9,-.8){\scriptsize e}

\put(6.9,-.7){\scriptsize f}
\put(7.9,0){\scriptsize g}
\put(9.5,0){\scriptsize h}
\put(10.8,-.9){\scriptsize i}

\put(11.6,-.8){\scriptsize j}
\put(12.7,0){\scriptsize k}
\put(14.5,0){\scriptsize l}
\put(16,-.9){\scriptsize m}

\put(.7,-.8){\scriptsize a}
\put(1.8,-.7){\scriptsize b}

\put(4,-6.5){\scriptsize  Figure 21 : Cycle of type $G_1$}

\end{picture}

\vspace{2.5cm}

We argue with the similar argument of Lemma \ref{lem3342-cycle} for the path of types $G_1, G_1'$ and so, every maximal path of types $G_1, G_1'$ is a cycle and non-contractible (by the similar argument of Lemma \ref{lem3342-non-con}). Observe that the cycles of type $G_1$ and $G_1'$ are mirror image of each other. So, these types define same type of cycle. (Similar argument of this is given details in Section \ref{3461} for the type $\{3^4, 6\}$.) Clearly, there are two cycles of type $G_{1}$ through each vertex of $M$. Let $uvw$ be a $3$-gon in $M$. Let $L_1(u,$ $uw)$, $L_2(w,$ $wv)$ and $L_3(v, vu)$ denote three cycles through $u,$ $w,$ and $v$ respectively where $L_1(u, uw)$ contains the edge $uw$, $L_2(w, wv)$ contains the edge $wv$ and $L_3(v, vu)$ contains the edge $vu$. We repeat the similar argument of Section \ref{33421} and define a $T(r, s, k)$ representation of the map $M$ for some $r, s, k$. In the process, we take first cut along $L_1(u,$ $uw)$, and then, second cut along $L_2(w,$ $wv)$ where the starting adjacent face to horizontal base cycle $L_{1}(u,$ $uw)$ is a $12$-gon. Let length$(L_1(u,$ $uw))) = r$, $s$ denote the number of homologous cycles of $L_1(u,$ $uw)$ of type $G_1$, and $k$ denote the distance of the starting vertex of upper horizontal cycle from the starting vertex $w$ in $L_1(u,$ $uw)$. By this process, we get a $T(r, s, k)$ representation of $M$. Now, we proceed with process of Section \ref{36361} for the map of type $\{3, 12^2\}$. We show that a map of type $\{3, 12^2\}$ contains at most three non-homologous cycles of type $G_{1}$ of different lengths.

\begin{lemma} \label{lem4.3.2} The map $M$ contains at most three cycles of type $G_{1}$ of different lengths.
\end{lemma}

\begin{proof}
We proceed as in the case of Lemma \ref{lem3636-hom-length} for the map of type $\{3, 12^2\}$. So, consider map of type $\{3, 12^2\}$ in place of $\{3, 6, 3, 6\}$, cycle of type $G_1$ in place of $X_1$ and $3$-gon in the proof of Lemma \ref{lem3636-hom-length}. Thus, we get three non-homologous cycles of type $G_{1}$ of different lengths. 
\end{proof} 

We define admissible relations among $r,$ $s,$ $k$ of $T(r, s, k)$ in $M$.

\begin{lemma}\label{lem3122-all} The maps of type $\{3, 12^2\}$ of the form $T(r, s, k)$ exist if and only if the following holds : (i) $s \geq  1$, (ii) $4 \mid r$, (iii) number of vertices of $T(r, s, k)$ = $\frac{3}{2}rs \geq  36$, (iv) $r \geq 24$ if $s = 1$, (v) $r \geq 16$ if $s = 2$, (vi) $r \geq 12$ if $s \geq 3$, (vii) $k \in \{ 4t+9 \colon 0 \leq t \leq \frac{r-20}{4}\}\setminus \{ 4 (\frac{r}{8}-3) + 9\}$ if $s = 1$, $k \in \{ 4t+5 \colon 0 \leq t \leq \frac{r-16}{4}\}$ if $s = 2 $ $\&$ $ k \in \{ 4t+1 \colon 0 \leq t \leq \frac{r-4}{4}\}$ if $s \geq 3$.
\end{lemma}

\begin{proof}
We proceed as in the case of proof of the Lemma \ref{lem3636-1}. We consider map of type $\{3, 12^2\}$ in place of $\{3, 6, 3, 6\}$ and different values of $r, s, k$. Thus, we get all the cases.
\end{proof} 

Let $T_i = T(r_i, s_i, k_i)$ be representations of $M_i$, $i = 1, 2$ on same number of vertices. Similarly as in Section \ref{36361}, let $b_{i, j}$ = length($L_{i, j}$), $j=1, 2, 3$. Similarly as in Section \ref{36361}, we have 

\begin{lemma}\label{lem3122-iso} The map $M_{1} \cong M_{2}$ if and only if $(b_{1, 1}, b_{1, 2}, b_{1, 3})$ = $(b_{2, t_1}, b_{2, t_2}, b_{2, t_3})$ for $t_1 \neq t_2 \neq t_3 \in \{1, 2, 3\}$.
\end{lemma}


As in Section \ref{33421}, by Lemmas \ref{lem3122-all}, \ref{lem3122-iso}, maps of type $\{3, 12^2\}$ can be classified up to isomorphism on different number of vertices. We have done the calculation for the vertices $\le 48$. We have listed the obtained objects in the form of their $T(r, s, k)$ representation in Table \ref{table4}.

\begin{table}
\tiny 
\caption{Maps of type $\{3, 12^2\}$}
\centering 
\begin{tabular}{c c c c} 
\hline\hline 
$n$&Equivalence classes& Length of cycles & $i(n)$\\ [0.3ex] 
\hline 
36 & T(24, 1, 13) & (24, 12, 8) & 1(36)\\
\hline
42 & T(28, 1, 9), T(28, 1, 13) & (28, 28, 28) & 1(42)\\
& T(28, 1, 17) & &\\
\hline
48 & T(32, 1, 9), T(32, 1, 21) & (32, 32, 16) & 2(48)\\
& T(16, 2, 5)& &\\
 & T(32, 1, 17) & (32, 32, 8)&\\
[1ex] 
\hline\hline 
\end{tabular}
\label{table4} 
\end{table}

\section{Maps of type $\{3^{4}, 6\}$}\label{3461}

Let $M$ be a semi-equivelar map of type $\{3^{4}, 6\}$ on the torus. Let $Q(i) := P(w_{i}, w_{i+1},$ $w_{i+2},$ $w_{i+3})$ be a path in $M$, where $lk(w_{i})=C(w_{i-1}, x_{2}, w_{i+1}, \textit{\textbf{w}}_{i+2},  \textit{\textbf{x}}_{3}, \textit{\textbf{x}}_{4},  x_{5}, x_{6})$, $lk(w_{i+1})=C(w_{i}, x_{2}, x_{7}, x_{8}, w_{i+2}, \textit{\textbf{x}}_{3},  \textit{\textbf{x}}_{4}, \textit{\textbf{x}}_{5})$, $lk(w_{i+2})=C(w_{i+1}, x_{8}, x_{9}, w_{i+3}, x_{3}, \textit{\textbf{x}}_{4},  \textit{\textbf{x}}_{5}, \textit{\textbf{w}}_{i},)$, $lk(w_{i+3})$ $= C(w_{i+2}, x_{9}, w_{i+4},  \textit{\textbf{w}}_{i+5}, \textit{\textbf{x}}_{10}, \textit{\textbf{x}}_{11},$ $x_{12}, x_{3})$.
We define two fixed types of paths $Y_1$ and $Y_1'$ in the edge graph of $M$.

\begin{defn}\label{defn346-1} Let $R_1 := P(\cdots,v_{i-1}, v_i, v_{i+1},\cdots)$ be a path in edge graph of $M$. We say $R_1$ of type $Y_1$ if $L_1 :=P(u_{t}, u_{t+1}, u_{t+2}, u_{t+3})$ is a sub-path of $R_1$ or in the extended path of $R_1$ then $ L_1\mapsto Q(i)$ by $u_{j}\mapsto w_{j}$,  $L_1 \mapsto Q(i+1)$ by $u_{j}\mapsto w_{j+1}$ or $L_1 \mapsto Q(i+2)$ by $u_{j}\mapsto w_{j+2}$ for $j \in \{t,$ $t+1,$ $t+2,$ $t+3\}$.
\end{defn}

\begin{defn}\label{defn346-2} Let $R_2 := P(\cdots,x_{i-1},x_i,x_{i+1},\cdots)$ be a path in edge graph of $M$. We say $R_2$ of type $Y_1'$ if $L_2 : = P(x_{t}, x_{t+1}, x_{t+2}, x_{t+3})$ is a sub-path of $R_2$ or in the extended path of $R_2$ then $L_2 \mapsto Q(i)$ by $x_{j}\mapsto w_{2t+3-j}$, $L_2 \mapsto Q(i+1)$ by $x_{j}\mapsto w_{2t+3-j}$ or $L_2 \mapsto Q(i+2)$ by $x_{j}\mapsto w_{2t+3-j}$ for $j \in \{t, t+1, t+2, t+3\}$.
\end{defn}

Let $P$ be a maximal path of type $Y_{1}$ or $Y_1'$. By the similar argument of Lemma \ref{lem3342-cycle}, path $P$ defines a cycle of type $Y_1$ or $Y_1'$, that is, there is an edge $e$ in $M$ such that $P \cup \{e\}$ is a cycle of type $Y_1$ or $Y_1'$. We show that the cycles of types $Y_1$ and $Y_1'$ define same type of cycle. Let $C_{1} :=C(u_{1}, u_{2}, \cdots, u_{r})$ of type $Y_{1}$ and $C_{2}(v_{1}, v_{2}, \cdots, v_{r})$ of type $Y_{1}'$ be two cycles of same length $r$. Let $P_1 :=P(u_{i-1}, u_{i}, u_{i+1})$ be a sub path of $C_{1}$ where adjacent $6$-gon lies on one side and all $3$-gons lie on the other side of $P_1$ at the vertex $u_{i}$. Similarly, let $P_2 :=P(v_{j-1}, v_{j}, v_{j+1})$ be a sub path of $C_{2}$ where adjacent $6$-gon lies on one side and all $3$-gons lie on the other side of $P_2$ at the vertex $v_{j}$. Define a map $f : V(C_1) \rightarrow V(C_2)$ by $f(u_{i}) = v_{j},$ $f(u_{i+1}) = v_{j-1},$ $f(u_{i+2}) = v_{j-2}, \cdots, f(u_{i-1}) = v_{j+1}$. Let $P(u_t, \cdots, u_k)$ be a sub path of $C_1$. Then, $P(u_t, \cdots, u_k)$ and $P(f(u_t), \cdots, f(u_k))$ divide the link of the vertices $u_i$ and $f(u_i)$ for $t \leq i \leq k$ into the same ratio. This is true for every sub path of $C_1$. Therefore, cycles $C_1$ and $C_2$ are of type $Y_1$, and hence, $Y_1 = Y_1'$.

Let $M$ be a map and $C$ be a cycle of type $Y_1$ in $M$. By the similar argument of Lemma \ref{lem3342-non-con}, the cycle $C$ is non-contractible. As in Section \ref{33421}, we get that the cycles of type $Y_{1}$ which are homologous to $C$ have same length by the similar argument of Lemma \ref{lem3342-all-length}. There are three cycles of type $Y_1$ through each vertex of $M$. Let $v \in V(M)$ and $L_{1}(v),$ $L_{2}(v),$ $L_{3}(v)$ be three cycles of type $Y_1$ through the vertex $v$. We repeat similar construction of $(r, s, k)$-representation of a map as in Section \ref{33421} for $M$. In this process, we take first cut along $L_1(v)$ and second cut along $L_{2}(v)$ where the starting adjacent face to horizontal base cycle $L_{1}$ is a $3$-gon. This gives a $T(r, s, k)$ representation of the map $M$. Thus, $T(r, s, k)$ representation exists for every $M$.

Now, we show that map $M$ of type $\{3^{4}, 6\}$ contains at most three non-homologous cycles of type $Y_{1}$ of different lengths in Lemma \ref{lem346-hom}.

\begin{lemma} \label{lem346-hom} The map $M$ contains at most three non-homologous cycles of type $Y_{1}$ of different lengths.
\end{lemma}

\begin{proof}
Let $v \in V(M)$ and $T(r, s, k)$ denote a $(r, s, k)$-representation of $M$. We have three cycles, namely, $C_{1}$, $C_{2}$ and $C_{3}$ through $v$ in $M$ of type $Y_{1}$. Cycles $C_{1}$, $C_{2}$ and $C_{3}$ are not identical as $C_{i}$ divides link $lk(v)$ into different ratio. Also, cycles are not disjoint as $v \in V(C_{i})\cap V(C_{j})$ for $i \neq j$ and $i, j \in \{1, 2, 3\}$. So, $C_1$, $C_2$ and $C_3$ are not homologous to each other. Let $ w \in V(M)$ and $v \neq w$. Consider cycles of type $Y_{1}$ at $w$ in $T(r, s, k)$ and denoted by $C_1',$ $ C_2',$ $C_3'$. Now, by the definition of cycle of type $Y_1$ and considering cylinder, $C_{i}$ and $C_{j}'$ are homologous for some $i, j \in \{1, 2, 3\}$ as we have seen the same idea in Lemma \ref{lem3636-hom-length}. This is hold for any vertex of $M$. Therefore, $M$ contains at most three non-homologous cycles of type $Y_{1}$. We proceed as in the case of proof of Lemma \ref{lem3342-all-length} to show that the homologous cycles of type $Y_1$ have same length. Thus, the map $M$ contains at most three non-homologous cycles of type $Y_{1}$ of different lengths.
\end{proof} 

We define admissible relations among $r$, $s$, $k$ of $T(r, s, k)$ such that representation $T(r, s, k)$ gives a map of type $\{3^{4}, 6\}$ after identifying their boundaries.

\begin{lemma}\label{lem346-all}  The maps of type $\{ 3^{4}, 6\}$ of the form $T(r, s, k)$ exist if and only if the following holds : (i) $s \geq  2$ even, (ii) $3 \mid r$, (iii) number of vertices of $T(r, s, k)$ = $rs \geq  18$, (iv) $r \geq 9$ if $s = 2$, (v) $r \geq 6$ if $s \ge 4$, (vi) $k \in \{3t+5 : 0 \leq t \leq \frac{r-9}{3}\}$ if $s = 2 $ $\&$ $ k \in \{ 2+3t ~:~ 0 \leq t \leq \frac{r-3}{3}\}$ if $s \geq 4$. 
\end{lemma}

\begin{proof} Let $T(r, s, k)$ be a representation of $M$. It has $s$ number of disjoint horizontal cycles of type $Y_{1}$ of length $r$ by the definition of $(r, s, k)$-representation. These cycles cover all the vertices of $M$. So, $n = rs$. By Euler's formula, $n - 5n/2+ 4n/3 + n/6 = 0$. So, the number of $6$-gons in $M$ is $n/6$ and which is an integer. This implies that $6 \mid n$. That is, $6 \mid rs$ as $n = rs$. So, $3 \mid r$ for $s = 2$. Let $s \ge 3$. If $s$ is an odd integer then $T(r, s, k)$ contains odd number of horizontal cycles of type $Y_{1}$. Consider a vertex $v$ of base horizontal cycle which is belongs to only triangles which is a contradiction. So, $2 \mid s$. Similarly, for $3 \mid r$, we get a vertex whose link does not follow the type $\{3^{4}, 6\}$. So, $n = rs$ where $6 \mid n$, $2 \mid s$ and $3 \mid r$.

For $r \geq 9$, we proceed with similar argument as it is done in the proof of Lemma \ref{lem3342-all-rep}. We also proceed as in the case of the proof of Lemma \ref{lem3342-all-rep} to show $r \geq 6$ and $3 \mid r$ and remaining all other cases. This completes the proof.
\end{proof}

\begin{table}
\tiny 
\caption{Maps of type $\{3^{4}, 6\}$}
\centering 
\begin{tabular}{c c c c} 
\hline\hline 
$n$&Equivalence classes& Length of cycles & $i(n)$\\ [0.3ex] 
\hline 
18 & T(9, 2, 5) & (9, 9, 9) & 1(18)\\
\hline
24 & T(12, 2, 5), T(12, 2, 8) & (12, 6, 12) & 2(24) \\
& T(6, 4, 2)& &\\

& T(6, 4, 5) & (6, 6, 6)& \\
\hline
30 & T(15, 2, 5), T(15, 2, 8) & (15, 15, 15) & 1(30)\\
& T(15, 2, 11) & &\\
\hline
36 & T(18, 2, 5), T(18, 2, 14) & (18, 9, 18) & 2(36)\\
&T(9, 4, 2)& &\\

& T(18, 2, 8), T(18, 2, 11) & (18, 6, 9)&\\
& T(9, 4, 5), T(9, 4, 8)& &\\
& T(6, 6, 2), T(6, 6, 5) & &\\
\hline
42 & T(21, 2, 5), T(21, 2, 8) & (21, 21, 21) & 1(42)\\
&T(21, 2, 11), T(21, 2, 14)& &\\
&T(21, 2, 17)& &\\
[1ex] 
\hline\hline 
\end{tabular}
\label{table5} 
\end{table}

Let $M_{1}$ and $M_{2}$ be two maps of type $\{3^{4}, 6\}$ on same number of vertices. Similarly as in Sections \ref{33421}, \ref{324341},  let $a_{i, j}$ = length($C_{i, j}$) where $C_{i, k}$ for $i = 1, 2, 3$ denote three non-homologous cycles of type $Y_{1}$ in $T(r_{i}, s_{i}, k_{i})$ of $M_i$. Then, 

\begin{lemma}\label{lem346-iso} The map $M_{1} \cong M_{2}$ if and only if $(a_{1, 1}, a_{1, 2}, a_{1, 3})$ = $(a_{2, t_1}, a_{2, t_2}, a_{2, t_3})$ for $t_1 \neq t_2 \neq t_3 \in \{1, 2, 3\}$.
\end{lemma}

\begin{proof}
Let $T(r_{i},s_{i},k_{i})$ be a representation of $M_{i}$. If $r = r_{1}, s = s_{1}, k = k_{1}$, then, consider horizontal cycles in $T(r_{i},s_{i},k_{i})$ of type $Y_{1}$. Proceed with similar argument as in the Lemma \ref{lem3342-iso}. So, we get a map which defines isomorphism between $T(r_{1},s_{1},k_{1})$ and $T(r_{2},s_{2},k_{2})$. So, $M_{1} \cong M_{2}$. Again, if $(r, s, k) \neq (r_{1},s_{1},k_{1})$ then proceed with similar argument as in the proof of the Lemma \ref{lem32434-iso}. Converse part of the lemma follows from similar argument of the converse part of Lemma \ref{lem3342-iso}. This completes the proof.
\end{proof} 

As in Section \ref{33421}, by Lemmas \ref{lem346-all}, \ref{lem346-iso}, the maps of type $\{3^{4}, 6\}$ can be classified up to isomorphism on different number of vertices. We have done calculation for the vertices $\le 42$. We have listed the obtained objects in the form of their $T(r, s, k)$ representation in Table \ref{table5}.

\section{Maps of type $\{4, 6, 12\}$}\label{46121}

Let $M$ be a semi-equivelar map of type $\{ 4, 6, 12\}$ on the torus. We define a fixed type of path $H_1$ in the edge graph of $M$. Let $Q(i) := P(u_{i}, u_{i+1}, u_{i+2}, u_{i+3}, u_{i+4}, u_{i+5}, u_{i+6})$ be a path in $M$ where \textit{$lk(u_{i})=C(u_{i-1}, \textbf{b}, \textbf{c}, \textit{\textbf{u}}_{i+2}, u_{i+1} \textbf{p}, q, \textbf{r}, \textbf{s}, \textbf{t}, \textbf{u}, \textbf{v}, \textit{\textbf{u}}_{i-5},$ $\textit{\textbf{u}}_{i-4}, \textit{\textbf{u}}_{i-3}, \textit{\textbf{u}}_{i-2})$,
$lk(u_{i+1})=C(u_{i}, \textit{\textbf{u}}_{i-1} \textbf{b}, \textbf{c}, u_{i+2},$  $\textit{\textbf{u}}_{i+3}, \textit{\textbf{u}}_{i+4}, \textit{\textbf{u}}_{i+5}, \textit{\textbf{u}}_{i+6}, \textbf{k}, \textbf{l}, \textbf{m}, \textbf{n}, \textbf{o}, p, \textbf{q})$,\newline
$lk(u_{i+2})=C(u_{i+1}, \textit{\textbf{u}}_{i}, \textit{\textbf{u}}_{i-1} \textbf{b}, {c},  \textbf{d}, u_{i+3},$  $\textit{\textbf{u}}_{i+4}, \textit{\textbf{u}}_{i+5}, \textit{\textbf{u}}_{i+6}, \textbf{k}, \textbf{l}, \textbf{m}, \textbf{n}, \textbf{o}, \textbf{p})$,
$lk(u_{i+3})=C(u_{i+2},$ $\textit{\textbf{c}}, d, \textbf{e}, \textbf{f}, \textbf{g}, u_{i+4} \textit{\textbf{u}}_{i+5} \textit{\textbf{u}}_{i+6}, \textbf{k}, \textbf{l}, \textbf{m},$ $ \textbf{n}, \textbf{o}, \textbf{p}, \textit{\textbf{u}}_{i+1})$,
$lk(u_{i+4})=C(u_{i+3}, \textit{\textbf{d}}, \textbf{e}, \textbf{f}, g, \textbf{h},$  $u_{i+5} \textit{\textbf{u}}_{i+6}, \textbf{k}, \textbf{l},$ $\textbf{m}, \textbf{n}, \textbf{o}, \textbf{p}, \textit{\textbf{u}}_{i+1}, \textit{\textbf{u}}_{i+2})$,
$lk(u_{i+5})=C(u_{i+4}, \textit{\textbf{g}}, h, \textbf{i}, \textit{\textbf{u}}_{i+8}, \textit{\textbf{u}}_{i+7}, u_{i+6}, \textbf{k}, \textbf{l}, \textbf{m}, \textbf{n}, \textbf{o}, \textbf{p}, \textit{\textbf{u}}_{i+1}, \textit{\textbf{u}}_{i+2},$ $\textit{\textbf{u}}_{i+3})$, and
$lk(u_{i+6})=C(u_{i+5}, \textit{\textbf{h}}, i, $ $\textit{\textbf{u}}_{i+8}, u_{i+7}, \textbf{j}, {k}, \textbf{l}, \textbf{m}, \textbf{n}, \textbf{o}, \textbf{p}, \textit{\textbf{u}}_{i+1}, \textit{\textbf{u}}_{i+2}, \textit{\textbf{u}}_{i+3}, \textit{\textbf{u}}_{i+4})$}.

\begin{defn}\label{defn4612-1} Let $P_1:=P(\dots,v_{i-1},v_{i},v_{i+1},\dots)$ be a path in edge graph of $M$. We say $P_1$ of type $H_1$ if $L_1 :=P(v_{t}, v_{t+1}, v_{t+2}, v_{t+3}, v_{t+4}, v_{t+5}, v_{t+6})$ is a subpath of $P_1$ or lies in the extended path of $P_1$ then either $L_1\mapsto Q(i)$ by $v_{j}\mapsto u_{j}$, $L_1 \mapsto Q(i+1)$ by $v_{j}\mapsto u_{j+1}$, $L_1 \mapsto Q(i+2)$ by $v_{j}\mapsto u_{j+2}$, $L_1 \mapsto Q(i+3)$ by $v_{j}\mapsto u_{j+3}$, $L_1 \mapsto Q(i+4)$ by $v_{j}\mapsto u_{j+4}$ or $L_1 \mapsto Q(i+5)$ by $v_{j}\mapsto u_{j+5}$ for $j \in \{t, t+1, t+2, t+3, t+4, t+5, t+6\}$.
\end{defn}

\begin{defn}\label{defn4612-2} Let $P_2 :=P(\dots,x_{i-1},x_{i},x_{i+1},\dots)$ be a path in edge graph of $M$. We say $P_2$ of type $H_1'$ if $L_2 : = P(x_{t}, x_{t+1}, x_{t+2}, x_{t+3}, x_{t+4}, x_{t+5}, x_{t+6})$ is a subpath of $P_2$ or lies in the extended path of $P_2$ then either $L_2 \mapsto Q(i)$ by $x_{j}\mapsto u_{2t+6-j}$, $L_2 \mapsto Q(i+1)$ by $x_{j}\mapsto u_{2t+6-j}$, $L_2 \mapsto Q(i+2)$ by $x_{j}\mapsto u_{2t+6-j}$, $L_1 \mapsto Q(i+3)$ by $v_{j}\mapsto u_{2t+6-j}$, $L_1 \mapsto Q(i+4)$ by $v_{j}\mapsto u_{2t+6-j}$ or $L_1 \mapsto Q(i+5)$ by $v_{j}\mapsto u_{2t+6-j}$ for $j \in \{t, t+1, t+2, t+3, t+4, t+5, t+6\}$.
\end{defn}

Let $P$ be a maximal path of type $H_{1}$ or $H_1'$ in $M$. By the similar argument of Lemma \ref{lem3342-cycle}, the map $M$ contains an edge $e$ which defines a cycle $P \cup \{e\}$ of type $H_{1}$ or $H_1'$. The cycle $C :=P \cup \{e\}$ is a non-contractible cycle (by the similar argument of Lemma \ref{lem3342-non-con}). Observe that the cycles of types $H_1$ and $H_1'$ are mirror image of each other. So, it follows that they define same type of cycles as in Section \ref{3461}. So, we consider only type $H_1$.  Let $C_{1},$ $C_{2},\cdots,C_{m}$ be a sequence of homologous cycles of type $H_{1}$ in $M$. We argue with the similar argument of Lemma \ref{lem3342-all-length} and we get that length($C_{i}$) = length($C_{j}$) $1 \le i, j \le m$. By Definition \ref{defn4612-1}, there are three cycles of type $H_{1}$ through each vertex of $M$. Let $v \in V(M)$ and $L_{1}(v),$ $L_{2}(v),$ $L_{3}(v)$ denote three cycles through $v$. Define a $T(r, s, k)$ representation of $M$ by the similar construction which is given in Section \ref{33421}. In this process, we cut $M$ along $L_1$ and then, take second cut along the cycle $L_{3}$ where the starting adjacent face to base horizontal cycle $L_{1}$ is a $6$-gon. So, every map $M$ has a $T(r, s, k)$ representation. In Lemma \ref{lem4612-hom}, we show that map of type $\{4, 6, 12\}$ contains at most three non-homologous cycles of type $H_{1}$ of different lengths.

\begin{lemma} \label{lem4612-hom} The map $M$ contains at most three non-homologous cycles of type $H_{1}$ of different lengths.
\end{lemma}

\begin{proof} We proceed as in the case of proof of the Lemma \ref{lem346-hom}. In this process, consider map $\{4, 6, 12\}$ in place of $\{3^{4}, 6\}$ and cycle of type $H_1$ in place of $Y_1$. Let $u \in V(M)$ and $T(r, s, k)$ denote a $(r, s, k)$-representation of $M$. Let $L_{1}$, $L_{2}$ and $L_{3}$ denote three cycles of type $H_{1}$ through $u$ in $M$. They are not identical as $L_{i}$ divides link $lk(u)$ into different ratio. Also, cycles are not disjoint as $u \in V(L_{i})\cap V(L_{j})$ for $i \neq j$. Therefore, cycles are not homologous to each other. Again, let $v \in V(M)$, $u \neq v$ and consider cycles of type $H_{1}$ at $v$ in $T(r, s, k)$. Let $L_1',$ $L_2',$ $L_3'$ denote three cycles through $v$ of type $H_1$. Then, by the definition of cycle of type $H_1$ and considering cylinder, $L_i$ and $L_{j}'$ are homologous for some $i, j \in \{1,$ $2,$ $3\}$. This is hold for any vertex $v$ of $M$. Thus, $M$ contains at most three non-homologous cycles of type $H_{1}$. We proceed as in the case of proof of Lemma \ref{lem3342-all-length} to show that the homologous cycles of type $H_1$ have same length. This completes the proof.
\end{proof} 

We define admissible relations among $r$, $s$, $k$ of $T(r, s, k)$ such that representation $T(r, s, k)$ gives a map of type $\{4, 6, 12\}$ after identifying their boundaries.

\begin{lemma}\label{lem4612-all} The maps of type $\{4, 6, 12\}$ of the form $T(r, s, k)$ exist if and only if the following holds : (i) $s \geq  2$ even, (ii) $6 \mid r$, (iii) number of vertices of $T(r, s, k)$ = $rs \geq  36$, (iv) $r \geq 18$ if $s = 2$, (v) $r \geq 12$ if $s \ge 4$, (vi) $k \in \{ 6t+9 \colon 0 \leq t \leq \frac{r-18}{6}\}$ if $s = 2 $ $\&$ $ k \in \{ 6t+3 \colon 0 \leq t \leq \frac{r-6}{6}\}$ if $s \geq 4$.
\end{lemma}

\begin{proof}
 We proceed as in the case of proof of the Lemma \ref{lem346-all}. We proof this lemma by considering link of some vertices in $T(r, s, k)$. We consider map of type $\{4, 6, 12\}$ in place of type $\{3^{4}, 6\}$ and different values of $r$, $s$ and $k$ in the proof of Lemma \ref{lem346-all}. Thus, we get all possible rages of $r$, $s$ and $k$ of $T(r, s, k)$. This completes the proof.
\end{proof}

\begin{table}
\tiny 
\caption{Maps of type $\{4, 6, 12\}$}
\centering 
\begin{tabular}{c c c c} 
\hline\hline 
$n$&Equivalence classes& Length of cycles & $i(n)$\\ [0.3ex] 
\hline 
36 & T(18, 2, 9) & (18, 18, 18) & 1(36)\\
\hline
48 & T(24, 2, 9), T(12, 4, 9) & (24, 12, 24) & 2(48)\\
   & T(24, 2, 15) & &\\

 & T(12, 4, 3) & (12, 12, 12)&\\
\hline
60 & T(30, 2, 9), T(30, 2, 15) & (30, 30, 30) & 1(60)\\
& T(30, 2, 21)& &\\
\hline\hline
\end{tabular}
\label{table6} 
\end{table}

Let $M_{1}$ and $M_{2}$ be two maps of type $\{4, 6, 12\}$ on same number of vertices. Let $T(r_{i}, s_{i}, k_{i})$ denote a $(r_i, s_i, k_i)$-representation of $M_i$. Similarly as in Sections \ref{33421}, \ref{324341}, by Lemma \ref{lem4612-hom}, let $b_{i, j}$ = length($L_{i, j}$) where $L_{i, j}$, $j = 1, 2, 3$ denote non-homologous cycles of type $H_{1}$ in $T(r_{i}, s_{i}, k_{i})$. Then

\begin{lemma}\label{lem4612-iso} The map $M_{1} \cong M_{2}$ if and only if $(b_{1, 1}, b_{1, 2}, b_{1, 3})$ = $(b_{2, t_1}, b_{2, t_2}, b_{2, t_3})$ for $t_1 \neq t_2 \neq t_3 \in \{1, 2, 3\}$.
\end{lemma}

\begin{proof}
We proceed as in the case of the proof of Lemma \ref{lem346-iso}. Let $r = r_{1}, s = s_{1}, k = k_{1}$. Consider horizontal cycles in $T(r_{i},s_{i},k_{i})$ of type $H_{1}$. We proceed with similar argument as in Lemma \ref{lem3342-iso}. So, we get $M_{1} \cong M_{2}$. Again, if $(r, s, k) \neq (r_{1},s_{1},k_{1})$ then we proceed as in the case of proof of Lemmas \ref{lem32434-iso} and \ref{lem3636-iso}. Converse part follows from similar argument of the converse part of Lemma \ref{lem3342-iso}. This completes the proof.
\end{proof} 

As in Section \ref{33421}, by Lemmas \ref{lem4612-all} and \ref{lem4612-iso}, the maps of type $\{4, 6, 12\}$ can be classified on different number of vertices. We have done calculation for the vertices $\le 60$. We have listed the obtained objects in the form of their $T(r, s, k)$ representation in Table \ref{table6}.

\section{Maps of type $\{3, 4, 6, 4\}$}\label{34641}

Let $M$ be a semi-equivelar map of type $\{3, 4, 6, 4\}$ on the torus. We define path $W_1$ in the edge graph of $M$ on the torus. Let \textit{$Q(i) := P(u_{i}, u_{i+1}, u_{i+2}, u_{i+3})$ be a path in $M$, where $lk(u_{i})=C(u_{i-1}, \textbf{a}, {b}, u_{i+1}, \textbf{i}, j, \textbf{k}, \textbf{l}, \textit{\textbf{u}}_{i-2})$, $lk(u_{i+1})=C(u_{i}, b, \textbf{c}, u_{i+2},$ $\textit{\textbf{u}}_{i+3}, \textbf{g}, \textbf{h}, {i}, \textbf{j})$, $lk(u_{i+2})=C(u_{i+1}, \textbf{b}, {c}, d, \textit{\textbf{e}}, u_{i+3}, \textbf{g}, \textbf{h}, \textbf{i})$ and $lk(u_{i+3})=C(u_{i+2}, \textbf{d}, {e}, u_{i+4}, \textit{\textbf{f}},$ $g, \textbf{h}, \textbf{i}, \textit{\textbf{u}}_{i+1})$}

\begin{defn}\label{defn3464-1} Let $P_1 :=P(\dots,v_{i-1},v_{i},v_{i+1},\dots)$ be a path in edge graph of $M$. We say $P_1$ of type $W_1$ if $L_1 := P(v_{t}, v_{t+1}, v_{t+2}, v_{t+3})$ is a subpath of $P_1$ or in a path containing $P_1$. In this case, either $L_1\mapsto Q(i)$ by $v_{j}\mapsto u_{j}$, $L_1 \mapsto Q(i+1)$ by $v_{j}\mapsto u_{j+1}$ or $L_1 \mapsto Q(i+2)$ by $v_{j}\mapsto u_{j+2}$ for $j \in \{t, t+1, t+2, t+3\}$.
\end{defn}

\begin{defn}\label{defn3464-2} Let $P_2 := P(\dots,x_{i-1},x_{i},x_{i+1},\dots)$ be a path in edge graph of $M$. We say $P_2$ of type $W_1'$ if $L_2 : = P(x_{t}, x_{t+1}, x_{t+2}, x_{t+3})$ is a subpath of $P_2$ or in a path containing $P_2$. In this case, either $L_2 \mapsto Q(i)$ by $x_{j}\mapsto u_{2t+3-j}$, $L_2 \mapsto Q(i+1)$ by $x_{j}\mapsto u_{2t+3-j}$ or $L_2 \mapsto Q(i+2)$ by $x_{j}\mapsto u_{2t+3-j}$ for $j \in \{t, t+1, t+2, t+3\}$.
\end{defn}

We consider only cycle of type $W_1$ as $W_1$ and $W_1'$ define same type of cycle (by the similar argument as in Section \ref{3461}). Repeat the similar argument which is given in Section \ref{46121} and define a $T(r, s, k)$ representation. In this process, we consider path of type $W_1$ in place of $H_1$. Now, by Definition \ref{defn3464-1}, there are three cycles through each vertex of $M$. Let $v \in V(M)$ and $L_1$, $L_2$ and $L_3$ be three cycles through $v$. We cut along $L_1$ and then, take second cut along $L_{3}$ where the starting adjacent face to base horizontal cycle $L_{1}$ is a $4$-gon. So, every map has this $T(r, s, k)$ representation. In Lemma \ref{lem3464-hom}, we show that map of type $\{3, 4, 6, 4\}$ contains at most three non-homologous cycles of type $W_{1}$ of different lengths.

\begin{lemma} \label{lem3464-hom} The map $M$ contains at most three non-homologous cycles of type $W_{1}$ of different lengths.
\end{lemma}

\begin{proof}
As above, proceed with similar argument of Lemma \ref{lem346-hom}. Consider map of type $\{3, 4, 6, 4\}$ in place of $\{3^{4}, 6\}$ and cycle of type $W_1$ in place of $Y_1$. Let $w_1 \neq w_2$ be two vertices of $M$. Let $J_1, J_2, J_3$ denote three cycles through $w_1$ and $J_1', J_2' J_3'$ denote three cycles through $w_2$. Then, by the definition of cycle of type $W_1$, we get a cylinder which is bounded by $J_{i}$ and $J_{j}'$. That is, $J_{i}$ and $J_{j}'$ are homologous for some $i, j \in \{1,$ $2,$ $3\}$. This is hold for arbitrary vertex of $M$. We proceed as in the case of proof of Lemma \ref{lem3342-all-length} to show that the homologous cycles of type $W_1$ have same length. Thus, the map $M$ contains at most three non-homologous cycles of type $W_{1}$ of different lengths. This completes the proof.
\end{proof} 

We define admissible relations among $r$, $s$, $k$ of $T(r, s, k)$ such that representation $T(r, s, k)$ gives a map of type $\{3, 4, 6, 4\}$ after identifying their boundaries.

\begin{lemma}\label{lem3464-all} The maps of type $\{3, 4, 6, 4\}$ of the form $T(r, s, k)$ exist if and only if the following holds : (i) $s \geq  2$ even, (ii) $3 \mid r$, (iii) number of vertices of $T(r, s, k)$ = $rs \geq  18$, (iv) $r \geq 9$ if $s = 2$, (v) $r \geq 6$ if $s \ge 4$, (vi) $k \in \{ 3t+4 \colon 0 \leq t \leq \frac{r-9}{3}\}$ if $s = 2 $ $\&$ $ k \in \{ 3t+1 \colon 0 \leq t \leq \frac{r-3}{3}\}$ if $s \geq 4$.
\end{lemma}

\begin{proof}
We follow similar argument of the proof of the Lemma \ref{lem346-all}. We proof this lemma by considering link of some vertices in $T(r, s, k)$ and by showing that link of those vertices are not cycle if we consider the values of $r$, $s$ and $k$ outside the given range in Lemma \ref{lem3464-all}. So, consider map of type $\{3, 4, 6, 4\}$ in place of type $\{3^{4}, 6\}$ and different ranges of $r$, $s$ and $k$ in the proof of Lemma \ref{lem346-all}. Thus, we get all the cases of this lemma. This completes the proof.
\end{proof} 

\begin{table}
\tiny 
\caption{Maps of type $\{3, 4, 6, 4\}$}
\centering 
\begin{tabular}{c c c c} 
\hline\hline 
$n$&Equivalence classes& Length of cycles & $i(n)$\\ [0.3ex] 
\hline 
18 & T(9, 2, 4) & (9, 9, 9) & 1(18)\\
\hline
24 & T(12, 2, 4), T(12, 2, 7), T(6, 4, 4) & (12, 6, 12) & 2(24) \\ 

& T(6, 4, 1)& (6, 6, 6)& \\
\hline
30 & T(15, 2, 4), T(15, 2, 7), T(15, 2, 10) & (15, 15, 15) & 1(30)\\
\hline
36 & T(18, 2, 4), T(18, 2, 13), T(9, 4, 7) & (18, 9, 18) & 2(36)\\

& T(18, 2, 7), T(18, 2, 10), T(9 , 4, 1) & (18, 9, 6)&\\
&T(9, 4, 4), T(6, 6, 4), T(6, 6, 1)& &\\
\hline
42 & T(21, 2, 4), T(21, 2, 7), T(21, 2, 10) & (21, 21, 21) & 1(42)\\
&T(21, 2, 13), T(21, 2, 16)& &\\
\hline
48 & T(24, 2, 4), T(24, 2, 7), T(24, 2, 16) & (24, 24, 12) & 3(48)\\
& T(24, 2, 19), T(12, 4, 4), T(12, 4, 10)& &\\

& T(24, 2, 10), T(24, 2, 13), T(6, 8, 4) & (24, 24, 6)&\\

& T(12, 4, 1), T(12, 4, 7), T(6, 8, 1) & (12, 12, 6) &\\
\hline
54 & T(27, 2, 4), T(27, 2, 13), T(27, 2, 22) & (27, 27, 27) & 3(54)\\

 & T(27, 2, 7), T(27, 2, 10), T(27, 2, 16) & (27, 27, 9)&\\
& T(27, 2, 19), T(9, 6, 4), T(9, 6, 7)& &\\

 & T(9, 6, 1) & (9, 9, 9)&\\
\hline\hline
\end{tabular}
\label{table7} 
\end{table}

Let $M_{i}$ for $i = 1, 2$ be maps of type $\{3, 4, 6, 4\}$ on same number of vertices and $T_i = T(r_{i}, s_{i}, k_{i})$ be $(r_i, s_i, k_i)$-representations of $M_i$. By Lemma \ref{lem3464-hom}, there are at most three non-homologous cycles of different lengths in $T_i$. So, let $c_{i, j}$ = length($N_{i, j}$) where $N_{i, j}$, $j = 1, 2, 3$ denote non-homologous cycles of type $W_{1}$ in $T_i$. Then,

\begin{lemma}\label{lem3464-iso} The map $M_{1} \cong M_{2}$ if and only if $(c_{1, 1}, c_{1, 2}, c_{1, 3})$ = $(c_{2, t_1}, c_{2, t_2}, c_{2, t_3})$ for $t_1 \neq t_2 \neq t_3 \in \{1, 2, 3\}$.
\end{lemma}

\begin{proof}
We proceed as in the case of proof of Lemma \ref{lem346-iso}. Let $r = r_{1}, s = s_{1}, k = k_{1}$. We consider horizontal cycles of $T(r_{i},s_{i},k_{i})$ of type $W_{1}$. Proceed with similar argument as in the Lemma \ref{lem3342-iso}. Thus, we get $M_{1} \cong M_{2}$. Again, if $(r, s, k) \neq (r_{1},s_{1},k_{1})$ then we proceed as in the case of the proof of Lemma \ref{lem32434-iso} and Lemma \ref{lem3636-iso}. Converse part follows from similar argument of the converse part of Lemma \ref{lem3342-iso}. This completes the proof.
\end{proof} 

As in Section \ref{33421}, by Lemmas \ref{lem3464-all} and \ref{lem3464-iso}, the maps of type $\{3, 4, 6, 4\}$ can be classified on different number of vertices. We have done calculation for the vertices $\le 54$. We have listed the obtained objects in the form of their $T(r, s, k)$ representation in Table \ref{table7}.

\section{Maps of type $\{4, 8^{2}\}$}\label{4881}

Let $M$ be a semi-equivelar map of type $\{ 4, 8^2\}$ on the torus. We define a fixed type of path $Z_1$ in $M$. Let $Q(i) := P(u_{i}, u_{i+1}, u_{i+2}, u_{i+3}, u_{i+4})$ be a path in $M$, where \textit{$lk(u_{i}) = C(u_{i-1}, \textbf{f}, \textbf{g}, \textbf{h}, \textbf{i}, \textit{\textbf{u}}_{i+2}, u_{i+1}, \textit{\textbf{s}}, a, \textit{\textbf{b}}, \textit{\textbf{c}}, \textit{\textbf{d}}, \textit{\textbf{u}}_{i-3}, \textit{\textbf{u}}_{i-2})$, $lk(u_{i+1}) = C(u_i, \textit{\textbf{u}}_{i-1}, \textit{\textbf{f}}, \textbf{g}, \textbf{h}, \textbf{i}, u_{i+2}, \newline \textit{\textbf{u}}_{i+3}, \textit{\textbf{u}}_{i+4}, \textbf{p},\textit{\textbf{q}}, \textbf{r}, s, \textit{\textbf{a}})$, $lk(u_{i+2})=C(u_{i+1}, \textit{\textbf{u}}_{i}, \textit{\textbf{u}}_{i-1}, \textbf{f}, \textbf{g}, \textbf{h}, i, \textbf{j}, u_{i+3}, \textit{\textbf{u}}_{i+4}, \textbf{p}, \textit{\textbf{q}},  \textbf{r}, \textbf{s})$, \newline  $lk(u_{i+3})=C(u_{i+2}, \textbf{i}, j, \textit{\textbf{k}},   \textbf{l}, \textbf{m}, \textit{\textbf{u}}_{i+6}, \textit{\textbf{u}}_{i+5}, u_{i+4}, \textbf{p},   \textit{\textbf{q}}, \textbf{r}, \textbf{s}, \textit{\textbf{u}}_{i+1})$,
$lk(u_{i+4})=C(u_{i+3}, \textbf{j}, \textit{\textbf{k}},$ $\textbf{l}, \textbf{m},$ $ \textit{\textbf{u}}_{i+6}, u_{i+5}, \textbf{o}, p, \textit{\textbf{q}}, \textbf{r}, \textbf{s}, \textit{\textbf{u}}_{i+1}, \textit{\textbf{u}}_{i+2})$}.

\begin{defn}\label{defn482-1} Let $R_1 := P(\dots,v_{i-1},v_i,v_{i+1},\dots)$ be a path in edge graph of $M$. We say $R_1$ of type $Z_1$ if $L_1 :=P(v_{t}, v_{t+1}, v_{t+2}, v_{t+3}, v_{t+4})$ is a subpath of $R_1$ or in the extended path of $R_1$ then either $L_1\mapsto Q(i)$ by $v_{j}\mapsto u_{j}$, $L_1 \mapsto Q(i+1)$ by $v_{j}\mapsto u_{j+1}$, $L_1 \mapsto Q(i+2)$ by $v_{j}\mapsto u_{j+2}$ or $L_1 \mapsto Q(i+3)$ by $v_{j}\mapsto u_{j+3}$ for $j \in \{t, t+1, t+2, t+3, t+4\}$.
\end{defn}

\begin{defn}\label{defn482-2} Let $R_2 := P(\dots,x_{i-1},x_i,x_{i+1},\dots)$ be a path in edge graph of $M$. We say $R_2$ of type $Z_1'$ if $L_2 : = P(x_{t}, x_{t+1}, x_{t+2}, x_{t+3}, x_{t+4})$ is a subpath of $R_2$ or in the extended path of $R_2$ then either $L_2 \mapsto Q(i)$ by $x_{j}\mapsto u_{2t+4-j}$, $L_2 \mapsto Q(i+1)$ by $x_{j}\mapsto u_{2t+4-j}$, $L_2 \mapsto Q(i+2)$ by $x_{j}\mapsto u_{2t+4-j}$ or $L_1 \mapsto Q(i+3)$ by $v_{j}\mapsto u_{2t+4-j}$ for $j \in \{t, t+1, t+2, t+3, t+4\}$.
\end{defn}

As in Section \ref{3461}, we consider path of type $Z_1$ as $Z_1$ and $Z_1'$ define same type of path (by the similar argument as in Section \ref{3461}). Let $P$ be a maximal path of type $Z_{1}$. We argue similar argument of Lemma \ref{lem3342-cycle} for $P$. So, we get an edge $e$ in $M$ such that $P \cup e$ is a cycle of type $Z_1$. Therefore, every maximal path of type $Z_1$ is a cycle. Let $C = P \cup e$. The cycle $C$ is of type $Z_1$ and non-contractible (by the similar argument of Lemma \ref{lem3342-non-con}). Let $C_1,$  $C_{2}, \cdots, C_{t}$ be a sequence of homologous cycles of type $Z_1$. Then, we proceed with the similar argument of Lemma \ref{lem3342-all-length}. Thus, length($C_{i}$) = length($C_{j})$ for $1 \leq i, j \leq t$. By Definition \ref{defn482-1}, there are two cycles of type $Z_{1}$ through each vertex of $M$. Let $v \in V(M)$ and $L_{1}(v),$ $L_{2}(v)$ be two cycles through $v$. We repeat the similar construction of $(r, s, k)$-representation of Section \ref{33421} for $M$. So, we get a $T(r, s, k)$ representation of $M$. In this processes, we take first cut along $L_1$ and then, second cut along $L_{2}$ where the starting adjacent face to the base horizontal cycle $L_{1}$ is a $4$-gon. By this construction, every map of type $\{4, 8^2\}$ on the torus has a $T(r, s, k)$ representation.

We have two cycles of type $Z_1$ through each vertex of $M$. Therefore, by Lemma \ref{lem482-hom}, map $M$ contains at most two non-homologous cycles of type $Z_{1}$ of different lengths.

\begin{lemma} \label{lem482-hom} The map $M$ contains at most two non-homologous cycles of type $Z_{1}$ of different lengths.
\end{lemma}

\begin{proof}
Let $v$ be a vertex in $M$. By the definition of cycle of type $Z_{1}$, we have two cycles namely, $C_{1}$ and $C_{2}$ through $v$. We proceed as in the case of proof of Lemma \ref{lem346-hom}. So, we get that the map $M$ contains at most two non-homologous cycles of type $Z_{1}$ of different lengths.
\end{proof} 

We claim in the Lemma \ref{lem482-1} that map of type $\{4, 8^{2}\}$ does not contain cycle of type $Z_{1}$ which has length four. We use this result to classify the maps of type $\{4, 8^{2}\}$ on the torus.

\begin{lemma} \label{lem482-1} The representation $T(r, s, k)$ does not contain cycle of type $Z_{1}$ of length four.
\end{lemma}

\begin{proof} Suppose, $T(r, s, k)$ has a cycle, say $C$ of length four which is of type $Z_{1}$. Let $C_{F_{8}}$ denote a boundary cycle of a $8$-gon $F_8$. By the definition of cycle of type $Z_{1}$, if $C\cap C_{F_{8}} \neq \emptyset$ then $C\cap C_{F_{8}}$ is a path of length three. So, let $C \cap C_{F_{8}} = P(u_{1}, u_{2}, u_{3}, u_{4})$. If $C$ is a cycle of length four then, $u_{1} = u_{4}$ which is a contradiction as $u_{1},$ $u_{4}$ are in $C_{F_{8}}$ and $C_{F_{8}}$ is a cycle without chord. Again, by the definition of cycle of type $Z_{1}$, $C$ must intersect a $8$-gon. Thus, length$(C) > 4$. So, a map $M$ of type $\{4, 8^{2}\}$ does not contain a cycle $C$ of type $Z_{1}$ of length four. This completes the proof.
\end{proof}

We define admissible relations among $r$, $s$, $k$ of $T(r, s, k)$.

\begin{lemma}\label{lem482-all}  The maps of type $\{4, 8^{2}\}$ of the form $T(r, s, k)$ exist if and only if the following holds : (i) $4 \mid r$ and $s \geq  1$, (ii) number of vertices of $T(r, s, k)$ = $rs \geq  20$, (iii) $r \geq 20$ if $s = 1$, (iv) $r \geq 16$ if $s = 2$, (v) $r \geq 8$ if $s \ge 3$, (vi) $k \in \{ 4t+6 \colon 0 \leq t \leq \frac{r-12}{4}\}$ if $s = 1$, $k \in \{ 4t+7 \colon 0 \leq t \leq \frac{r-16}{4}\}$ if $s = 2 $ $\&$ $ k \in \{ 4t-1 (mod~r) \colon 0 \leq t \leq \frac{r-4}{4}\}$ if $s \geq 3$. 
\end{lemma}

\begin{proof}
We proceed as in the case of proof of Lemma \ref{lem346-all} and use Lemma \ref{lem482-1} to prove this lemma. Considering link of some vertices in $T(r, s, k)$, map of type $\{4, 8^{2}\}$ in place of type $\{3^{4}, 6\}$, and different values of $r$, $s$ and $k$ in the proof of Lemma \ref{lem346-all}. So, we get all possible rages of $r$, $s$ and $k$ of $T(r, s, k)$. This completes the proof.
\end{proof} 

\begin{table}
\tiny 
\caption{Maps of type $\{4, 8^{2}\}$}
\centering 
\begin{tabular}{c c c c} 
\hline\hline 
$n$&Equivalence classes& Length of cycles & $i(n)$\\ [0.3ex] 
\hline 
20 & T(20, 1, 6), T( 20, 1, 14) & (20, 20) & 1(20) \\
\hline
24 & T(24, 1, 6), T(24, 1, 18) & (24, 8) & 2(24) \\
&T(8, 3, 2), T(8, 3, 6)& &\\
 & T(24, 1, 14), T(24 , 1, 10) & (24, 24)&\\
[1ex] 
\hline\hline 
\end{tabular}
\label{table8} 
\end{table}

Let $M_{i}$ for $i = 1, 2$ be maps of type $\{4, 8^{2}\}$ on same number of vertices and $T_i = T(r_{i}, s_{i}, k_{i})$ represent of $M_i$. By Lemma \ref{lem482-hom}, there are at most two non-homologous cycles of different lengths in $T_i$. So, let $a_{i, j}$ = length($C_{i, j}$) where $C_{i, j}$, $j = 1, 2$ denote non-homologous cycles of type $Z_{1}$ in $T_i$. Then, 

\begin{lemma}\label{lem482-iso} The map $M_{1} \cong M_{2}$ if and only if $(a_{1, 1}, a_{1, 2})$ = $(a_{2, t_1}, a_{2, t_2})$ for $t_1 \neq t_2 \in \{1, 2\}$.
\end{lemma}

\begin{proof}
We proceed as in the case of proof of Lemma \ref{lem346-iso}. Let $r = r_{1}, s = s_{1}, k = k_{1}$. Consider horizontal cycles in $T(r_{i},s_{i},k_{i})$ of type $Z_{1}$. Proceed with similar argument as in Lemma \ref{lem3342-iso}. So, we get $M_{1} \cong M_{2}$. Again, if $(r, s, k) \neq (r_{1},s_{1},k_{1})$ then we proceed as in the case of the proof of Lemmas \ref{lem32434-iso} and \ref{lem3636-iso}. Converse part follows from the similar argument of the converse part of Lemma \ref{lem3342-iso}.
\end{proof} 

As in Section \ref{33421}, by Lemmas \ref{lem482-all}, \ref{lem482-iso}, the maps of type $\{4, 8^{2}\}$ can be classified up to isomorphism on different number of vertices. We have done calculation for the vertices $\le 24$. We have listed the obtained objects in the form of their $T(r, s, k)$ representation in Table \ref{table8}.

\section{Semi-equivelar maps}\label{result4}

\begin{proof}[Proof of Theorem \ref{thm-main}] The proof of the Theorem \ref{thm-main} follows from the Sections \ref{33421}, \ref{324341}, \ref{36361}, \ref{31221}, \ref{3461}, \ref{46121}, \ref{34641}, \ref{4881}. Let $M$ be a map on $n$ vertices of type $\{3^{3}, 4^{2}\}$ on the torus. We consider all admissible $T(r, s, k)$ representations of $M$ by Lemma \ref{lem3342-all-rep}. We calculate length of the cycles of types $A_1$, $A_2$, $A_3$ and $A_4$. We classify them by Lemma \ref{lem3342-iso}. In Table \ref{table1}, we have classified upto $22$ vertices. Similarly, we consider maps of types $\{3^{2}, 4, 3, 4\}$, $\{3, 6, 3, 6\}$, $\{3, 12^{2}\}$, $\{3^{4}, 6\}$, $\{4, 6, 12\}$, $\{3,4,6,4\}$, $\{4, 8^{2}\}$ on the torus. That is, we consider cycles of type $B_1$ in maps of type $\{3^{2}, 4, 3, 4\}$ and classify by Lemma \ref{lem32434-iso}. Table \ref{table2} contains the classified maps upto $32$ vertices. Consider cycles of type $X_1$ in maps of type $\{3, 6, 3, 6\}$ and classify by Lemma \ref{lem3636-iso}. Table \ref{table3} contains the maps upto $30$ vertices. Consider cycles of type $G_1$ in maps of type $\{3, 12^{2}\}$ and classify by Lemma \ref{lem3122-iso}. Table \ref{table4} contains the maps upto $48$ vertices. Consider cycles of type $Y_1$ in maps of type $\{3^{4}, 6\}$ and classify by Lemma \ref{lem346-iso}. Table \ref{table5} contains the maps upto $42$ vertices. Consider cycles of type $H_1$ in maps of type $\{4, 6, 12\}$ and classify by Lemma \ref{lem4612-iso}. Table \ref{table6} contains the maps upto $60$ vertices. Consider cycles of type $W_1$ in maps of type $\{3,4,6,4\}$ and classify by Lemma \ref{lem3464-iso}. Table \ref{table7} contains the classified maps upto $54$ vertices. Consider cycles of type $Z_1$ in maps of type $\{4, 8^{2}\}$ and classify by Lemma \ref{lem482-iso}. Table \ref{table8} contains the maps upto $24$ vertices. This completes the proof.
\end{proof} 

\section{Acknowledgement}

The authors are grateful to the anonymous referee whose comments led to a substantial
improvement in the paper. Work of second author is partially supported by SERB, DST grant No. SR/S4/MS:717/10.

\end{document}